\documentclass{svjour3}   

\usepackage{microtype}
\usepackage{algorithmic}
\usepackage[ruled,vlined]{algorithm2e}
\usepackage{graphicx,array}
\usepackage{subfigure,tabularx}
\usepackage{booktabs} 
\usepackage{amssymb,amsmath,amsfonts,mathrsfs}
\usepackage{etoolbox,verbatim,color}
\usepackage{url}
\usepackage{hyperref}
\usepackage{xcolor}    
\usepackage{nicefrac} 

\newcommand{\mL}{\mathcal L}

\newtheorem{assumption}{Assumption}

\usepackage{multicol}
\usepackage{multirow}
\catcode`~=11 \def\UrlSpecials{\do\~{\kern -.15em\lower .7ex\hbox{~}\kern .04em}} \catcode`~=13 

\newcommand{\lrpar}[1]{\left(#1\right)}

\newcommand{\bbE}{\mathbb{E}}

\newcommand{\bbR}{\mathbb{R}}

\DeclareMathAlphabet{\mathbsf}{OT1}{cmss}{bx}{n}
\DeclareMathAlphabet{\mathssf}{OT1}{cmss}{m}{sl}

\DeclareSymbolFont{bsfletters}{OT1}{cmss}{bx}{n}  
\DeclareSymbolFont{ssfletters}{OT1}{cmss}{m}{n}
\DeclareMathSymbol{\bsfGamma}{0}{bsfletters}{'000}
\DeclareMathSymbol{\ssfGamma}{0}{ssfletters}{'000}
\DeclareMathSymbol{\bsfDelta}{0}{bsfletters}{'001}
\DeclareMathSymbol{\ssfDelta}{0}{ssfletters}{'001}
\DeclareMathSymbol{\bsfTheta}{0}{bsfletters}{'002}
\DeclareMathSymbol{\ssfTheta}{0}{ssfletters}{'002}
\DeclareMathSymbol{\bsfLambda}{0}{bsfletters}{'003}
\DeclareMathSymbol{\ssfLambda}{0}{ssfletters}{'003}
\DeclareMathSymbol{\bsfXi}{0}{bsfletters}{'004}
\DeclareMathSymbol{\ssfXi}{0}{ssfletters}{'004}
\DeclareMathSymbol{\bsfPi}{0}{bsfletters}{'005}
\DeclareMathSymbol{\ssfPi}{0}{ssfletters}{'005}
\DeclareMathSymbol{\bsfSigma}{0}{bsfletters}{'006}
\DeclareMathSymbol{\ssfSigma}{0}{ssfletters}{'006}
\DeclareMathSymbol{\bsfUpsilon}{0}{bsfletters}{'007}
\DeclareMathSymbol{\ssfUpsilon}{0}{ssfletters}{'007}
\DeclareMathSymbol{\bsfPhi}{0}{bsfletters}{'010}
\DeclareMathSymbol{\ssfPhi}{0}{ssfletters}{'010}
\DeclareMathSymbol{\bsfPsi}{0}{bsfletters}{'011}
\DeclareMathSymbol{\ssfPsi}{0}{ssfletters}{'011}
\DeclareMathSymbol{\bsfOmega}{0}{bsfletters}{'012}
\DeclareMathSymbol{\ssfOmega}{0}{ssfletters}{'012}

\DeclareMathOperator*{\argmin}{argmin} 
\DeclareMathOperator*{\dist}{dist}

\definecolor{brightpink}{rgb}{1.0, 0.0, 0.5}

\definecolor{brightpink}{rgb}{1.0, 0.0, 0.5}

\newcommand{\hien}[1]{{{\color{black} #1}}}

\newcommand{\revise}[1]{{{\color{black} #1}}}

\begin{document}


\title{Inertial Alternating Direction Method of Multipliers for \\ 
Non-Convex Non-Smooth Optimization\thanks{L. T. K. Hien and D. N. Phan contributed equally to this work.\\
L. T. K. Hien finished this work when she was at the University of Mons, Belgium.}} 
\titlerunning{iADMM for nonconvex nonsmooth optimization}    
\author{Le Thi Khanh Hien         \and 
        Duy Nhat Phan \and
        Nicolas Gillis  
}

\authorrunning{L. T. K. Hien, D. N. Phan, N. Gillis} 

\institute{L. T. K. Hien \at Huawei Belgium Research Center, 3001 Leuven, Belgium\\
              \email{let.hien@huawei.com}  \\ 
              D. N. Phan \at Dynamic Decision Making Laboratory, Carnegie Mellon University, USA \\
  \email{dnphan@andrew.cmu.edu}               \\
               N. Gillis \at Department of Mathematics and Operational Research,  
Facult\'e Polytechnique, Universit\'e de Mons,
Rue de Houdain 9, 7000 Mons, Belgium\\
              \email{nicolas.gillis@umons.ac.be}       
}


\maketitle

\begin{abstract}
In this paper, we propose an algorithmic framework, dubbed inertial alternating direction methods of multipliers (iADMM), for solving a class of nonconvex nonsmooth multiblock composite optimization problems with linear constraints. Our framework employs the general minimization-majorization (MM) principle to update each block of variables so as to not only unify the convergence analysis of previous ADMM that use specific surrogate functions in the MM step, but also lead to new efficient ADMM schemes. To the best of our knowledge, in the \emph{nonconvex nonsmooth} setting, ADMM used in combination with the MM principle to update each block of variables, and ADMM combined with \emph{inertial terms for the primal variables} have not been studied in the literature.  Under standard  assumptions, we prove the subsequential convergence and global convergence for the generated sequence of iterates. We illustrate the effectiveness of iADMM on a class of nonconvex low-rank representation problems. 
\end{abstract}



\section{Introduction}
\label{sec:Intro}
In this paper, we consider the following nonconvex minimization problem with linear constraints  
\begin{equation}
\label{model}
 \min_{x,y}  \; F(x_1,\ldots,x_s) + h(y) 
 \quad 
\text{ such that}  \quad \sum_{i=1}^s \mathcal A_i x_i + \mathcal  By = b,
\end{equation}
where $y \in \mathbb R^{\mathbf q}$, $x_i \in \mathbb R^{\mathbf n_i}$, $x:=[x_1;\ldots; x_s] \in \mathbb R^{\mathbf n}$, $\mathbf n=\sum_{i=1}^s \mathbf n_i$, $\mathcal A_i$ is a linear map from $\mathbb R^{\mathbf n_i}$ to $\mathbb R^{\mathbf m}$, $\mathcal  B$ is a linear map from $\mathbb R^{\mathbf q}$ to $\mathbb R^{\mathbf m}$, $b\in \mathbb R^{\mathbf m} $, $h:\mathbb R^{\mathbf q} \to \mathbb R$ is a  differentiable function, and $F(x)=f(x)+\sum_{i=1}^s g_i(x_i)$, where $f: \mathbb R^{\mathbf n} \to \mathbb R$ is a nonconvex nonsmooth function and $g_i:\mathbb R^{\mathbf n_i} \to \mathbb R\cup \{+\infty\}$ are  proper lower semi-continuous functions for $i=1,2,\dots,s$. We assume that $F$ satisfies 
$\partial F(x) = \partial_{x_1} F(x) \times \ldots \times \partial_{x_s} F(x)$, where $\partial F$ denote the limiting  subdifferential of $F$ (see the definition in Appendix~\ref{sec:pre}). Note that this condition is satisfied when $f$ is a sum of a continuously differentiable function and a block separable function;  see~\cite[Proposition 2.1]{Attouch2010}. 

\paragraph{Notation.} We denote $[s]:=\{1,\ldots,s\}$. For the $\mathbf p$-dimensional Euclidean space $\mathbb R^{\mathbf p}$, we use $\langle\cdot,\cdot\rangle$ to denote the inner product, and $\|\cdot\|$ to denote the corresponding  induced norm. 
For a linear map $\mathcal  M$, $\mathcal  M^*$  denotes the adjoint linear map with
respect to the inner product, and $\|\mathcal  M\|$ is the induced operator norm of $\mathcal  M$.  We use $\mathcal I$ to denote the identity map. For a positive definite self-adjoint operator $\mathcal  Q$, we denote $\|x\|^2_{\mathcal Q}:=\langle x, \mathcal  Qx\rangle $.  We denote the smallest eigenvalue of a symmetric linear self-map (that is, $\mathcal  M=\mathcal  M^*$) by $\lambda_{\min}(\mathcal  M)$. We use $Im(\mathcal  B)$ to denote the image of $\mathcal  B$.

\subsection{Nonconvex low-rank representation problem}
\label{sec:intro-LRR}
\vspace{-0.05in}

\revise{Low-rank matrix approximations play a central role in various fields of computer science and applied mathematics, and are used in many applications, e.g., 
recommender systems~\cite{koren2009matrix}, 
topic modeling~\cite{lee1999learning}, 
system identification~\cite{markovsky2012low}, 
graph clustering~\cite{von2007tutorial}, 
compression and denoising~\cite{UHZB14}, to cite a few; 
see also below for other examples.  
Given a data matrix, $D$, the goal of low-rank matrix approximations is to find a nearby low-rank matrix, $X$. The low-rank assumption is valid in many applications as there are typically redundancy and correlations within large data sets; see, e.g.,
~\cite{recht2010guaranteed,udell2019big} and the references therein.} 

\revise{In this paper, we will illustrate the use of~\eqref{model} on the} following generalized nonconvex low-rank representation problem:  given a data matrix $D \in \mathbb{R}^{d \times n}$, solve   
\vspace{-0.05in}
\begin{equation}\label{LRR}
\min_{X,Y,Z}  
\,\, \sum_{i=1}^{\min(m,n)} r_1(\sigma(X)) + r_2(Y) + r_3(Z) 
\quad \text{such that}   
\quad D = A_1 X + Y A_2 +  Z,
\end{equation}
where $X\in\mathbb{R}^{m\times n}, Y\in\mathbb{R}^{d\times q}, Z\in \mathbb{R}^{d\times n}$,  
$A_1\in\mathbb{R}^{d\times m}$,  
$A_2\in\mathbb{R}^{q\times n}$,  
\revise{$\sigma(X)$ is the vector of singular values of $X$}, 
$r_1$ is an increasing concave function to promote $X$ 
to be of low rank \revise{(by promoting the sparsity of $\sigma(X)$)}, 
 $r_2$ is a regularization function, 
 and $r_3$ is a function that models some noise; 
 e.g., taking $r_3(Z) = \frac{1}{2}\|Z\|_F^2$ when $Z$ models Gaussian noise.   
\revise{Problem~\ref{LRR} generalize low-rank matrix approximations, taking $A_1$ as the identity matrix and $A_2 = 0$, so that $D = X+Z$ where $X$ is low rank, and $Z$ models the noise.

 In particular, Problem~\eqref{LRR} generalizes the following machine learning problems: } 
\begin{itemize}
\item[(i)] Let $r_1(t) = t^\chi$ with $0< \chi \leq 1$, $r_2(Y) = \sum_{i=1}^{q-1}\|Y_{i} - Y_{i+1}\|$ where $Y_{i}$ is the 
$i$-th column of $Y$, and let $A _1$ and $A_2$ be the identity matrices so that Problem~\eqref{LRR} decomposes the data matrix $D$ 
into the sum of three components, $X$, $Y$ and~$Z$. 
An application is video surveillance where 
each column of $D$ is a vectorized image of a video frame, 
$X$ is a low-rank matrix that plays the role of the background, $Y$ is the foreground that has small variations between its columns (such as slowly moving objectives), and $Z$ represents some 
noise~\cite{Wang2019}. 

\item[(ii)]   When $A _1$ and $A_2$ are identity matrices, $r_1(t) = t$, and  $r_2(Y) = \lambda\|Y\|_1$ for some constant $\lambda > 0$, 
Problem~\eqref{LRR} recovers the robust principal component analysis (robust PCA) model, see, e.g., \cite{Cande2011}. 
Robust PCA decomposes the input matrix $D$ as the sum of a low-rank matrix $X$, a sparse matrix $Y$ modeling gross corruptions and outliers, and an additional noise matrix  $Z$ (e.g., $r_3(Z)$ is a multiple of $\|Z\|_F^2$  to model Gaussian noise). 
Robust PCA is also used for foreground-background separation in video surveillance. 

\item[(iii)]  When  $r_1(t) = t$ and $r_2(Y) = \|Y\|_*$,  Problem~\eqref{LRR} is the latent low-rank representation problem \cite{Liu2011}. 
In~\cite{Liu2011}, authors used 
$A_1 = DP_1$ and $A_2 = P_2^*D$, where $P_1$ and $P_2$ are computed by orthogonalizing the columns of $D^*$ and $D$, respectively. {We will use this application to illustrate the effectiveness of our proposed framework, iADMM, in Section~\ref{sec:numerical}.} 
\end{itemize}
Other applications of Problem~\eqref{model} include statistical learning, see, e.g., \cite{Bach2011,WangTSP2011}, 
 and minimization on compact manifolds, 
 see, e.g., \cite{Lai2014,Wen2010}. 

\subsection{Motivation and related works} 
Let  $\mathcal  A:=[\mathcal  A_1 \ldots \mathcal A_s]$ and $\mathcal A x:=\sum_{i=1}^s \mathcal A_i x_i\in \mathbb R^{\mathbf m}$.  The augmented Lagrangian for Problem~\eqref{model} is 
\begin{equation}
\label{Lagrangian}
\mathcal L(x, y, \omega):= F(x) + h(y)
 +  \langle \omega,\mathcal A x + \mathcal By-b \rangle + \frac{\beta}{2} \|\mathcal A x + \mathcal By-b\|^2,
\end{equation}
where $\beta>0$ is a penalty parameter.  ADMM was first introduced by \cite{Glowinski1975} and \cite{Gabay1976}. It has recently become popular because of its efficacy in solving  emerging large-scale problems in machine learning and computer vision \cite{Boyd2011,Scheinberg2010,Yang2009,Yang2017,Yin2008}. 
For simplicity, let us describe the iteration scheme of a classical ADMM for solving Problem~\eqref{model} with 2 blocks $x$ and $y$: 
\begin{subequations}
\label{ADMM-classic}
\begin{align}
x^{k+1} &\in\argmin_{x} \mL(x,y^{k},\omega^k),\label{subproblem_x}\\
y^{k+1} &\in\argmin_{y} \mL(x^{k+1},y,\omega^k),\label{subproblem_y}\\
\omega^{k+1} &=\omega^k+ \beta (\mathcal A x^{k+1} +\mathcal  B y^{k+1} -b). \label{subproblem_w} 
\end{align}
\end{subequations}
For a multi-block problem, with $s>1$, the scheme is similar, 
see, e.g., \cite{Wang2019}. 
The update of $x$ in~\eqref{subproblem_x} (a similar discussion is applicable to~\eqref{subproblem_y}) can be rewritten as 
$
x^{k+1} \; \in \; \argmin_{x} F(x) + \varphi^k(x),
$
where 
\begin{equation}
\label{eq:varphi}
\varphi^k(x)=\frac{\beta}{2}\|\mathcal A x + \mathcal By^k -b \|^2 + \langle \omega^k,\mathcal A x +  \mathcal By^k -b \rangle. 
\end{equation}
Solving the subproblem~\eqref{subproblem_x} is usually very expensive especially when $F$ is not smooth. 
A remedy is minimizing a suitable surrogate function of $\mL(\cdot, y^k,\omega^k)$ that allows a more efficient update for $x$. For example, since $\varphi^k(x)$ is upper bounded by 
\begin{equation}
\label{eq:varphihat}
\hat \varphi(x)= \varphi^k(x^k)+ \langle \nabla \varphi^k(x^k),x-x^k\rangle + \frac{\kappa\beta}{2} \|x-x^k\|^2,
\end{equation}
 where $\kappa\geq \|\mathcal A^* \mathcal A \|$ (because $\nabla \varphi^k(x)$ is $\beta \|\mathcal A^* \mathcal A \|$-Lipschitz continuous), $x$ can be updated by
$
x^{k+1} \; \in \; \argmin_{x} F(x) + \hat{\varphi}(x),
$
which leads to the linearized ADMM method, see   \cite{Lin2011,XuWu2011}. This update 
has a closed form for some nonsmooth $F$; see \cite{Boyd2014}. 
When $F= f + g$ 	and $f$ is $L_f$-smooth then  we can also use the upper bound $ \hat F(x)=f(x^k)+ \langle \nabla f(x^k),x-x^k\rangle + \frac{L_f}{2} \|x-x^k\|^2 + g(x)$ of $ F$ to obtain 
$
x^{k+1} \; \in \; \argmin_{x} \hat F(x) + \hat{\varphi}(x).
$
This leads to the proximal linearized ADMM method, see  \cite{Bot2020,LiuShen2011}. 
We note that $\mL(\cdot, y^k,\omega^k)$ is always upper bounded by  $\mL(\cdot, y^k,\omega^k) +\mathbf D_\phi(x,x^k)$, where $\mathbf D_\phi$ is the Bregman distance associated with a continuously differentiable convex function $\phi$ on $\mathbb R^n$: 
\begin{equation}
\label{def:Bregman}
\mathbf D_\phi(a,b):=\phi(a)-\phi(b)-\langle \nabla \phi(b),a-b\rangle, \forall a,b\in\mathbb R^n.
\end{equation}
For example, if $\phi(x)=\|x\|_{\mathcal Q}^2=\langle x, \mathcal Q x\rangle$ then $\mathbf D_\phi(a,b)=\|a - b\|^2_{\mathcal Q}$. 
This upper bound leads to proximal ADMM, 
see \cite{DengYin2012,Li2015}. 
The above mentioned upper bound functions are specific examples of surrogate functions for $\mL(\cdot, y^k,\omega^k)$ (see Definition~\ref{def:surrogate} \revise{at}  page~\pageref{def:surrogate} \revise{for the definition of a surrogate function}) while each method of updating $x$ corresponds to a majorization-minimization (MM) step \revise{that minimizes the corresponding majorizer/surrogate function (see \cite{Sun2017} for more specific examples of the MM procedure)}. In the convex setting (that is, $f(\cdot,\cdot)$ is convex), \cite{Lu2016} and \cite{HongADMM2020} use the MM principle to unify and generalize the convergence analysis of many ADMM for multi-blocks problems (that is, $s>1$). However, ADMM with the MM principle has not been studied for the \emph{nonconvex} problem~\eqref{model}, 
to the best of our knowledge.  

When the linear coupling constraint is absent, the block coordinate descent (BCD) method is a standard approach to solve~\eqref{model}.
Razaviyayn et al.~\cite{Razaviyayn2013} proposed the block successive upper-bound minimization (BSUM) framework that employs the MM principle in each block update. By employing suitable surrogate functions in each block update, BSUM recovers the typical BCD methods, for example of  \cite{GRIPPO20001,Hildreth,Powell1973,Tseng2001,Beck2013,Bolte2014,Tseng2009}. In the non-convex setting, BCD methods with inertial terms\footnote{We use in this paper the terminology ``inertial" to mean that an inertial term that involves the current iterate and the previous iterates is added to the objective of the subproblem to  update each block, see \cite{Titan2020}.} have also been studied, and have showed significant improvement in their practical performance; \hien{see, e.g.,~\cite{Ochs2019} for inertial BCD methods with heavy-ball acceleration, \cite{Xu2013,Xu2017} for inertial BCD methods with Nesterov-type acceleration, and \cite{Pock2016,Hien_ICML2020} for inertial BCD methods that use two extrapolation points.}  
Recently, the authors in~\cite{Titan2020} proposed a general inertial block MM framework for solving~\eqref{model} 
without the linear coupling constraint. To the best of our knowledge, inertial ADMM with \hien{\emph{Nesterov-type acceleration for the primal variables}} have not been studied in the nonconvex case of~\eqref{model} although some variants of ADMM with inertial terms for the primal variables have been analysed in the convex case (\hien{that is, when both $F$ and $h$ are convex}); see e.g.,  \cite{Buccini2020,LiLin2019,Ouyang2015}.   
  
\hien{Recently, \cite{Sun2019} proposes ADMM with inertial term for the dual variable, see the description in~\cite[Expression (17)]{Sun2019}. We would like to remark that we realize a gap\footnote{Specifically, the second equality of \cite[Expression (51)]{Sun2019} is not correct.} in the proof of \cite[Lemma 5]{Sun2019}.} 
Let us also mention stochastic ADMM methods for solving Problem~\eqref{model} in which the objective is in expectation formulation, see, e.g., \cite{Huang19a,Huang2016}, which is out of the scope of this paper.

\vspace{-0.05in}
\subsection{Contribution and Organization}
In this paper, we propose iADMM, a framework of inertial alternating direction methods of multipliers, for solving the nonconvex nonsmooth problem~\eqref{model}. When no extrapolation is used, iADMM becomes a general ADMM framework that employs the minimization-majorization principle in each block update. For the first time in the \emph{nonconvex} nonsmooth setting of Problem~\eqref{model}, we study ADMM and its inertial version combined with the MM principle when updating each block of variables. 	
Moreover, our framework allows to use an over-relaxation parameter $\alpha \in (0,2)$ to set $\alpha \beta$ as the constant stepsize for updating the dual variable $\omega$. Note that $\alpha=1$, see, e.g.,  \cite{HongADMM2020,Li2015,Wang2019}, or  $ \alpha\in \big( 0, \frac{1+\sqrt{5}}{2}\big)$, see, e.g.,~\cite{Fazel2013,Yang2017}, are  the standard choices in the nonconvex setting.
Recently, \cite{Bot2020} proposed proximal ADMM that use $\alpha\in (0,2)$ for solving a special case of the nonconvex Problem~\eqref{model} with $s=1$ and $\mathcal A=-\mathcal I$.   
 
Under standard assumptions and $\alpha\in (0,2)$, we analyse the subsequential convergence for the generated sequence of iADMM and ADMM. 
When $F(x) + h(y)$ satisfies the Kurdyka-{\L}ojasiewicz (K{\L}) property and $\alpha=1$, we prove the global convergence \hien{and provide the convergence rate for the generated sequence. We would like to emphasize that although proving convergence towards a critical point has become a typical task when considering the nonconvex nonsmooth Problem~\eqref{model}, see e.g.,  \cite{Bot2020,Li2015,Wang2019}, the techniques to accomplish this task heavily depend on the considered algorithms and the involved assumptions. As far as we are aware of, this has not been done for ADMM used in combination with the MM principle and inertial terms for the primal variables.  
} 

Finally, we apply the proposed framework to solve a class of \revise{nonconvex low-rank representation to illustrate the efficacy of iADMM. More specifically, in order to illustrate the effect of MM procedure in Algorithm~\ref{alg:iADMM}, we use suitable surrogate functions such that each block of variables has a close-form update rule (thus, we do not need to use an outer optimization solver to find a solution for the corresponding subproblem), see details in Section \ref{sec:numerical_surrogate}. In order to illustrate the acceleration effect of Algorithm~\ref{alg:iADMM}, we also employ inertial terms and the extrapolation parameters are appropriately chosen to guarantee \emph{a global convergence}, see details in Section \ref{sec:numerical_parameter}. Indeed, the numerical results presented in Section~\ref{sec:numerical_experiment} (see also Appendix~\ref{sec:addexp} and Appendix~\ref{appendix:NMF}) empirically show the significant acceleration effect of using inertial terms. 
} 

\hien{The paper is organized as follows. In the next section, we describe the proposed method, iADMM, and analyse its convergence properties. 
In Section~\ref{sec:numerical}, 
we \revise{report the numerical results of iADMM} on a class of nonconvex low-rank representation problems.  
We conclude the paper in Section~\ref{sec:conclusion}.  
All the technical proofs are presented in Appendix~\ref{app:proofs}. }

\section{An inertial ADMM framework}
\label{sec:algorithms}

In this section, we describe the iADMM framework and prove its subsequential and global convergence. 
Throughout the paper, we make the following assumptions that are standard for studying Problem~\eqref{model} and the convergence of ADMMs in the nonconvex setting, see for example  \cite{Wang2019,Bot2020,Li2015}. 
\begin{assumption}
\label{assume}
(i) $\sigma_{\mathcal B}:=\lambda_{\min}(\mathcal B \mathcal B^*) >0 $.

(ii) 
$F(x) + h(y)$ is lower bounded. 

(iii) The function $h$ is $L_h$-smooth, that is, $\nabla h$ is $L_h$-Lipschitz continuous. 
\end{assumption}

\subsection{Description of iADMM} 
\label{sec:iADMM-description}

Let us first formally define a surrogate function. Some examples were given in the introduction. More examples can be found in~\cite{Mairal_ICML13,Razaviyayn2013,Titan2020}. 

\begin{definition}[Surrogate function]
\label{def:surrogate}
Let $\mathcal X\subseteq \mathbb R^{\mathbf n}$. 
A function $u:\mathcal X \times \mathcal X \to \mathbb R$ is called a surrogate function of a function $f$ on $\mathcal X$  
if the following two conditions are satisfied: 
\[
(a) \; u(z,z) = f(z) \text{ for all } z\in \mathcal X, 
 \quad \text{ and } \quad (b) \; 
 u(x,z) \geq f(x) \text{ for all }  x,z\in\mathcal X. 
\] 
\end{definition}
\vspace*{-0.05in}
As we are considering multi-block problems, we need the following definition of a block surrogate function, which is a generalization of Definition~\ref{def:surrogate}. 

\begin{definition}[Block surrogate function]
\label{def:surrogate-block} Let $\mathcal X_i\subseteq \mathbb R^{\mathbf n_i}$, $\mathcal X\subseteq \mathbb R^{\mathbf n}$. 
A function $u_i:\mathcal X_i \times \mathcal X \to \mathbb R  $ is called a block $i$ surrogate function of $f$ on $\mathcal X$
if the following conditions are satisfied:
\vspace*{-0.05in}
\begin{itemize}
    \item[(a)] $u_i(z_i,z) = f(z)$ for all $z\in \mathcal X$, 
    \item[(b)] $u_i(x_i,z) \geq f(x_i,z_{\ne i})$ for all $x_i\in\mathcal X_i$ and $z\in \mathcal X$, 
\end{itemize}
where  
    $(x_i,z_{\ne i})$ denotes $(z_1,\ldots,z_{i-1},x_i,z_{i+1},\ldots,z_s).$ 
The block approximation error is defined as $ e_i(x_i,z):=u_i(x_i,z) - f(x_i,z_{\ne i}).$
\end{definition}
A separability condition is necessary in~\cite[Def.~3]{Lu2016} 
for the surrogate function of $f$ (i.e., when fixing $z$, the surrogate function $u$ of $f$ satisfies $\hat u(x) = \sum_{i=1}^s \hat u_i (x_i)$, where $\hat u(x)= u(x,z)$ and $\hat u_i(x_i) = u_i(x_i,z)$) while our upcoming analysis does not require such a  condition. 
 
\begin{algorithm}[ht!]
\caption{iADMM, a general framework for solving Problem~\eqref{model}} 
\label{alg:iADMM} 
\begin{algorithmic} 
\STATE Choose $x^{0}= x^{-1}$, $y^0=y^{-1}$, $\omega^0$. Let $u_i$, $i\in [s]$, be a block $i$ surrogate functions of $f(x)$ on $\mathbb R^{\mathbf n}$.
\STATE For the choice of the extrapolation parameters, $\zeta_i^k$ and $\delta_k$, and of the parameters $\kappa_i$, $\alpha$ and $\beta$, see the paragraph ``Choosing parameters for iADMM" (page~\pageref{choosingparam}).  
\FOR{$k=0,\ldots$}
\FOR{ $i = 1,...,s$}
   \STATE\label{step-update} 
   Compute $\bar x_i^k=x_i^k + \zeta_i^k (x_i^k - x_i^{k-1})$, 
   and update block $x_i$ as follows 
   \begin{equation}
\label{eq:xi-updatenew} 
\begin{split}
 x_i^{k+1}\in\argmin_{x_i}\Big\{ u_i(x_i,x^{k,i-1})+g_i(x_i) 
         & \,+ \langle \mathcal A_i^*\big( \omega^k+\beta(\mathcal  A \bar x^{k,i-1} +\mathcal  By^k-b) \big),x_i\rangle\\
         &\,+\frac{\kappa_i\beta}{2}\|x_i-\bar x_i^k\|^2\Big\},
\end{split}
\end{equation}
   where $\kappa_i\geq \|\mathcal A_i^* \mathcal A_i\|$, and    
   $\bar x^{k,i-1}=(x^{k+1}_1,\ldots,x^{k+1}_{i-1},\bar x^{k}_i, x^{k}_{i+1},\ldots,x^k_s)$.
   \ENDFOR
   \STATE Compute $\hat y^k= y^k + \delta_k (y^k - y^{k-1})$, and update $y$ as follows  
   \begin{equation}
   \label{eq:y_update}
   \begin{split}
       y^{k+1}\in\argmin_{y}\Big\{ \langle \mathcal B^*\omega^k +\nabla h(\hat y^k),y \rangle + \frac{\beta}{2} \|\mathcal  A x^{k+1} + \mathcal  By-b\|^2+\frac{L_h}{2}\|y-\hat y^k\|^2\Big\}.
        \end{split}
   \end{equation}
   \STATE Update $\omega$ as follows  
      \begin{equation}
   \label{eq:omega_update}
   \omega^{k+1}=\omega^k  +\alpha\beta (\mathcal  A x^{k+1}+ \mathcal B y^{k+1} -b). 
   \end{equation}
\ENDFOR
\end{algorithmic}
\end{algorithm} 

The inertial alternating direction method of multipliers (iADMM) framework is described in Algorithm~\ref{alg:iADMM}.
iADMM cyclically update the blocks $x_1,\ldots,x_s$ and $y$. 
We use 
$x^{k,i}$ to denote $(x^{k+1}_1, \ldots,x^{k+1}_{i},x^{k}_{i+1},\ldots,x^{k}_s)$, let $x^{k,0}=x^k$ and $x^{k+1} = x^{k,s}$, where $k$ is the outer iteration index, and $i$ the cyclic inner iteration index ($i \in [s]$). 
The update of  block $x_i$ in~\eqref{eq:xi-updatenew} (note that $x^{k+1}_{ i}=x^{k,i}_{ i} $) means that iADMM chooses a surrogate function for $x_i\mapsto \mL(x_i,x^{k,i}_{\ne i},y^k,\omega^k)$, which is formed by summing a surrogate function of $x_i\mapsto f(x_i,x^{k,i}_{\ne i}) + g_i(x_i)$ and a surrogate function of $x_i\mapsto \varphi^k(x_i,x^{k,i}_{\ne i})$ where $\varphi^k$ is defined in  \eqref{eq:varphi}, then apply extrapolation to the latter surrogate function\footnote{\label{footnote_general}It is important noting that it is possible to embed the general inertial term $\mathcal G_i^k$ to the surrogate of $x_i\mapsto \mL(x_i,x^{k,i}_{\ne i},y^k,\omega^k)$ as in~\cite{Titan2020}. This inertial term  may also lead to the extrapolation for the block surrogate function of $ f(x)$ or for both the two block surrogates. However, to simplify our analysis, we only consider here the effect of the inertial term for the block surrogate of $\varphi^k(x)$.}. To update block $y$, as $h(y)$ is $L_h$-smooth, we apply Nesterov type acceleration on $h$ as in \eqref{eq:y_update}.     
Together with Assumption~\ref{assume}, we make the following standard assumption for $u_i$ throughout the paper.

\begin{assumption}
\label{assump:Lipschitz_ui}
(i) The block surrogate function $u_i(x_i,z)$ is continuous. 

(ii) Given $z\in \mathbb R^{\mathbf n}$, for $i\in [s]$, there exists a function $x_i\mapsto \bar e_i(x_i,z)$ such that $ \bar e_i(\cdot,z)$ is continuously differentiable at $z_i$, $\bar e_i(z_i,z)=0$, $\nabla_{x_i} \bar e_i(z_i,z)=0$, and the block approximation error $ e_i$  satisfies  \begin{equation}
\label{lemma:h_property} 
e_i(x_i,z) \leq \bar e_i(x_i,z) \;  \text{ for all } \;  x_i.
\end{equation} 
\end{assumption}
 Assumption~\ref{assump:Lipschitz_ui} (ii) is satisfied when we simply choose $u_i(x_i,z) = f(x_i,z_{\ne i})$ (i.e., $f(x_i,z_{\ne i})$ is a surrogate function of itself), or when $e_i(\cdot,z)$ is continuously differentiable at $z_i$ and $\nabla_{x_i} e_i(z_i,z)=0$, or when $e_i(x_i,z)\leq c \|x_i-z_i\|^{1+\epsilon} $ for some $\epsilon>0$ and $c>0$; see \cite[Lemma 3]{Titan2020}.  \hien{In the following, we provide some examples of block surrogate functions satisfying Assumption~\ref{assump:Lipschitz_ui}. 
\begin{itemize}
\item The block proximal surrogate function,  
see, e.g.,  \cite{Attouch2009,Attouch2013,Hien_ICML2020}, has the following form 
\begin{equation*}
\begin{array}{ll}
u_i(x_i,z) = f(x_i,z_{\ne i}) + \frac{\rho_i}{2}\|x_i - z_i\|^2,
\end{array}
\end{equation*}
where  $\rho_i>0$ is a scalar. We have $e_i(x_i,z)= \frac{\rho_i}{2} \|x_i-z_i\|^2$. In this case, $\bar{e}_i = e_i$. 
\item The Lipschitz gradient surrogate function, 
see, e.g.,~\cite{Xu2013,Xu2017,Hien_ICML2020}, 
has the form 
\begin{equation}
\label{Lipschitz_grad_surrogate}
\begin{array}{ll}
u_i(x_i,z) = f(z) + \langle\nabla_i f(z), x_i- z_i\rangle + \frac{\kappa_i L^{(z)}_i}{2}\|x_i - z_i\|^2,
\end{array}
\end{equation}
where $\kappa_i\geq 1$ and we assume $x_i\mapsto f(x_i,z_{\ne i})$ is differentiable and $\nabla_i f(x_i,z_{\ne i})$ is $L_i^{(z)}$-Lipschitz continuous (we note that $L^{(z)}_i$ may depend on $z$).
 We have 
$$
\nabla_{x_i} e_i(x_i,z)= L^{(z)}_i (x_i-z_i) + \nabla_i f(z)- \nabla_i f(x_i,z_{\ne i}). 
$$
Hence $\nabla_{x_i} e_i(z_i,z)=0$. In this case $\bar e_i= e_i$.  
\item The quadratic surrogate, see e.g.,~\cite{Emilie2016,Ochs2019}, has the following form
\begin{equation}
u_i(x_i,z) =  f(z) + \langle\nabla_i f(z), x_i- z_i\rangle + \frac{\kappa_i}{2}(x_i- z_i)^T H^{(z)}_i(x_i- z_i),
\end{equation}
where $\kappa_i\geq 1$ and we assume $f$ is twice differentiable,  $H^{(z)}_i$ is a positive definite matrix such that $(H^{(z)}_i - \nabla_i^2 f(x_i,z_{\ne i}))$ is positive definite (we note that $H^{(z)}_i$ may depend on $z$). Similarly, we also have $\bar e_i= e_i$ in this case. 
\end{itemize}
}

\paragraph{\textbf{Choosing parameters for iADMM.}}
\label{choosingparam} 

The parameters of iADMM include: $\alpha$ in~\eqref{eq:omega_update}, $\kappa_i$ and the extrapolation parameters $\zeta_i^k$ in the update~\eqref{eq:xi-updatenew} of block $x_i$, the extrapolation parameter $\delta_k$ in the update~\eqref{eq:y_update} of $y$, and the penalty parameter $\beta$. 
In the next section,  Proposition~\ref{prop:NSDP} provides the formulas for $\eta_i$ and $\gamma_i^k$ that involve $\beta$, $\kappa_i$, and $\zeta_i^k$, while  Proposition~\ref{prop:NSDP-fory} provides the formulas for $\eta_y$ and $\gamma_y^k$ that involve $\beta$ and $\delta_k$. 
To guarantee a subsequential convergence, we choose $\alpha\in (0,2)$, and the parameters  $\eta_y$, $\gamma_y^k$, $\eta_i$ and $\gamma_i^k$ satisfying the conditions of Proposition~\ref{prop:subsequential}; see Theorem~\ref{thrm:subsequential}. To guarantee a global convergence, we choose $\alpha=1$, use no extrapolation for $y$, and choose the other parameters to satisfy~\eqref{requirement-global}; see Theorem~\ref{thrm:global}. It is important noting that the convexity of  $x_i\mapsto u_i(x_i,z) + g_i(x_i)$  allows larger extrapolation parameters in the update of $x_i$ (Proposition~\ref{prop:NSDP}),  while the convexity of $h$ allows larger extrapolation parameters in the update of $y$ (Proposition~\ref{prop:NSDP-fory}).    

\begin{remark} 
As we target Nesterov-type acceleration to update of $y$ ($h$ is assumed to be $L_h$-smooth), we analyse the update rule as in~\eqref{eq:y_update} for $y$. Updating $y$  using  
      $y^{k+1} 
      \in 
      \argmin_{y} \mL(x^{k+1},y,\omega^k)$ 
 would work as well, and the convergence analysis of iADMM would be simplified by using the same rationale  to obtain subsequential as well as global convergence. We hence omit this case in our analysis. 
\end{remark}   

\subsection{Convergence analysis}

\paragraph{Assumptions.} Throughout the paper we assume Assumption~\ref{assume} and Assumption~\ref{assump:Lipschitz_ui} hold, and $\alpha\in (0,2)$. 
 
  Let $x^{k,i}$, $y^k$ and $\omega^k$ be the iterates generated by iADMM. We define some additional notations as follows. We denote $\Delta x^k_i= x^{k}_i -x^{k-1}_i$, $\Delta y^k=y^{k} -y^{k-1}$, $\Delta \omega^k =\omega^{k} -\omega^{k-1}$, $\alpha_1=\frac{|1-\alpha|}{\alpha  \sigma_{\mathcal B}(1-|1-\alpha|)}$, $ \alpha_2=\frac{3\alpha}{\sigma_{\mathcal B}(1-|1-\alpha|)^2}$ and $\mL^k=\mL(x^{k},y^k,\omega^k )$. We let $\nu_i$, $i\in [s]$, and $\nu_y$ be arbitrary constants in $(0,1)$. We take the following convention in the notation that allows us to analyse iADMM and its non-inertial version in parallel: 
\begin{itemize} 
\item If $\zeta_i^k = 0$ (i.e., there is no extrapolation in the update of $x_i^k$), then $\zeta_i^k/\nu_i=0$ and  $\nu_i=0$.
\item If $\delta_k =0$ (i.e., there is no extrapolation in the update of $y$), then $\delta_k/\nu_y=0$ and $\nu_y=0$. 
\end{itemize}
Now we present our main convergence results; see the proofs in Appendix~\ref{app:proofs}. 

As iADMM allows to use extrapolation in the update of $x_i^k$ and $y^k$, the Lagrangian is not guaranteed to decrease at each iteration. 
Instead, it has the following nearly sufficiently decreasing property as stated in the following Propositions~\ref{prop:NSDP} and~\ref{prop:NSDP-fory}.

\begin{proposition}
\label{prop:NSDP}

 (i)   Considering the update in~\eqref{eq:xi-updatenew},  \revise{in general (when $x_i\mapsto u_i(x_i,z) + g_i(x_i)$ can be nonconvex),}
 we choose $\kappa_i> \|A_i^*A_i\|$. Denote $a_i^k= \beta \zeta_i^k(\kappa_i+\| \mathcal A_i^* \mathcal A_i \|)$. Then we have
\begin{equation}
\label{NSDP}
\mL(x^{k,i},y^k,\omega^k)+\eta_i \|\Delta x^{k+1}_i\|^2
\leq
\mL(x^{k,i-1},y^k,\omega^k)+ \gamma_i^k \|\Delta x^k_i\|^2,
\end{equation}
where 
\begin{equation}
\eta_i = \frac{(1-\nu_i)(\kappa_i-\|\mathcal A_i^*\mathcal A_i\|)\beta}{2}, \quad
\gamma_i^k= \frac{(a^k_i)^2}{2\nu_i(\kappa_i-\|\mathcal A_i^*\mathcal A_i\|)\beta}. 
\end{equation}
(ii) If  $x_i\mapsto u_i(x_i,z) + g_i(x_i)$ is convex,  we take $\kappa_i =\|\mathcal A_i^* \mathcal A_i\|$ \revise{(note that if $\|\mathcal A_i^* \mathcal A_i\|=0$ then we can choose $\kappa_i$ as in case (i))}.  Inequality~\eqref{NSDP} is then satisfied with 
\begin{equation}
\label{eq:gammastrong}
\gamma^{k}_i =\frac{\beta \|\mathcal A_i^* \mathcal A_i\| (\zeta_i^k)^2}{2} , \quad \eta_i=\frac{\beta \|\mathcal A_i^* \mathcal A_i\|}{2}.
\end{equation}  
\end{proposition}

\begin{proposition}
\label{prop:NSDP-fory} 
Considering the update in~\eqref{eq:y_update}, we have 
\begin{equation*}
\label{eq:NSDP-y}
\mL(x^{k+1},y^{k+1},\omega^k)+ \eta_y \|\Delta y^{k+1}\|^2 \leq \mL(x^{k+1},y^k,\omega^k) + \gamma_y^k\|\Delta y^{k}\|^2,
\end{equation*}
where $\eta_y=\frac{(1-\nu_y)(\beta\lambda_{\min}(\mathcal B^* \mathcal B)+L_h) }{2}$ and $\gamma_y^k=\frac{2L_h^2 \delta_k^2}{\nu_y (\beta\lambda_{\min}(\mathcal B^* \mathcal B)+L_h)} $ when $h(y)$ is nonconvex, and $\eta_y=\frac{L_h }{2}$ and $\gamma_y^k= \frac{L_h \delta_k^2}{2}$ when $h(y)$ is convex. 
\end{proposition}

From Proposition~\ref{prop:NSDP} and Proposition~\ref{prop:NSDP-fory}, we obtain the following recursion for $\{\mL^k\}$. 

\begin{proposition}
\label{prop:recursive}
We have 
\begin{equation}
\label{eq:recursive2}
\begin{split}
&\mL^{k+1}+  \eta_y\|\Delta y^{k+1}\|^2 +\sum_{i=1}^s \eta_i \|\Delta x^{k+1}_i\|^2
\\ 
&\leq \mL^k +\sum_{i=1}^s \gamma_i^k  \|\Delta x^{k}_i\|^2 
+\gamma_y^k \|\Delta y^{k}\|^2 + \frac{\alpha_1}{\beta }  (\|B^* \Delta \omega^{k}\|^2-  \|B^* \Delta \omega^{k+1} \|^2)
\\
&\quad+  \frac{\alpha_2}{\beta}L_h^2 \|\Delta y^{k+1}\|^2+ \frac{\alpha_2}{\beta}\big( \bar \delta_kL_h^2   \|\Delta y^{k}\|^2 + 4L_h^2\delta_{k-1}^2\|\Delta y^{k-1}\|^2\big),
\end{split}
\end{equation}
where $\bar \delta_k = 2$ if $\delta_k=0$ for all $k$ and $4(1+\delta_k)^2$ otherwise.
\end{proposition}
Now we characterize the chosen parameters for Algorithm~\ref{alg:iADMM} in the following proposition.
\begin{proposition} 
\label{prop:subsequential}
Let $\eta_y$,  $\gamma_y^k$, $\eta_i$, and $\gamma_i^k$, $i\in [s]$, be defined in Proposition~\ref{prop:NSDP} and Proposition \ref{prop:NSDP-fory}. 
Denote 
$\mu=  \eta_y-\frac{\alpha_2L_h^2}{\beta}.$
 For $k\geq 1$, suppose the parameters are chosen such that $\mu>0$, $\eta_i>0$, and the following conditions are satisfied for some constants $0<C_x,C_y<1$:
\begin{equation}
\label{requirement}
\begin{split}
 \gamma_i^k \leq C_x \eta_i, \quad  \frac{ 4\alpha_2 L_h^2\delta_{k-1}^2}{\beta}\leq C_2 \mu,
\quad \frac{ \alpha_2 L_h^2\bar \delta_k}{\beta} + \gamma_y^k \leq  C_1 \mu,
\end{split}
\end{equation}
where 
$\begin{cases} C_1=C_y \,\text{and} \,C_2=0 &\mbox{if}\, \delta_k=0 \,\forall\,k,
\\
0<C_1<C_y \, \text{and}\, C_2=C_y-C_1 & \mbox{otherwise},
\end{cases} $  and $\bar \delta_k$ is defined in Proposition~\ref{prop:recursive}. Furthermore, suppose we use one of the following methods:
\begin{itemize}
\item we choose $\delta_k=0$ for all $k$, that is, there is no extrapolation in the update of $y$,
\item we use extrapolation in the update of $y$ and choose the parameters such that 
 \begin{equation}
\label{requirement2}
\beta \geq \frac{4L_h \alpha}{\sigma_{\mathcal{ B}} (1-|1-\alpha| )}, \quad
\beta \geq\frac{6\alpha L_h^2 }{\mu  \sigma_{\mathcal B} (1-|1-\alpha|)}  \max\Big\{1,\frac{12\delta_k^2}{1-C_1}\Big \}.
\end{equation}  
\end{itemize}

(i) For $K>1$ we have 
\begin{equation}
\label{eq:upperbounded}
\begin{split}
&\mL^{K+1} +  \mu \|\Delta y^{K+1} \|^2 + \sum_{i=1}^s\eta_i \| \Delta x^{K+1}_i \|^2 
+ \frac{\alpha_1}{ \beta }  \|\mathcal B^* \Delta w^{K+1}\|^2 + (1-C_1)\mu \|\Delta y^{K}\|^2\\
&\, + \sum_{k=1}^{K-1}\big[(1-C_y)\mu \|\Delta y^{k}\|^2  + (1-C_x)\sum_{i=1}^s\eta_i \|\Delta x^{k+1}_i\|^2 \big] \\
&\leq \mL^1+ \frac{\alpha_1}{\beta  }  \|\mathcal B^* \Delta \omega^{1}\|^2 + C_x\sum_{i=1}^s \eta_i  \|\Delta x^{1}_i\|^2 +  \mu \|\Delta y^{1}\|^2 + C_2 \mu \|\Delta y^{0}\|^2.
\end{split}
\end{equation}

(ii) The sequences $\{\Delta y^{k}\}$, $\{\Delta x^{k}_i\}$ and $\{\Delta \omega^{k}\}$  converge to 0. 
\end{proposition}

We will assume that Algorithm~\ref{alg:iADMM} generates a bounded sequence in our subsequential and global convergence results. Let us provide a sufficient condition that guarantees this boundedness assumption. 

\begin{proposition}
\label{prop:bounded-sequence}
If $b+ Im(\mathcal A) \subseteq Im(\mathcal B)$, $\lambda_{\min}(\mathcal B^* \mathcal B)>0 $ and $F(x) + h(y)$ is coercive over the feasible set $\{(x,y): \mathcal  Ax + \mathcal By = b\}$ then the sequences $\{x^k\}$, $\{y^k\}$  and $\{\omega^k\}$ generated by Algorithm~\ref{alg:iADMM} are bounded.
\end{proposition}

It is important noting that the coercive condition of $F(x) + h(y)$ over the feasible set is weaker than the coercive condition of $F(x) + h(y)$ over  $x\in\mathbb R^{\mathbf n}, y\in \mathbb R^{\mathbf q}$. \hien{Let us now present the subsequential convergence, the global convergence of the generated sequence and its convergence rate}. 

\begin{theorem}[Subsequential convergence]
\label{thrm:subsequential}
Suppose the parameters of Algorithm~\ref{alg:iADMM} are chosen satisfying the conditions in  Proposition~\ref{prop:subsequential}. If the generated sequence of Algorithm~\ref{alg:iADMM} is bounded, then every limit point of the generated sequence is a critical point of $\mL$. 
\end{theorem}

\begin{theorem}[Global convergence]
\label{thrm:global}
Suppose we do not use extrapolation to update $y$, that is, $\delta_k=0$ for all $k$ (note that extrapolation to update $x_i$ is still applicable), and we take $\alpha=1$. Then the conditions in  \eqref{requirement} become 

\begin{equation}
\label{requirement-global}
\gamma_i^k \leq C_x \eta_i, \quad
\frac{2 \alpha_2 L_h^2}{\beta}  \leq  C_y \mu, \, \text{for all}\, k \geq 0, i \in [s] 
\end{equation}
for some constants $0< C_x,C_y <1$. 
Furthermore, we assume that (i) for any $x,z\in \mathbb R^{\mathbf n}$, $x_i\in {\rm dom} (g_i)$, 
we have
\begin{equation}
\label{assume:partial}
\begin{array}{ll}
\partial_{x_i} \big(f(x) + g_i(x_i)\big)=\partial_{x_i} f(x) + \partial_{x_i} g_i(x_i), \\
\partial_{x_i} \big(u_i(x_i,z) + g_i(x_i) \big)= \partial_{x_i} u_i(x_i,z) + \partial_{x_i} g_i(x_i),
\end{array}
\end{equation} 
 and (ii) for any $x,z$ in a bounded subset of $\mathbb R^{\mathbf n}$, if $\mathbf s_i\in\partial_{x_i} u_i(x_i,z)$, there exists $\xi_i\in \partial_{x_i} f(x)$ such that

\begin{equation}
\label{eq:l1}
\|\xi_i - \mathbf s_i\| \leq  L_i\|x-z\| \, \text{for some constant}\, L_i.  
\end{equation}
If the generated sequence of Algorithm~\ref{alg:iADMM} is bounded and $F(x) + h(y)$ has the K{\L} property (see Appendix~\ref{sec:pre}),  then the whole generated  sequence converges to a critical point of $\mL$.
\end{theorem}

We refer the readers to \cite[Corollary 10.9]{RockWets98} for a sufficient condition for~\eqref{assume:partial} (see Appendix~\ref{sec:pre} for more details). Some specific examples that satisfy \eqref{assume:partial} include: (i) $g_i=0$, (ii) the functions $x_i\mapsto f(x)$ and $x_i\mapsto u_i(x_i,z)$ are strictly differentiable (see \cite[Exercise 10.10]{RockWets98}), (iii) the functions $x_i\mapsto f(x)$ and $x_i\mapsto u_i(x_i,z)$ are convex and the relative interior qualification conditions are satisfied: ${\rm ri} ({\rm dom} (f(\cdot,x_{\ne i})) \cap {\rm ri} ({\rm dom} g_i) \ne \emptyset $ and ${\rm ri} ({\rm dom} (g(\cdot,z)) \cap {\rm ri} ({\rm dom} g_i) \ne \emptyset $, where ${\rm ri}$ is short for relative interior. We note that  although the condition in \eqref{eq:l1} is necessary for our convergence proof, the constant $L_i$ does not influence how to choose the parameters in our framework. The condition in \eqref{eq:l1} is satisfied when both $u_i$ and $f$  are twice continuously differentiable and $\nabla_{x_i} e_i(x_i,x)=0$  for all $x$ (which implies that  $\nabla_{x_i} u_i(x_i,x)= \nabla_{x_i} f(x)$ for all $x$). Indeed, in this case we have $$\| \nabla_{x_i} u_i(x_i,z)-\nabla_{x_i} f(x) \|=\| \nabla_{x_i} u_i(x_i,z)-\nabla_{x_i} u_i(x_i,x)\| \leq L_i \|x-z\|$$ for some $L_i$ because $\nabla_{x_i} u_i (x_i,z)$ is continuously differentiable and thus is Lipschitz continuous over any bounded subset. We note that all the examples given after Assumption~\ref{assump:Lipschitz_ui} in Section~\ref{sec:algorithms}  satisfy the condition in \eqref{eq:l1} when $f$ is twice continuously differentiable.

 \paragraph{Convergence rate} A convergence rate for the generated sequence  of iADMM can be derived using the same technique as in~\cite[Theorem~2]{Attouch2009}.  To the best of our knowledge, in the nonconvex setting, the convergence rate for block coordinate methods (including inertial as well as non-inertial algorithms) appears to be the same in different papers in the literature since all papers use the technique in \cite{Attouch2009}.  As it is similar to establish the rate for iADMM, we omit the details. Instead, we refer the readers to \cite[Theorem~2.9]{Xu2013} and  \cite[Theorem~3]{Hien_ICML2020} for some examples of using this technique to establish the convergence rate. 
 The type of the convergence rate depends on the value of the K{\L}  exponent, which is the coefficient $\mathbf a$ such that $\Upsilon(t)$ in Definition~\ref{def:KL} 
 (see Appendix~\ref{sec:pre}) equals $c t^{1-\mathbf a}$, where $c$ is a constant. 
 Specifically, when $\mathbf a=0$, the algorithm converges after a finite number of  steps, when $\mathbf a\in (0,1/2]$, the algorithm has linear convergence, and when $\mathbf a \in (1/2, 1)$, the algorithm has sublinear convergence. 
 Determining the value of the K{\L}  exponent is out of the scope of this paper, and is an active and challenging topic. 

\vspace{-0.05in}
\section{Numerical results}
\label{sec:numerical}

\revise{
In this section, we apply iADMM to solve a latent low-rank representation problem. We consider Problem~\eqref{LRR} with 
\begin{itemize}
    \item $r_1(t) = \lambda_1 t$ to promote $X$ to be of low-rank, 
    since $r_1( \sigma(X) ) =\lambda_1 \sum_i \sigma_i(X) 
    = \lambda_1 \|X\|_*$ is the nuclear norm~\cite{recht2010guaranteed}.  
    
    \item $r_2(Y) = \lambda\sum_{i=1}^q\phi(\|Y_i\|_2)$, where $\phi(t) = 1 - \exp(-\theta t)$ is concave, $\theta > 0$ is a parameter, and $Y_i$ is the $i$-th column of $Y$. 
    This is a \emph{nonconvex} regularization that promotes $Y$ to be column sparse, that is, it promotes $Y$ to have many columns equal to the zero vector~\cite{brafea}. 
    In fact,  $\phi(t) = 0$ when $t = 0$, while 
    $\phi(t)$ quickly goes to 1 as $t$ increases. 
    
    \item  $r_3(Z) = \frac{1}{2}\|Z\|^2$  to model Gaussian noise. 
    
    \item $A_1 = DP_1$ and $A_2 = P_2^*D$, where $P_1$ and $P_2$ are computed by orthogonalizing the columns of $D^*$ and $D$, respectively, as proposed in~\cite{Liu2011}. In \cite{Liu2011}, the authors showed that this resulting problem is a simpler equivalent form of the one in which $D$ is considered as a dictionary, i.e. $A_1 = A_2 = D$. Hence, it can be scaled for data sets with a large number of observations. 
    
\end{itemize} 
In this scenario, Problem~\eqref{LRR} takes the form of \eqref{model} with  $B$ being the identity operator, $b$ being the data set $D$, $x_1$ and $x_2$ being the matrices $X$ and $Y$, $y$ being the matrix $Z$, $g_i = 0$, $h(Z) =  \frac{1}{2}\|Z\|^2$ and $f(X,Y) = \lambda_1\|X\|_* +  r_2(Y)$.  
}

\subsection{Surrogate functions and iADMM updates.}
\label{sec:numerical_surrogate}
We choose  $u_1(X,X^k,Y^k) =  \lambda_1\|X\|_* +  r_2(Y^k)$, and $u_2(Y,X^{k+1},Y^k) = r_2(Y^k) +  \sum_{i=1}^q\varsigma_i^k(\|Y_i\|_2-\|Y^k_i\|_2) +   \lambda_1\|X^{k+1}\|_*$, 
where $\varsigma_i^k = \lambda\nabla\phi(\|Y_i^k\|_2)$. 
The function $u_1$ satisfies Assumption \ref{assump:Lipschitz_ui}, 
and $u_2$ satisfies Assumption \ref{assump:Lipschitz_ui}~(i). 
Since $\phi$ is continuously differentiable with Lipschitz gradient on $[0,+\infty)$, and the Euclidean norm is Lipschitz continuous, it follows from Section 4.5 of \cite{Titan2020} that  $u_2$ also satisfies 
Assumption~\ref{assump:Lipschitz_ui}~(ii). We derive  from~\cite[Corollary 5Q]{Rock1981} that the condition in~\eqref{eq:l1} is satisfied.  
According to the update \eqref{eq:xi-updatenew}, $X^{k+1}$ is computed by solving the following nuclear norm problem\vspace*{-0.05in}
\begin{equation}\label{subLRR}
\min_X  \lambda_1\|X\|_* + \bigg\langle A_1^* \Big(\beta(A_1\bar X^k+ Y^kA_2 + Z^k  - D)
 + W^k\Big), X\bigg\rangle 
 +   \frac{\kappa_1\beta}{2}\|X - \bar X^k\|^2,
\end{equation}
where $\kappa_1 \geq \|A_1^*A_1\|$ and $\bar X^k = X^k + \zeta_1^k(X^k - X^{k-1})$. Let $\text{diag}(u)$ denote a diagonal matrix whose diagonal elements are the entries of $u$, and $[.]_+$ denote the projection onto the nonnegative orthant.  Problem \eqref{subLRR} has a closed-form solution given by $ X^{k+1} = US_{\lambda_1/(\kappa_1\beta)} V^T,
$
where $USV^T$ is the SVD of $\bar X^k  - A_1^*(A_1\bar X^k+ Y^kA_2 + Z^k - D+ W^k)/(\kappa_1\beta)$ and $S_{\lambda_1/(\kappa_1\beta)}  = \text{diag}([S_{ii} - \lambda_1/(\kappa_1\beta)]_+)$. 
Letting $\kappa_2 \geq \|A_2A_2^*\|$ and $\bar Y^k = Y^k + \zeta_2^k(Y^k - Y^{k-1})$, 
the update \eqref{eq:y_update} for $Y$ is 

\begin{align*}
Y^{k+1} &\in\arg\min_Y  \sum_{i=1}^q\varsigma_i^k\|Y_i\|_2+ \langle (W^k + \beta(A_1X^{k+1} + \bar Y^kA_2\\
&\qquad\qquad +  Z^k - D))A_2^*, Y\rangle 
+ \frac{\kappa_2\beta}{2}\|Y  - \bar Y^k \|^2. 
\end{align*}
It has a closed-form solution given by
$Y^{k+1}_{i}=
\left[ 
\|P_{i}^k\| - \varsigma_i^k/(\kappa_2\beta) \right]_+\frac{P_{i}^k}{\|P_{i}^k\|}, 
$
where $P_{i}^k$ is the $i$-th column of $\bar Y^k - (A_1X^{k+1} + \bar Y^kA_2 + Z^k - D)/\kappa_2 - W^k/(\kappa_2\beta)$. \\ 
The update \eqref{eq:y_update}  for  $Z$ is $Z^{k+1} = -(W^k +\beta(A_1X^{k+1} + Y^{k+1}A_2 - D))/(1+\beta)$, and the update \eqref{eq:omega_update} for $W$ is $W^{k+1} = W^k + \alpha\beta(A_1X^{k+1} + Y^{k+1}A_2 + Z^{k+1}- D)$.

\subsection{Choosing parameters}
\label{sec:numerical_parameter}
We have $L_h = 1$, $\sigma_{\mathcal B} = 1$, and $\delta_k = 0$.  As $h(Z)$ is convex and we do not apply extrapolation for $Z$, by Proposition~\ref{prop:NSDP-fory}, 
$\eta_y=\frac12$ and $\gamma_y^k=0$. Since $\|X\|_*$ and $\sum_{i=1}^q\varsigma_i^k\|Y_i\|_2$ are convex, we choose  $\kappa_1 = \|A_1^*A_1\|$, $\kappa_2 = \|A_2A_2^*\|$, and the conditions in~\eqref{requirement-global} become
$\zeta_i^k  \leq \sqrt{C_x}$ ($i=1,2$) and  $\frac{(2+C_y)\alpha_2}{\beta}  
\leq \frac{C_y}{2}$. 
We take 
$C_x= 1 - 10^{-15}$, 
$\alpha = 1$, 
{$C_y = 1 - 10^{-6}$}, 
$\beta = 2(2 + C_y)\alpha_2/C_y$, $a_0=1$, 
\mbox{$a_k=\frac12(1+\sqrt{1+4a_{k-1}^2})$},  and 
$\zeta_i^k =\min\Big\{\frac{a_{k-1}-1}{a_k},\sqrt{C_x }\Big\}$. We set $\alpha = 1$ as we target global convergence. We have also conducted experiments with other values of $\alpha$ (namely $0.5, 1.4$ and $1.8$); see Appendix~\ref{sec:addexp}. 
\begin{figure*}[h!]
\begin{center}
\begin{tabular}{cc}
\includegraphics[width=0.49\textwidth]{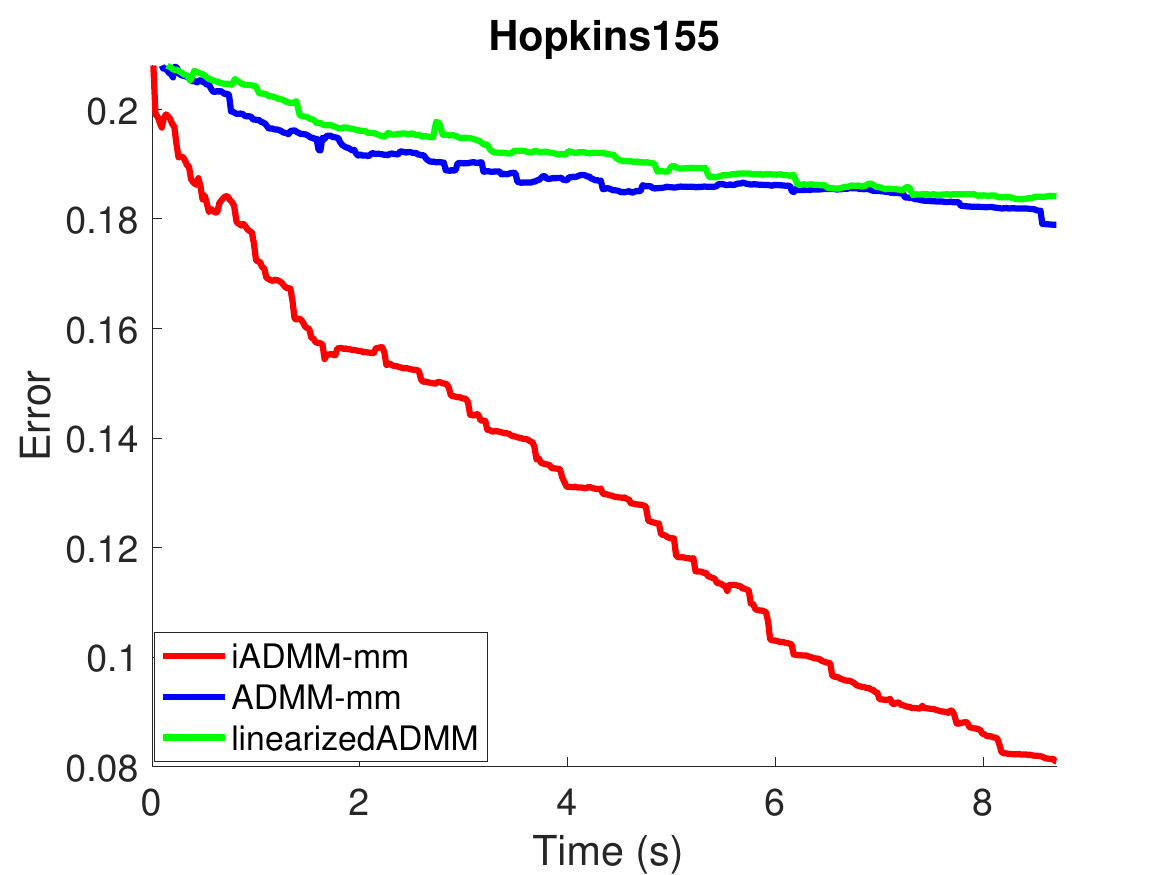}  & 
\includegraphics[width=0.49\textwidth]{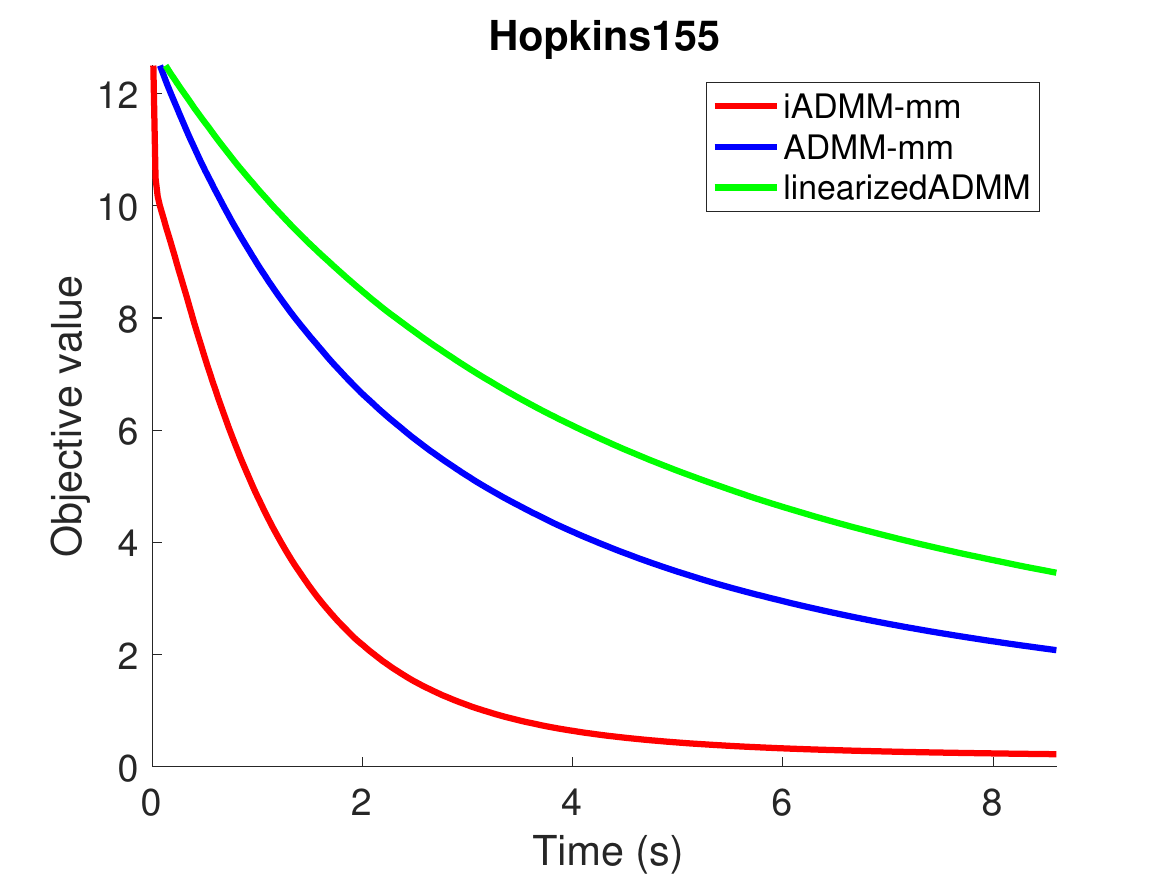} \\
\includegraphics[width=0.49\textwidth]{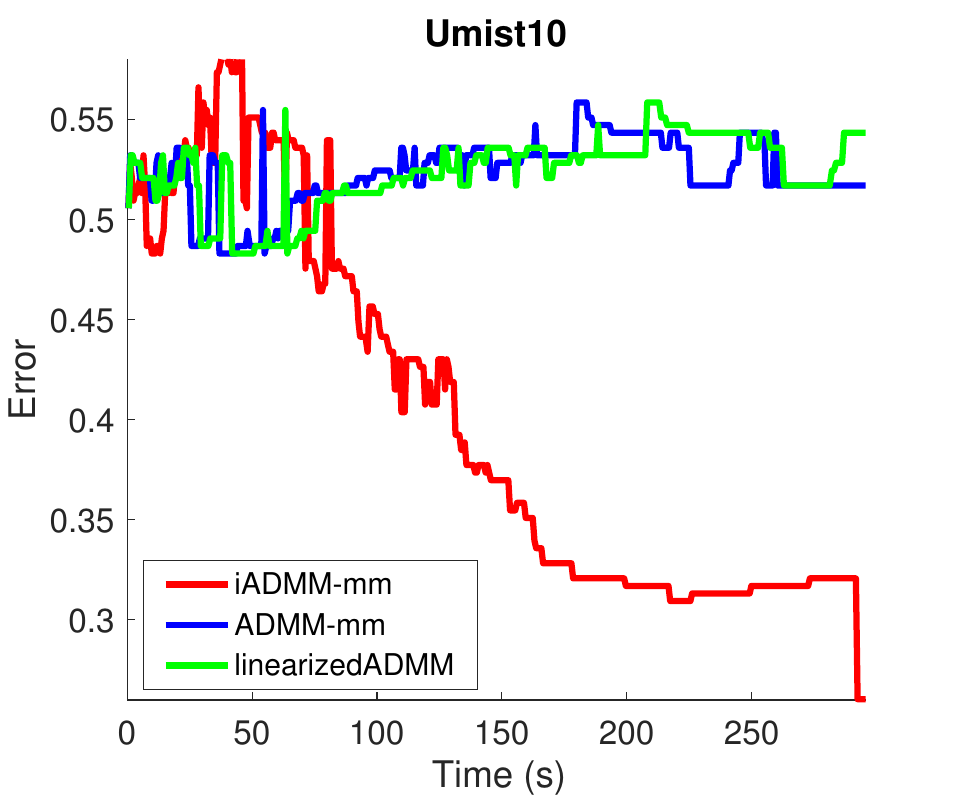}  & 
\includegraphics[width=0.49\textwidth]{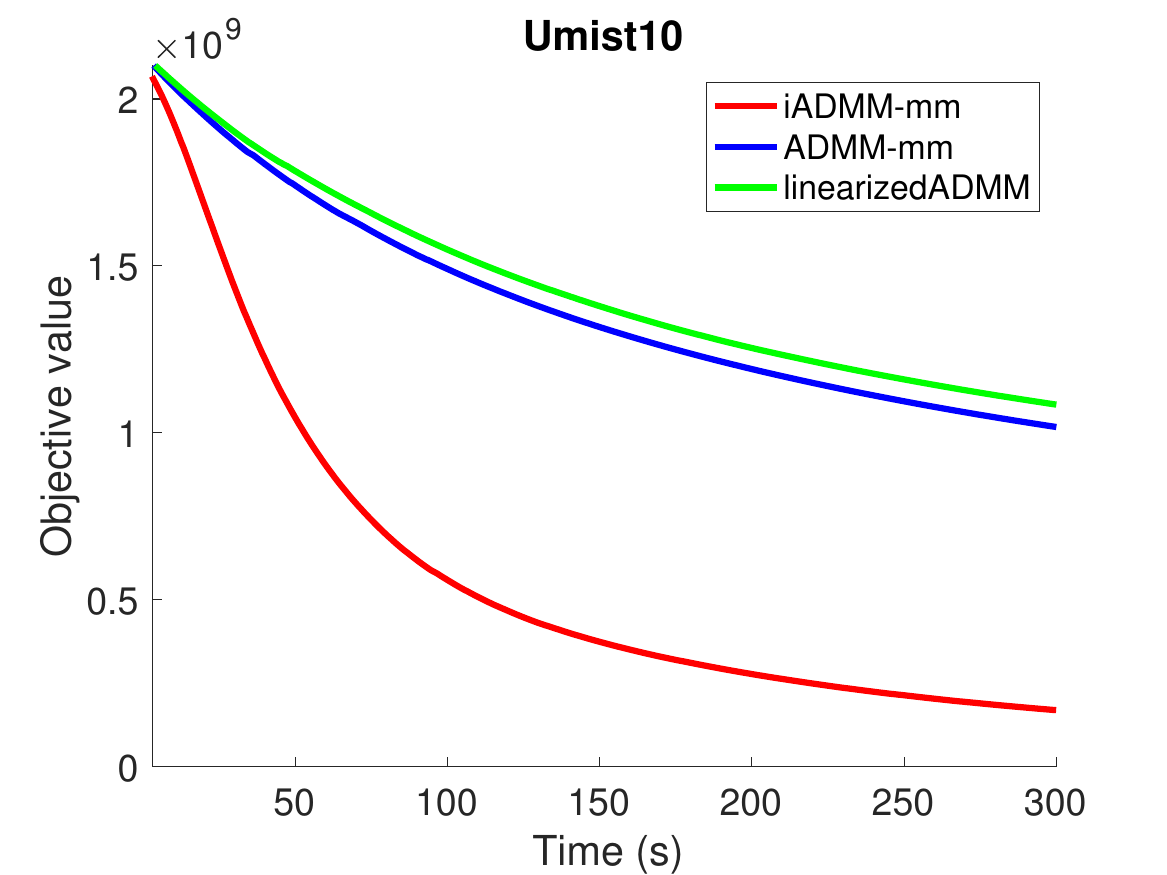} \\
\includegraphics[width=0.49\textwidth]{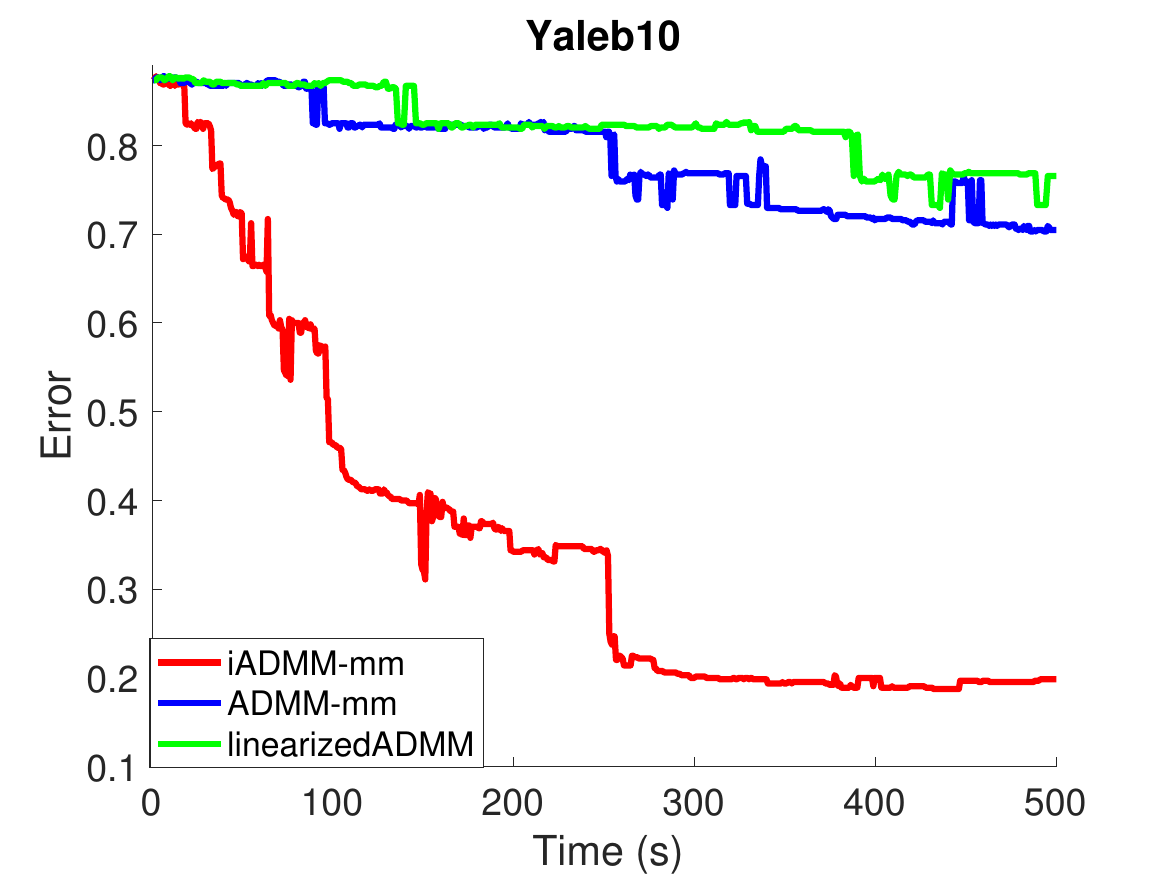}  & 
\includegraphics[width=0.49\textwidth]{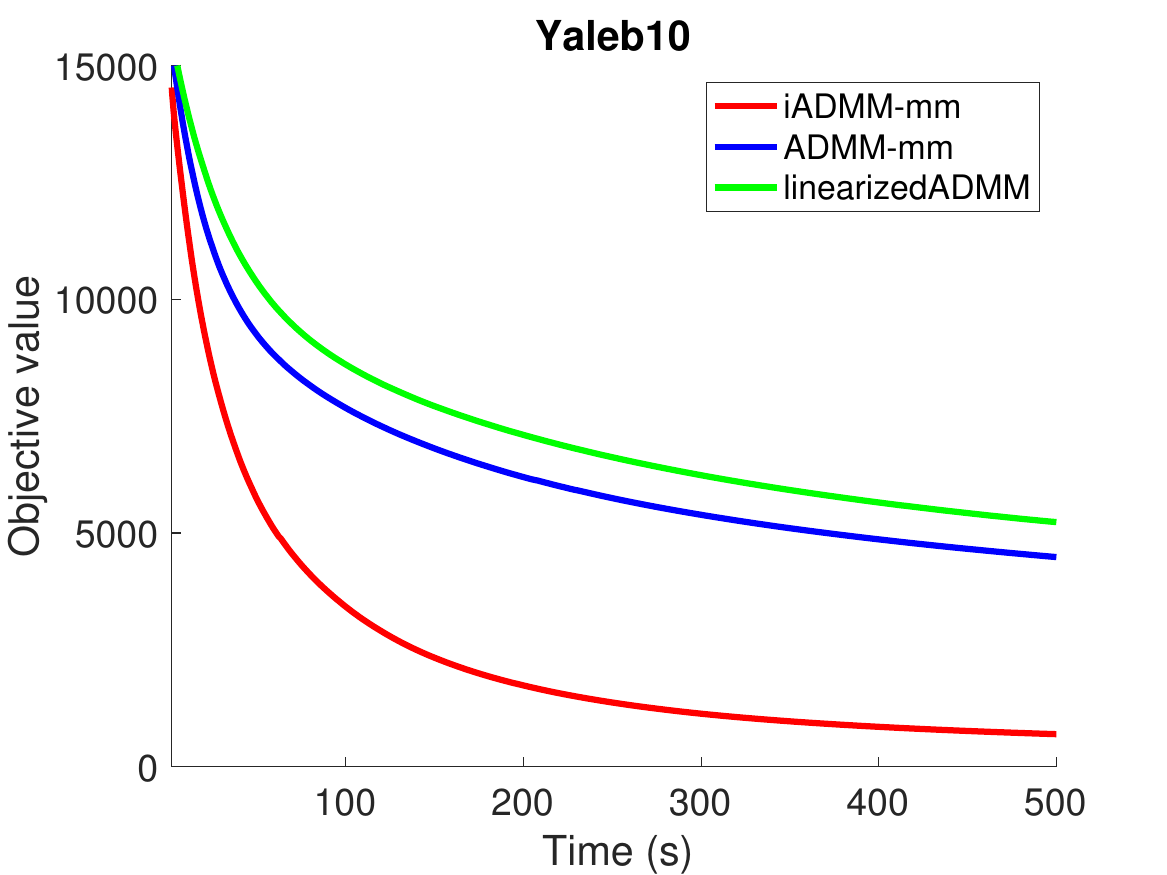}
\end{tabular}
\caption{Evolution of the segmentation error rate and the objective function value with respect to time. For Hopkins155, the results are the average values over 156 sequences.
\label{fig:model1}} 
\end{center}
\end{figure*} 

\subsection{Experiments}
\label{sec:numerical_experiment} 
We compare the following three methods: 
(1) ADMM-mm: iADMM without extrapolation, 
(2) iADMM-mm: iADMM with extrapolation,
(2) linearizedADMM: a linearized ADMM which is  different from ADMM-mm for updating $Y$. linearizedADMM updates $Y$ by solving 
$\min -\lambda\exp(\|Y_i\|_2) + \frac{\kappa_2\beta}{2}\|Y_i - V_i^k\|^2$,
where $V^k_i$ is the $i$-the column of $X^k - (W^k + \beta(A_1X^{k+1} + \bar Y^kA_2 +  Z^k - D))A_2^*/(\kappa_2\beta)$. Since these sub-problems do not have closed-form solutions, we employ {an} MM scheme to solve them. 
To examine the performance of the three algorithms, we consider subspace segmentation tasks. 
After obtaining a solution $X^*$, we follow the setting in~\cite{Lui2013} to construct the affinity matrix $Q$ by $Q_{ij} = (\tilde{U}\tilde{U}^T)_{ij}$, where $\tilde{U}$ is formed by $U^*(\Sigma^*)^{1/2}$ with normalized rows and $U^*\Sigma^*(V^*)^T$ being the SVD of $X^*$. Finally, we apply the Normalized Cuts \cite{Shi2000} on $Q$ to cluster the data into groups. The experiments are run on three data sets: Hopkins 155, extended Yale B and Umist. Hopskins 155 consists of 156
sequences, each of which has from 39 to 550 vectors drawn from two or three motions (one motion corresponds to one subspace). Each sequence
is a sole segmentation task and thus there are 156
clustering tasks in total. Yale B contains 2414 frontal face
images with 38 classes, and Umist contains
564 images with 20 classes. To avoid computational issues when computing the segmentation error rate, we construct clustering tasks by using the first 10 classes of these two data sets~\cite{Lu2015}. 
All tests are preformed using Matlab
R2019a on a PC 2.3 GHz Intel Core i5
of 8GB RAM. \\  
In our experiments, we choose $\theta = 5$, $\lambda_1 = \lambda  = 0.01$ for Hopkins 155, and $\lambda_1 = \lambda  = 1$ for the two other data sets. We set the initial points to zero, that is, $X^0=0$, $Y^0=0$, $Z^0=0$, $W^0 = 0$.  
We do not optimize numerical results by tweaking the parameters and initial points as this is beyond the scope of this work. 
It is important noting that we evaluate the algorithms on the same models with the same initializations.   
We run each algorithm 10, 300, and 500 seconds for each sequence of Hopkins 155, Umist10, and Yaleb10, respectively. 
Figure~\ref{fig:model1} displays the values 
of the segmentation error rate and the objective function  
versus the training time, and Table~\ref{results} 
reports the final values. Since there are 156 sequences (data sets) in Hopkins 155, we plot the average values, and report the final average results and standard deviation over  these sequences. 
We observe that iADMM-mm converges the fastest on all the data sets, providing a significant acceleration of ADMM-mm. iADMM-mm achieves not only the best final objective function values but also the best segmentation error rates. 
This illustrates the usefulness of the acceleration technique.  
In addition, ADMM-mm outperforms linearizedADMM which illustrates the usefulness of properly choosing a proper  surrogate function. The conclusions are the same for other values of $\alpha$; 
see Appendix~\ref{sec:addexp}. 

\begin{table}[h!]
\centering
\caption{Comparison of segmentation error rate and final objective function values obtained within the allotted time. 
Bold values indicate the best results.}\label{results}
\begin{small}
\begin{tabular}{@{}clll@{}}
\toprule
     & \multirow{2}{*}{Method} & Error & Obj. value \\
& & mean $\pm$ std & mean $\pm$ std \\ 
\midrule
   \multirow{3}{*}{\rotatebox[origin=c]{90}{Hopkins}} & linearizedADMM & 0.1579 $\pm$ 0.1550 & 3.0254 $\pm$ 2.4189 \\
    & ADMM-mm & 0.1472 $\pm$  0.1513  &  1.8081  $\pm$  1.6674 \\
  & iADMM-mm  & $\mathbf{0.0562  \pm  0.1006}$ & $\mathbf{0.2023  \pm   0.1062}$ \\
   \hline
   \multirow{3}{*}{\rotatebox[origin=c]{90}{Umist}} &  linearizedADMM &  0.5170  &  1.0838$\times10^9$ \\
  & ADMM-mm  & 0.5170 & 1.0167$\times10^9$ \\
  &  iADMM-mm  & \textbf{0.2604} & \textbf{0.1694}$\times10^9$ \\
  \hline
  \multirow{3}{*}{\rotatebox[origin=c]{90}{Yaleb}} &  linearizedADMM &  0.7656  &  5.2317$\times10^3$ \\
  & ADMM-mm  & 0.7047 & 4.4829$\times10^3$ \\
  &  iADMM-mm  & \textbf{0.1984} & \textbf{0.6951}$\times10^3$ \\
   \hline
\end{tabular}
\end{small}
\end{table}

\vspace{-0.1in}
\section{Conclusion}
\label{sec:conclusion}
We have proposed and analysed iADMM, a framework of inertial alternating direction methods of multipliers, for solving a class of nonconvex nonsmooth optimization problem with linear constraints. The preliminary computational results in solving a class of nonconvex low-rank representation problems not only show the efficacy of using inertial terms for ADMM but also show the advantage of using suitable block surrogate functions that provide closed-form solutions in the block update of ADMM. We conclude the paper by mentioning two important questions that we consider as a future research directions: 
(i) 
Can we extend the cyclic update rule of iADMM to randomized/non-cyclic setting? 
(ii)  To guarantee the global convergence, iADMM does not allow extrapolation in the update of $y$; see Theorem~\ref{thrm:global}. Can we extend the analysis to allow the extrapolation in the update of $y$? 

\section*{Declarations}

\paragraph{Availability of data and material, and Code availability} The data and code are available from \url{https://github.com/nhatpd/iADMM}.

\paragraph{Funding} LTKH and NG acknowledge the support by the European Research Council (ERC starting grant no 679515), and by the Fonds de la Recherche Scientifique - FNRS and the Fonds Wetenschappelijk Onderzoek - Vlaanderen (FWO) under EOS Project no O005318F-RG47. NG also acknowledges the Francqui Foundation.

\paragraph{Conflicts of interest/Competing interests} Not applicable.


\bibliography{iADMM}
\bibliographystyle{spmpsci} 

\newpage 

\appendix
\begin{center}
\textbf{APPENDIX
}\end{center}
\vspace{-0.25in}

\section{Preliminaries of non-convex non-smooth optimization}
\label{sec:pre} 

\label{sec:prelnnopt}

In this appendix, we recall some basic definitions and results, namely 
directional derivative and subdifferentials in Definition~\ref{def:dd}, 
critical point in Definition~\ref{def:type2}, 
the subdifferential of a sum of function in Proposition~\ref{propRock}, 
and K{\L} functions in Definition~\ref{def:KL}.

Let $g: \bbE\to \bbR\cup \{+\infty\} $ be a proper lower semicontinuous function.  
\begin{definition} \cite[Definition 8.3]{RockWets98}
\label{def:dd}
\begin{itemize}
\item[(i)] For any $x\in{\rm dom}\,g,$  and $d\in\bbE$, we denote the directional derivative of $g$ at $x$ in the direction $d$ by 
\[g'\lrpar{x;d}=\liminf_{\tau \downarrow 0}\frac{g(x+\tau d)-g(x)}{\tau}. 
\] 
\item[(ii)] For each $x\in{\rm dom}\,g,$ we denote $\hat{\partial}g(x)$ as
the Frechet subdifferential of $g$ at $x$ which contains vectors
$v\in\mathbb{E}$ satisfying 
\[
\liminf_{y\ne x,y\to x}\frac{1}{\left\Vert y-x\right\Vert }\left(g(y)-g(x)-\left\langle v,y-x\right\rangle \right)\geq 0.
\]
If $x\not\in{\rm dom}\:g,$ then we set $\hat{\partial}g(x)=\emptyset.$  
\item[(iii)] The limiting-subdifferential $\partial g(x)$ of $g$ at $x\in{\rm dom}\:g$
is defined as follows: 
\[
\partial g(x) := \left\{ v\in\mathbb{E}:\exists x^{(k)}\to x,\,g\left(x^{(k)}\right)\to g(x),\,v^{(k)}\in\hat{\partial}g\left(x^{(k)}\right),\,v^{(k)}\to v\right\} .
\]
\item[(iv)] The horizon subdifferential $\partial^{\infty} g(x)$ of $g$ at $x$ is defined as follows:
\begin{align*}
\partial^{\infty} g(x)
&:= \Big\{ v\in\mathbb{E}:\exists \lambda^{(k)}\to 0,  \lambda^{(k)}\geq0, \lambda^{(k)} x^{(k)}\to x,\,g(x^{(k)})\to g(x),\\
&\qquad\,v^{(k)}\in\hat{\partial}g(x^{(k)}),\,v^{(k)}\to v\Big\} .
\end{align*}
\end{itemize}
\end{definition}
\begin{definition}
\label{def:type2}
We call $x^{*}\in \rm{dom}\,F$ a critical point of $F$ if $0\in\partial F\left(x^{*}\right).$ 
\end{definition}

\begin{definition}\cite[Definition 7.5]{RockWets98}
A function $f:\mathbb R^{\mathbf n} \to \mathbb R \cup \{+\infty\}$ is called subdifferentially regular at $\bar x$ if $f(\bar x)$ is finite and the epigraph of $f$ is Clarke regular at $ (\bar x, f(\bar x))$ as a subset of $ \mathbb R^{\mathbf n} \times \mathbb R$ (see \cite[Definition 6.4]{RockWets98} for the definition of Clarke regularity of a set at a point).
\end{definition}
\begin{proposition} \cite[Corollary 10.9]{RockWets98} \label{propRock}
Suppose $f=f_1 +\cdot + f_m$ for proper lower semi-continuous function $f_i:\mathbb R^{\mathbf n}\to \mathbb R\cup \{+\infty\}$ and let $\bar x \in {\rm dom} f$. Suppose each function $f_i$ is subdifferential regular at $\bar x$, and the condition that the only combination of vector $\nu_i \in \partial^{\infty} f_i(\bar x) $ with $ \nu_1 + \ldots \nu_m=0$ is $\nu_i=0$ for $i\in [m]$. Then we have 
$$ \partial f(\bar x) =  \partial f_1(\bar x) + \ldots  \partial f_m(\bar x).
$$
\end{proposition}
To obtain a global convergence, we need the following Kurdyka-{\L}ojasiewicz (K{\L}) property for $F(x) + h(y)$. 
\begin{definition}
\label{def:KL}  
A function $\phi(\cdot)$ is said to have the K{\L} property
at $\bar{\mathbf x}\in{\rm dom}\,\partial\, \phi$ if there exists $\varsigma\in(0,+\infty]$,
a neighborhood $U$ of $\bar{\mathbf x}$ and a concave function $\Upsilon:[0,\varsigma)\to\mathbb{R}_{+}$
that is continuously differentiable on $(0,\varsigma)$, continuous at
$0$, $\Upsilon(0)=0$, and $\Upsilon'(t)>0$ for all $t\in(0,\eta),$ such that for all
$\mathbf x\in U\cap[\phi(\bar{\mathbf x})<\phi(\mathbf x)<\phi(\bar{\mathbf x})+\varsigma],$ we have
\begin{equation}
\label{ieq:KL}
\Upsilon'\left(\phi(\mathbf x)-\phi(\bar{\mathbf x})\right) \, \dist\left(0,\partial\phi(\mathbf x)\right)\geq1,
\end{equation}
where $\dist\left(0,\partial\phi(\mathbf x)\right)=\min\left\{ \|\mathbf z\|:\mathbf z\in\partial\phi(\mathbf x)\right\}$.
If $\phi(\mathbf x)$ has the K{\L} 
property at each point of ${\rm dom}\, \partial\phi$ then $\phi$ is a K{\L} 
  function.  
  
  When $\Upsilon(t) = c t^{1-\mathbf a}$, where $c$ is a constant, we call $\mathbf a$ the {K\L} coefficient. 
\end{definition}
 Many non-convex non-smooth functions in practical applications belong to the class of K{\L} functions, for examples, real analytic functions, semi-algebraic functions, and locally strongly convex functions, see for example~\cite{Bochnak1998,Bolte2014}.

\vspace{-0.05in}
\section{Proofs} \label{app:proofs} 
In this appendix, we provide the proofs of all propositions and theorems of our paper. 
Before that, let us give some preliminary results.  We use $x,z$ to denote vectors in $\mathbb R^n$. 
\begin{lemma}\cite[Lemma 2.8]{Titan2020}
\label{lem1}
If the function  $x_i\mapsto \Theta(x_i,z)$ is $\rho$-strongly convex, differentiable at $z_i$, and $\nabla_{x_i} \Theta(z_i,z)=0$   then we have $$\Theta(x_i,z) \geq \frac{\rho}{2}\|x_i-z_i\|^2.$$
\end{lemma}
 We recall the notation $(x_i,z_{\ne i}) = (z_1,\ldots,z_{i-1},x_i,z_{i+1},\ldots,z_s)$. 
 Suppose we are trying to solve 
 $$\min_x \Psi(x):=\Phi(x) + \sum_{i=1}^s g_i(x_i).$$
 
\begin{proposition}\cite[Theorem 2.7]{Titan2020}
\label{prop:NSDP-titan}
 Suppose $\mathcal G^k_i: \mathbb R^{\mathbf n_i} \times \mathbb R^{\mathbf n_i} \to \mathbb R^{\mathbf n_i}$ be some extrapolation operator that satisfies $\mathcal G^k_i(x^{k}_i, x^{k-1}_i)\leq a_i^k\|x^{k}_i - x^{k-1}_i\|$. Let  $u_i(x_i,z)$ is a block surrogate function of $\Phi(x)$. We assume one of the following conditions holds:
 \begin{itemize}
 \item $x_i\mapsto u_i(x_i,z) + g_i(x_i)$ is $\rho_i$-strongly convex,
 \item the approximation error $ \Theta(x_i,z):=u_i(x_i,z)-\Phi(x_i,z_{\ne i}) $  satisfying $\Theta(x_i,z)\geq \frac{\rho_i}{2} \|x_i-z_i\|^2$
  for all $x_i$. 
\end{itemize}  
Note that $\rho_i$ may depend on $z$. Let 
  $$ x_i^{k+1}=\argmin_{x_i} u_i(x_i,x^{k,i-1}) + g_i(x_i)- \langle  \mathcal G^k_i(x^{k}_i, x^{k-1}_i),x_i\rangle.$$
  Then we have 
  \begin{equation}
  \label{NSDP-titan}
  \Psi(x^{k,i-1}) + \gamma_i^k \|x_i^k-x_i^{k-1} \|^2 \geq \Psi(x^{k,i}) + \eta_i^k \|x_i^{k+1}-x_i^{k} \|^2,
 \end{equation}
where 
$$
\begin{array}{ll}
\gamma^{k}_i=\frac{(a^k_i)^2}{2\nu \rho_i } , \qquad\eta^{k}_i = \frac{(1-\nu)\rho_i}{2},
\end{array}
$$ 
and $0<\nu<1$ is a constant. If we do not apply extrapolation, that is $a_i^k=0$, then \eqref{NSDP-titan} is satisfied with $\gamma_i^k=0$ and $\eta_i^k = \rho_i/2$. 
\end{proposition}
The following proposition is derived from \cite[Remark 3]{Hien_ICML2020} and \cite[Lemma 2.1]{Xu2013}. 
\begin{proposition}
\label{prop:Lsmooth-strong}
Suppose $x_i\mapsto \Phi(x)$ is a $L_i$-smooth convex function and $g_i(x_i)$ is convex. Define $\bar x^{k,i-1}=(x^{k+1}_1,\ldots,x^{k+1}_{i-1},\bar x^{k}_i, x^{k}_{i+1},\ldots,x^k_s)$, $\hat x_i^k=x_i^k + \alpha_i^k (x_i^k-x_i^{k-1})$ and $\bar x_i^k=x_i^k + \beta_i^k (x_i^k-x_i^{k-1})$. Let 
$x_i^{k+1}=\argmin_{x_i} \langle \nabla \Phi(\bar x^{k,i-1}),x_i\rangle + g_i(x_i)+ \frac{L_i}{2}\|x_i -\hat x_i^k\|^2. 
$
Then we have Inequality \eqref{NSDP-titan} is satisfied with 
$$ \gamma^{k}_i=\frac {L_i}{2} \big((\beta_i^k)^2 + \frac{(\gamma_i^k-\alpha_i^k)^2}{\nu} \big), \qquad\eta^{k}_i = \frac{(1-\nu)L_i}{2 }.
$$
If $\alpha_i^k=\beta_i^k$ then we have Inequality \eqref{NSDP-titan} is satisfied with
$$\gamma^{k}_i=\frac {L_i}{2} (\beta_i^k)^2  , \qquad\eta^{k}_i = \frac{L_i}{2 }. $$ 
\end{proposition} 

\vspace*{-0.1in}
\subsection{Proof of Proposition~\ref{prop:NSDP}}


(i) Suppose we are updating $x_i^k$. Let us recall that $$\mathcal L(x, y, \omega):= f(x)+\sum_{i=1}^s g_i(x_i) + h(y)+  \varphi(x, y, \omega),$$ where 
\begin{equation}
\label{hatvarphi}
\varphi(x, y, \omega)=\frac{\beta}{2}\|\mathcal A x + \mathcal By -b \|^2 + \langle \omega,\mathcal A x + \mathcal  By -b \rangle. 
\end{equation}
Denote $\mathbf u_i(x_i,z,y,\omega)=   u_i(x_i,z)+ h(y)  + \hat\varphi_i(x_i,z,y,\omega),$
where 
$$\hat\varphi_i(x_i,z,y,\omega) =  \varphi(z, y, \omega) + \langle \mathcal A_i^*\big( \omega +\beta(\mathcal A z + \mathcal By-b) \big),x_i-z_i\rangle+\frac{\kappa_i\beta}{2}\|x_i-z_i\|^2.  
$$

We see that $\hat\varphi_i(x_i,z,y,\omega)$ is a block surrogate function of $x\mapsto\varphi(x, y, \omega)$ with respect to block $x_i$, and $\mathbf u_i(x_i,z,y,\omega)$  is a block surrogate function of $ x\mapsto f(x) + h(y) + \varphi(x, y, \omega)$ with respect to block $x_i$. The update in~\eqref{eq:xi-updatenew} can be rewritten as follows.
\begin{equation}
\label{eq:xi_update}
x_i^{k+1}=\argmin_{x_i} \mathbf u_i(x_i,x^{k,i-1},y^k,\omega^k) + g_i(x_i) - \langle  \mathcal G^k_i(x^{k}_i, x^{k-1}_i),x_i\rangle,
\end{equation}
where 
\begin{equation}
\label{eq:chooseG}
\begin{split}
\mathcal G^k_i(x^{k}_i, x^{k-1}_i)&= \beta \mathcal A_i^* \mathcal A \big(x^{k,i-1} - \bar x^{k,i-1})\big) + \kappa_i \beta \zeta_i^k (x_i^k - x_i^{k-1}).
\end{split}
\end{equation} 
 The block approximation error function between  $\mathbf u_i(x_i,z,y,\omega)$ and $x\mapsto f(x) + h(y) + \varphi(x, y, \omega)$ is defined as 
\begin{equation}
\label{eq:errore}
\begin{split}
&\mathbf e_i(x_i,z,y,\omega)=\mathbf u_i(x_i,z,y,\omega)-\big(f(x_i,z_{\ne i}) + h(y) + \varphi((x_i,z_{\ne i}), y, \omega)\big)\\
&=u_i(x_i,z) - f(x_i,z_{\ne i}) + \hat\varphi_i(x_i,z,y,\omega)  - \varphi((x_i,z_{\ne i}), y, \omega)\\
&\geq \theta_i(x_i,z,y,\omega):= \\
&\,\varphi(z, y, \omega) - \varphi((x_i,z_{\ne i}), y, \omega) + \langle \mathcal A_i^*\big( \omega+\beta(\mathcal A z + \mathcal By-b) \big),x_i-z_i\rangle+\frac{\kappa_i\beta}{2}\|x_i-z_i\|^2.
\end{split}
\end{equation}
We have $ \nabla_{x_i}\theta_i(x_i,z,y,\omega)=\kappa_i\beta (x_i - z_i) +\nabla_{x_i} \varphi(z, y, \omega) - \nabla_{x_i} \varphi((x_i,z_{\ne i}), y, \omega) $. So  $ \nabla_{x_i}\theta_i(z_i,z,y,w)=0$. On the other hand, note that $ x_i \mapsto \varphi((x_i,z_{\ne i}), y^k, \omega^k)$ is $\beta \| \mathcal A_i^* \mathcal A_i\|$ - smooth. So, $x_i\mapsto \theta_i(x_i,z,y,\omega)$ is a $ \beta(\kappa_i - \|\mathcal A_i^* \mathcal A_i\|) $ - strongly convex function. From Lemma~\ref{lem1} we have  $\theta_i(x_i,z,y,w)\geq \frac{ \beta(\kappa_i - \|\mathcal A_i^* \mathcal A_i\|) }{2} \|x_i-z_i\|^2$. The result follows from~\eqref{eq:xi_update}, \eqref{eq:errore} and Proposition~\eqref{prop:NSDP-titan}.

(ii) When $x_i\mapsto u_i(x_i,z)+g_i(x_i)$ is convex and we apply the update as in 
\eqref{eq:xi-updatenew}, it follows from Proposition~\ref{prop:Lsmooth-strong} (see also \cite[Remark 4.1]{Titan2020}) that 
\begin{equation}
\label{eq:strongcase}
\begin{split}
& u_i(x_i^{k},x^{k,i-1}) + g_i(x_i^k)+\varphi(x^{k,i-1},y^k,\omega^k)+ \frac{\beta\|\mathcal  A_i^* \mathcal A_i\|}{2} (\zeta_i^k)^2\|x_i^{k}-x^{k-1}_i\|^2 \\
&\quad \geq  u_i(x_i^{k+1},x^{k,i-1}) + g_i(x_i^{k+1})+ \varphi(x^{k,i},y^k,\omega^k) + \frac{\beta\|\mathcal A_i^* \mathcal A_i\|}{2}\|x_i^{k+1}-x^k_i\|^2. 
 \end{split}
\end{equation}
On the other hand, note that $ u_i(x_i^{k},x^{k,i-1}) = f(x^{k,i-1})$ and  $u_i(x_i^{k+1},x^{k,i-1})\geq  f(x^{k,i})$.
The result follows then. 

\vspace*{-0.1in}
\subsection{Proof of Proposition~\ref{prop:NSDP-fory}}
Denote $$\hat h(y,y') = h(y') + \langle \omega, \mathcal Ax+ \mathcal B y'-b\rangle +  \langle \mathcal B^*\omega + \nabla h(y'), y-y'\rangle + \frac{L_h}{2} \|y-y'\|^2.$$
 Then we have $\hat h(y,y') +\frac{\beta}{2}\|\mathcal A x +\mathcal  By -b \|^2 $  is a surrogate function of $y\mapsto h(y) + \varphi(x,y,\omega) $. Note that the function $y\mapsto \hat h(y,y') +\frac{\beta}{2}\|\mathcal A x +\mathcal  By -b \|^2  $ is $(L_h + \beta\lambda_{\min}(\mathcal B^*\mathcal B))$-strongly convex. The result follows from Proposition~\ref{prop:NSDP-titan} (see also \cite[Section 4.2.1]{Titan2020}).
 
Suppose $h(y)$ is convex. We note that $y\mapsto \frac{\beta}{2}\|\mathcal A x + \mathcal By -b \|^2$ is also convex and plays the role of $g_i$ in Proposition~\ref{prop:Lsmooth-strong}. The result follows from Proposition~\ref{prop:Lsmooth-strong}.

\vspace*{-0.1in}  
\subsection{Proof of Proposition~\ref{prop:recursive}}
Note that 
\begin{equation}
\label{eq:L1}
\mL(x^{k+1},y^{k+1},\omega^{k+1})= \mL(x^{k+1},y^{k+1},\omega^k)
+ \frac{1}{\alpha\beta}\langle \omega^{k+1}-\omega^k, \omega^{k+1}-\omega^k \rangle
\end{equation}
From the optimality condition of \eqref{eq:y_update} we have 
$$\nabla h(\hat y^k) + L_h(y^{k+1}-\hat y^k) + \mathcal B^*\omega^k+\beta \mathcal B^*(\mathcal A x^{k+1} + \mathcal B y^{k+1}-b)=0.$$ 
Together with~\eqref{eq:omega_update} we obtain
\begin{equation}
\label{eq:nablah}
\nabla h(\hat y^k) + L_h (\Delta y^{k+1} -\delta_k \Delta y^{k} )+  \mathcal  B^*\omega^{k}+\frac{1}{\alpha}\mathcal B^*(w^{k+1}-w^k)=0.
\end{equation}
 Hence, 
\begin{equation}
\label{Bwk}
\mathcal B^*w^{k+1}=(1-\alpha)\mathcal B^* \omega^{k}- \alpha (\nabla h(\hat y^k) + L_h (\Delta y^{k+1} -\delta_k \Delta y^{k} ) ),
\end{equation}
 which implies that 
\begin{equation}
\label{eq:deltaw}
\mathcal B^*\Delta w^{k+1} = (1-\alpha)\mathcal B^* \Delta w^{k} - \alpha \Delta z^{k+1},
\end{equation}
where $ \Delta z^{k+1} = z^{k+1} - z^k$ and $z^{k+1}= \nabla h(\hat y^k) + L_h (\Delta y^{k+1} -\delta_k \Delta y^{k} ) $. We now consider 2 cases.

Case 1: $0<\alpha \leq 1$. From the convexity of $\|\cdot\|$ we have 
\begin{equation}
\label{eq:alpha1}
\|\mathcal  B^*\Delta w^{k+1}\|^2 \leq (1-\alpha) \|\mathcal B^* \Delta w^{k} \|^2 + \alpha \| \Delta z^{k+1}\|^2
\end{equation}

Case 2: $1<\alpha < 2$. We rewrite~\eqref{eq:deltaw} as
$
 \mathcal B^*\Delta w^{k+1} = - (\alpha -1) \mathcal B^* \Delta w^{k} -  \frac{\alpha}{2-\alpha}  (2-\alpha)\Delta z^{k+1}.
$
Hence
\begin{equation}
\label{eq:alpha2}
\|\mathcal B^*\Delta w^{k+1} \|^2 \leq  (\alpha -1)\|\mathcal B^* \Delta w^{k}\|^2+ \frac{\alpha^2}{(2-\alpha)} \|\Delta z^{k+1}\|^2  
\end{equation}
Combine \eqref{eq:alpha1} and \eqref{eq:alpha2} we obtain 
\begin{equation}
\label{eq:alpha3}
\| \mathcal B^*\Delta w^{k+1}\|^2 \leq  |1-\alpha|\|\mathcal B^* \Delta w^{k}\|^2+ \frac{\alpha^2}{1-|1-\alpha|} \|\Delta z^{k+1}\|^2, 
\end{equation}
which implies
\begin{equation}
\label{eq:alpha4}
(1-|1-\alpha|)\|\mathcal B^*\Delta w^{k+1} \|^2 
\leq  |1-\alpha|(\|\mathcal B^* \Delta w^{k}\|^2-  \|\mathcal B^* \Delta w^{k+1} \|^2)+ \frac{\alpha^2}{1-|1-\alpha|} \|\Delta z^{k+1}\|^2.
\end{equation}
On the other hand, when we use extrapolation for the update of $y$ we have 
\begin{equation}
\label{eq:z-recursive}
\begin{split}
\|\Delta z^{k+1}\|^2& =\|\nabla h(\hat y^k) - \nabla h(\hat y^{k-1})+ L_h (\Delta y^{k+1} -\delta_k \Delta y^{k} ) - L_h (\Delta y^{k} -\delta_{k-1} \Delta y^{k-1} ) \|^2\\
&\leq 3 L_h^2 \|\hat y^k -\hat y^{k-1}\|^2 + 3 L^2_h  \|\Delta y^{k+1}\|^2  + 3 \|(1+\delta_k) L_h\Delta y^{k}-L_h\delta_{k-1} \Delta y^{k-1}\|^2  \\
&\leq 6 L_h^2 \big[  (1+\delta_k)^2\|\Delta y^{k}\|^2 + \delta_{k-1}^2 \|\Delta y^{k-1} \|^2\big]+ 3 L^2_h  \|\Delta y^{k+1}\|^2  \\
&\quad+ 6(1+\delta_k)^2 L_h^2 \|\Delta y^{k}\|^2 + 6 L_h^2 \delta_{k-1}^2\|\Delta y^{k-1}\|^2\\
&=3 L^2_h\|\Delta y^{k+1}\|^2 + 12(1+\delta_k)^2 L_h^2 \|\Delta y^{k}\|^2 + 12 L^2_h \delta_{k-1}^2\|\Delta y^{k-1}\|^2.
\end{split}
\end{equation}
If we do not use extrapolation for $y$ then we have 
\begin{equation}
\begin{split}
&\|\Delta z^{k+1}\|^2 =\|\nabla h(y^k) - \nabla h(y^{k-1})    +  L_h \Delta y^{k+1} - L_h\Delta y^{k}\|^2\\
&\leq 3 L_h^2 \|\Delta y^{k}\|^2 + 3 L^2_h  \|\Delta y^{k+1}\|^2  + 3 L_h^2 \|\Delta y^{k}\|^2= 6 L_h^2 \|\Delta y^{k}\|^2+ 3 L^2_h  \|\Delta y^{k+1}\|^2.
\end{split}
\end{equation}
Furthermore, note that $\sigma_{\mathcal B}\|\Delta w^{k+1}\|^2 \leq \|\mathcal B^* \Delta w^{k+1}\|^2$.
Therefore, it follows from \eqref{eq:alpha4} that
\begin{equation}
\label{eq:alpha5}
\begin{split}
&\|\Delta w^{k+1}\|^2 \leq \frac{|1-\alpha|}{\sigma_{\mathcal B}(1-|1-\alpha|)} (\|\mathcal B^* \Delta w^{k}\|^2-  \|\mathcal B^* \Delta w^{k+1} \|^2) \\
&+ \frac{\alpha^2 3 L^2_h}{\sigma_{\mathcal B}(1-|1-\alpha|)^2}( \|\Delta y^{k+1}\|^2 + \bar\delta_k   \|\Delta y^{k}\|^2 + 4\delta_{k-1}^2\|\Delta y^{k-1}\|^2).
\end{split}
\end{equation} 
The result is obtained from~\eqref{eq:alpha5}, \eqref{eq:L1} and Proposition~\ref{prop:NSDP}.

\vspace*{-0.1in}
\subsection{Proof of Proposition~\ref{prop:subsequential}}
(i) From Inequality~\eqref{eq:recursive2}  and the conditions in~\eqref{requirement},
\begin{equation}
\label{eq:recursive3-1}
\begin{split}
&\mL^{k+1} + \mu\|\Delta y^{k+1}\|^2  +\sum_{i=1}^s\eta_i \| \Delta x^{k+1}_i \|^2  + \frac{\alpha_1}{\beta }  \|\mathcal  B^* \Delta w^{k+1}\|^2 \\
&\leq \mL^{k}+  C_1\mu\|\Delta y^{k}\|^2  + C_2\mu\|\Delta y^{k-1}\|^2+  C_x\sum_{i=1}^s\eta_i \|\Delta x^{k}_i\|^2  +  \frac{\alpha_1}{ \beta }  \|\mathcal B^* \Delta w^{k}\|^2.
\end{split}
\end{equation}
By summing  from $k=1$ to $K$ Inequality~\eqref{eq:recursive3-1} and noting that $C_1+C_2=C_y$ we obtain~\eqref{eq:upperbounded}.

(ii) Let us prove $\{\Delta y^k\}$ and $\{\Delta x_i^k \}$ converge to 0.  
Let us first prove the second situation, that is we use extrapolation for the update of $y$ and Inequality~\eqref{requirement2} is satisfied. 
 From~\eqref{Bwk} we have
$
\alpha \mathcal B^*w^{k+1}=-(1-\alpha) \mathcal B^* \Delta \omega^{k+1}- \alpha z^{k+1} ,
$
 where $z^{k+1}= \nabla h(\hat y^k) + L_h (\Delta y^{k+1} -\delta_k \Delta y^{k} ) $. 
 Using the same technique that derives Inequality~\eqref{eq:alpha3}, we obtain the following 
 \begin{equation}
 \label{eq:alpha55}
\alpha \sigma_{\mathcal B}\|w^{k+1}\|^2 \leq \alpha \|\mathcal B^*w^{k+1}\|^2 \leq 
|1-\alpha| \|\mathcal B^* \Delta \omega^{k+1}\|^2 + \frac{\alpha^2}{1-|1-\alpha|} \|z^{k+1}\|^2.  
 \end{equation}
On the other hand, we have 
  \begin{align*}
 \mL^{k} &=F(x^{k}) +h(y^{k}) +\frac{\beta}{2}\|\mathcal Ax^{k}+\mathcal By^{k}-b+\frac{\omega^{k}}{\beta}\|^2 -\frac{1}{2\beta} \| \omega^{k}\|^2\geq F(x^{k}) +h(y^{k}) -\frac{1}{2\beta} \| \omega^{k}\|^2. 
 \end{align*}
 Together with~\eqref{eq:alpha55} and 
\begin{align*}
 \|z^{k}\|^2& =
 \|\nabla h(\hat y^{k-1})- \nabla h(y^{k})+\nabla h(y^{k})   +  L_h (\Delta y^{k} -\delta_{k-1}\Delta y^{k-1} ) \|^2
  \\  &\leq   4 \|\nabla h(\hat y^{k-1})- \nabla h(y^{k}) \|^2 + 4 \| \nabla h(y^{k}) \|^2 + 4L_h^2\| \Delta y^{k}\|^2 + 4 L_h^2\delta_{k-1}^2 \|\Delta y^{k-1}\|^2
  \\  &\leq    12L_h^2 \|\Delta y^{k}\|^2 + 12\delta_{k-1}^2 \|\Delta y^{k-1}\|^2  + 4 \| \nabla h(y^{k})\|^2. 
\end{align*}
 we obtain
  \begin{equation}
  \label{temp2}
  \begin{split}
 \mL^{k} &\geq F(x^{k}) +h(y^{k}) -  \frac{1}{2\alpha\beta\sigma_{\mathcal B}}\big( |1-\alpha| \|B^* \Delta \omega^{k}\|^2 + \frac{\alpha^2}{1-|1-\alpha|} \|z^{k}\|^2 \big)
 \\ 
 &\geq F(x^{k}) +h(y^{k}) -  \frac{|1-\alpha|}{2\alpha\beta\sigma_{\mathcal B}}   \|B^* \Delta \omega^{k}\|^2\\
 &\qquad  -   \frac{\alpha}{2\beta\sigma_{\mathcal B}(1-|1-\alpha|)} \big(12L_h^2 \|\Delta y^{k}\|^2 + 12\delta_{k-1}^2 \|\Delta y^{k-1}\|^2  + 4 \| \nabla h(y^{k})\|^2\big)
 \end{split}
 \end{equation}
Since $h(y)$ is $L_h$-smooth, for all $y\in \mathbb R^q$ and $\alpha_L>0$ we have, (see \cite{Nesterov2004}) 
$$h(y-\alpha_L \nabla f(y)) \leq h(y) - \alpha_L(1-\frac{L_h \alpha_L}{2}) \|\nabla h(y)\|^2. 
$$ 
Let us choose $\alpha_L$ such that $\alpha_L(1-\frac{L_h \alpha_L}{2})=\frac{4\alpha}{2\beta\sigma_{\mathcal B}(1-|1-\alpha|)}$. Note that this equation always has a positive solution when $\beta \geq \frac{4L_h \alpha}{\sigma_{\mathcal{ B}} (1-|1-\alpha| )} $. 
Then we have 
  $$ h(y^{k}) - \frac{4\alpha}{2\beta\sigma_{\mathcal B}(1-|1-\alpha|)} \|\nabla h(y^{k})\|^2 \geq h(y^k-\alpha_L \nabla f(y^k)).
  $$
Together with~\eqref{temp2} we get
\begin{equation}
\label{temp3}
\begin{split}
\mL^{k} &\geq F(x^{k}) + h(y^k-\alpha_L \nabla f(y^k)) -  \frac{|1-\alpha|}{2\alpha\beta\sigma_{\mathcal B}}   \|B^* \Delta  \omega^{k}\|^2 \\
&\quad-   \frac{\alpha}{2\beta\sigma_{\mathcal B}(1-|1-\alpha|)} ( 12L_h^2 \|\Delta y^{k} \|^2+ 12\delta_{k-1}^2 \|\Delta y^{k-1}\|^2). 
\end{split}
\end{equation} 
 So from $\frac{\alpha_1}{\beta}\geq \frac{|1-\alpha|}{2\alpha\beta\sigma_{\mathcal B}} $, $\mu \geq \frac{\alpha 12L_h^2}{2\beta\sigma_{\mathcal B}(1-|1-\alpha|)} $, $  (1-C_1)\mu\geq \frac{\alpha 12L_h^2 12\delta_{k}^2}{2\beta\sigma_{\mathcal B}(1-|1-\alpha|)} $ we have
  \begin{equation}
  \label{temp7}
\begin{split}
&\mL^{K+1} +  \mu \|\Delta y^{K+1} \|^2  
+ \frac{\alpha_1}{ \beta }  \|B^* \Delta w^{K+1}\|^2 + (1-C_1)\mu \|\Delta y^{K}\|^2 \\
&  \geq F(x^{K+1}) + h(y^{K+1}-\alpha_L \nabla f(y^{K+1})). 
\end{split}
\end{equation}
Hence $
\mL^{K+1} +  \mu \|\Delta y^{K+1} \|^2  
+ \frac{\alpha_1}{ \beta }  \|B^* \Delta w^{K+1}\|^2 + (1-C_1)\mu \|\Delta y^{K}\|^2   
$
is lower bounded. 
%
%

Furthermore, since $\eta_i$ and $\mu$ are positive numbers we derive from  Inequality~\eqref{eq:upperbounded} that $\sum_{k=1}^\infty \|\Delta y^k \|^2<+\infty $ and $\sum_{k=1}^\infty \|\Delta x_i^k\|^2 <+\infty $. Therefore, $\{\Delta y^k\}$ and $\{\Delta x_i^k \}$ converge to 0.

Let us now consider the first situation when $\delta_k=0$ for all $k$. 

From Inequality~\eqref{eq:recursive2}  and the conditions in~\eqref{requirement} we have
\begin{equation}
\label{eq:recursive3}
\begin{split}
&\mL^{k+1} + \mu\|\Delta y^{k+1}\|^2  +\sum_{i=1}^s\eta_i \| \Delta x^{k+1}_i \|^2  + \frac{\alpha_1}{\beta }  \|B^* \Delta w^{k+1}\|^2 \\
&\leq \mL^{k}+ C_y\mu\|\Delta y^{k}\|^2  + C_x\sum_{i=1}^s\eta_i \|\Delta x^{k}_i\|^2  +  \frac{\alpha_1}{ \beta }  \|B^* \Delta w^{k}\|^2.
\end{split}
\end{equation}
By summing  Inequality~\eqref{eq:recursive3} from $k=1$ to $K$  we obtain 
\begin{equation}
\begin{split}
&\mL^{K+1} + C_y \mu \|\Delta y^{K+1} \|^2 + C_x\sum_{i=1}^s\eta_i \| \Delta x^{K+1}_i \|^2 
+ \frac{\alpha_1}{ \beta }  \|B^* \Delta w^{K+1}\|^2 \\
&\qquad\qquad+ \sum_{k=1}^{K}\big[(1-C_y)\mu \|\Delta y^{k+1}\|^2  + (1-C_x)\sum_{i=1}^s\eta_i \|\Delta x^{k+1}_i\|^2 \big] \\
&\leq \mL^1+ \frac{\alpha_1}{\beta  }  \|B^* \Delta \omega^{1}\|^2 +\sum_{i=1}^s \eta_i^0  \|\Delta x^{1}_i\|^2  +C\mu  \|\Delta y^{1}\|^2.
\end{split}
\end{equation}
Denote the value of the right side of Inequality~\eqref{eq:recursive3}  by $\hat \mL^k$. Note that $0<C_x,C_y<1$, then from~\eqref{eq:recursive3} we have the sequence $\{\hat \mL^{k}\}$ is non-increasing. It follows from \cite[Lemma 2.9]{melo2017} that $\hat\mL^k\geq \vartheta$ for all $k$, where $\vartheta$ is is the lower bound of $F(x^{k}) +h(y^{k}) $. For completeness, let us provide the proof in the following. We have
\begin{equation}
\label{temp1}
\begin{split}
\hat \mL^k&\geq \mL^k =F(x^{k}) +h(y^{k})  +\frac{\beta}{2}\|Ax^{k}+By^{k}-b\|^2 +\frac{1}{\alpha\beta}\langle \omega^k, \omega^{k}-\omega^{k-1}\rangle \\
&\geq \vartheta + \frac{1}{2\alpha\beta}(\|\omega^k\|^2-\|\omega^{k-1}\|^2+\|\Delta\omega^k\|^2)\geq \vartheta + \frac{1}{2\alpha\beta}(\|\omega^k\|^2-\|\omega^{k-1}\|^2),
\end{split}
\end{equation}
Assume that there exists $k_0$ such that $\hat \mL^k < \vartheta$ for all $k\geq k_0$. As $\hat \mL^k$ is non-increasing we have
$$\sum_{k=1}^K (\hat \mL^k - \vartheta) \leq \sum_{k=1}^{k_0} (\hat \mL^k -\vartheta) + (K-k_0) (\hat \mL^k -\vartheta) 
$$
Hence $\sum_{k=1}^\infty (\hat \mL^k - \vartheta)= -\infty$. However, from \eqref{temp1} we have
$$ \sum_{k=1}^K (\hat \mL^k - \vartheta) \geq \sum_{k=1}^K\frac{1}{2\alpha\beta}\|\omega^k\|^2 - \frac{1}{2\alpha\beta}\|\omega^{k-1}\|^2\geq \frac{1}{2\alpha\beta}(-\|\omega^{0}\|^2), 
$$
which gives a contradiction. 

Since $\hat\mL^K\geq \vartheta$ and $\eta_i$ and $\mu$ are  positive numbers we derive from  Inequality~\eqref{eq:upperbounded} that $\sum_{k=1}^\infty \|\Delta y^k \|^2<+\infty $ and $\sum_{k=1}^\infty \|\Delta x_i^k\|^2 <+\infty $. Therefore, $\{\Delta y^k\}$ and $\{\Delta x_i^k \}$ converge to 0.

Now we prove $\{\Delta \omega^k\}$ goes to 0. Since $\sum_{k=1}^\infty \|\Delta y^k \|^2<+\infty $, we derive from~\eqref{eq:z-recursive} that $\sum_{k=1}^\infty \|\Delta z^k \|^2<+\infty $.  Summing up Equality~\eqref{eq:alpha3} from $k=1$ to $K$ we have 
\begin{align*}
(1-|1-\alpha|) \sum_{k=1}^K \|\mathcal B^*\Delta \omega^k \|^2 +  \|\mathcal B^*\Delta \omega^{K+1} \|^2  
\leq \|\mathcal B^*\Delta \omega^1 \|^2 + \frac{\alpha^2}{1-|1-\alpha|} \sum_{k=1}^K \|\Delta z^{k+1} \|^2,
\end{align*} 
which implies that $\sum_{k=1}^\infty \|\mathcal B^*\Delta \omega^k \|^2 <+\infty$. Hence, $\|\mathcal B^*\Delta \omega^k \|^2\to 0$. Since $\sigma_{\mathcal B}>0$ we have $\{\Delta \omega^k\}$ goes to 0.

\vspace*{-0.05in}
\subsection{Proof of Proposition~\ref{prop:bounded-sequence}}

We remark that we use the idea in the proof of \cite[Lemma 6]{Wang2019} to prove the proposition. However, our proof is more complicated since in our framework $\alpha\in (0,2)$, the function $h$ is linearized and we use extrapolation for $y$.

  Note that as $\sigma_{\mathcal B}>0$ we have $\mathcal B$ is a surjective. Together with the assumption $b+ Im(\mathcal A) \subseteq Im(\mathcal B)$ we have there exist $\bar y^k$ such that $\mathcal Ax^k+\mathcal B\bar y^k-b =0 $. 
  
 Now we have
 \begin{equation}
 \label{temp4}
 \begin{split}
 \mL^k &=F(x^{k})+ h(y^{k}) +\frac{\beta}{2}\|\mathcal Ax^{k}+\mathcal By^{k}-b\|^2 +\langle \omega^k,\mathcal Ax^{k}+\mathcal By^{k}-b\rangle\\
 &=F(x^{k}) + h(y^{k}) +\frac{\beta}{2}\|Ax^{k}+\mathcal By^{k}-b\|^2 + \langle \mathcal B^*\omega^k,y^{k}-\bar y^k\rangle
 \end{split}
 \end{equation}
From~\eqref{eq:nablah} we have 
\begin{align*}
\langle \mathcal B^*\omega^k,y^{k}-\bar y^k\rangle
&=\big\langle\nabla  h(\hat y^k) + L_h (\Delta y^{k+1} -\delta_k \Delta y^{k} )+   \frac{1}{\alpha}\mathcal B^*(w^{k+1}-w^k), \bar y^k-y^{k}\big\rangle\\
&\geq \langle\nabla h(y^k) , \bar y^k-y^{k}\rangle-\big(\|\nabla h(y^k)-\nabla h(\hat y^k)\| + L_h \| \Delta y^{k+1}\| + L_h \delta_k \| \Delta y^{k}\|\\
&\qquad + \frac{1}{\alpha}\|\mathcal B^*\Delta\omega^{k+1}\|\big) \| \bar y^k-y^{k}\|.
\end{align*}   
Therefore, it follows from~\eqref{temp4} and $L_h$-smooth property of $h$ that
\begin{equation}
\label{temp5}
\mL^k \geq F(x^{k}) + h(\bar y^k) - \frac{L_h}{2}\|y^k-\bar y^k\|^2- \big(2L_h\delta_{k}\|\Delta y^k\| + L_h \| \Delta y^{k+1}\|  + \frac{1}{\alpha}\|\mathcal B^* \Delta\omega^{k+1}\|\big) \| \bar y^k-y^{k}\|. 
\end{equation}
On the other hand, we have 
\begin{equation}
\label{temp6}
\| \bar y^k-y^{k}\|^2 \leq \frac{1}{\lambda_{\min}(\mathcal B^*\mathcal B)} \|\mathcal B(\bar y^k-y^{k})\|^2= \frac{1}{\lambda_{\min}(\mathcal B^*\mathcal B)}\|\mathcal A x^k + \mathcal By^{k} -b\|^2
=\frac{1}{\lambda_{\min}(\mathcal B^*\mathcal B)} \big\|\frac{1}{\alpha\beta } \Delta \omega^k\big\|^2.
\end{equation}
We have proved in Proposition~\ref{prop:subsequential} that $\|\Delta \omega^k\|$,  $\|\Delta x^k\|$ and $\|\Delta y^k\|$ converge to 0. Furthermore, from Proposition~\ref{prop:subsequential} we have $\mL^k$ is upper bounded. Therefore, from \eqref{temp5}, \eqref{temp6} and \eqref{eq:upperbounded} we have $ F(x^{k}) + h(\bar y^k)$ is upper bounded. So $\{x^k\}$ is bounded. Consequently, $\mathcal Ax^k $ is bounded.  

Furthermore, we have 
$$ \|y^k\|^2\leq \frac{1}{\lambda_{\min}(\mathcal B^*\mathcal B)} \|\mathcal By^k\|^2 = \frac{1}{\lambda_{\min}(\mathcal B^*\mathcal B)}\big \|\frac{1}{\alpha\beta } \Delta \omega^k-\mathcal Ax^k -b \big\|^2.
$$
Therefore, $\{y^k\}$ is bounded, which implies that $\|\nabla h(\hat y^k)\|$ is also bounded. Finally, from ~\eqref{eq:nablah} and the assumption $\lambda_{\min}(\mathcal B\mathcal B^*)>0 $ we also have $\{\omega^k\}$ is bounded. 

\vspace*{-0.1in}
\subsection{Proof of Theorem~\ref{thrm:subsequential}}
Suppose $(x^{k_n},y^{k_n},\omega^{k_n})$ converges to $(x^*,y^*,\omega^*)$.   
Since $\Delta x_i^k $ goes to 0, we have $x_i^{k_n+1}$ and $x_i^{k_n-1}$ also converge to $x_i^*$ for all $i\in [s]$. From \eqref{eq:xi_update},  for all $x_i$, 
\begin{align}
\label{Sub4}
 \mathbf u_i(x_i^{k+1},x^{k,i-1},y^k,\omega^k)+ g_i(x_i^{k+1})   \leq \mathbf u_i(x_i,x^{k,i-1},y^k,\omega^k) + g_i(x_i) -   \langle  \mathcal G^k_i(x^{k}_i, x^{k-1}_i),x_i-x^{k+1}_i\rangle.
\end{align}
Choosing $x_i=x_i^*$ and $k=k_n-1$ in \eqref{Sub4} and noting that $ \mathbf u_i(x_i,z)$ is continuous by Assumption~\ref{assump:Lipschitz_ui} (i), we have 
$
\limsup_{n\to\infty}  \mathbf u_i(x_i^*,x^*,y^*,\omega^*) + g_i(x_i^{k_n}) \leq  \mathbf u_i(x_i^*,x^*,y^*,\omega^*)+ g_i(x_i^*). 
$
On the other hand, as $g_i(x_i)$ is lower semi-continuous. Hence, $ g_i(x_i^{k_n})$ converges to $g_i(x_i^*)$. Now we choose $k=k_n\to\infty$ in~\eqref{Sub4} for all $x_i$ we obtain 
\begin{equation}
\label{eq:blah}
\begin{split}
L_0(x^*,y^*,\omega^*) + g_i(x_i^*) &\leq \mathbf u_i(x_i,x^*,y^*,\omega^*) +g_i(x_i)\\
&= L_0(x_i,x^*_{\ne i},y^*,\omega^*) + \mathbf e_i(x_i,x^*,y^*,\omega^*) + g_i(x_i),
\end{split}
\end{equation} 
where $L_0(x,y,\omega)=f(x) + h(y) + \varphi(x, y, \omega)$ and $\mathbf e_i$ is the approximation error defined in \eqref{eq:errore}. We have 
\begin{align*}
\mathbf e_i(x_i,x^*,y^*,\omega^*) &= u_i(x_i,x^*) - f(x_i,x^*_{\ne i}) + \hat\varphi_i(x_i,x^*,y^*,\omega^*)  - \varphi((x_i,x^*_{\ne i}), y^*, \omega^*)\\
&\leq  \bar e_i(x_i,x^*) + \hat\varphi_i(x_i,x^*,y^*,\omega^*)  - \varphi((x_i,x^*_{\ne i}), y^*, \omega^*). 
\end{align*}
Note that $\bar e_i(x^*_i,x^*)=0$ by Assumption~\ref{assump:Lipschitz_ui}. From~\eqref{eq:blah} we have $x_i^*$  is a solution of 
$$\min_{x_i} L(x_i,x^*_{\ne i},y^*,\omega^*)+ \bar e_i(x_i,x^*) + \hat\varphi_i(x_i,x^*,y^*,\omega^*)  - \varphi((x_i,x^*_{\ne i}), y^*, \omega^*). 
$$
Writing the optimality condition for this problem we obtain $0 \in \partial_{x_i} \mL(x^*,y^*,\omega^*)$. Totally similarly we can prove that $0 \in \partial_{y} \mL(x^*,y^*,\omega^*) $. On the other hand, we have 
$$\Delta \omega^k= \omega^{k} - \omega^{k-1}= \alpha\beta (\mathcal A x^k + \mathcal B y^k -b)\to 0.$$ 
Hence, $\partial_\omega \mL(x^*,y^*,\omega^*) = \mathcal A x^* + \mathcal B y^* -b=0.
$

As we assume $\partial F(x)=\partial_{x_1} F(x) \times \ldots \times \partial_{x_s} F(x)$, we have 
\begin{align*}
\partial \mL(x,y,\omega) &= \partial F(x)+ \nabla \Big(h(y) +  \langle \omega,\mathcal A x +\mathcal By-b \rangle + \frac{\beta}{2} \|\mathcal A x + \mathcal By-b\|^2\Big)\\
&=\partial_{x_1} \mL(x,y,\omega) \times\ldots\times \partial_{x_s} \mL(x,y,\omega) \times \partial_{y} \mL(x,y,\omega)\times \partial_{\omega} \mL(x,y,\omega).
\end{align*}
So $ 0\in \partial \mL(x^*,y^*,\omega^*)$.
\vspace*{-0.1in}
\subsection{Proof of Theorem~\ref{thrm:global}}
Note that we assume the generated sequence of Algorithm~\ref{alg:iADMM} is bounded. The following analysis is considered in the bounded set that contains the generated sequence of Algorithm~\ref{alg:iADMM}. We first prove some preliminary results.

(A) The optimality condition of~\eqref{eq:xi_update} gives us 
\begin{equation}
\label{Sub1}
\begin{split}
&\mathcal G_i^k(x_i^k - x_i^{k-1}) - \mathcal A_i^*\big(\omega^k+\beta(\mathcal A x^{k,i-1} + \mathcal By^k-b) \big)
-\kappa_i\beta(x^{k+1}_i-x_i^k) \\
&\qquad\in \partial_{x_i} \big(u_i(x_i^{k+1},x^{k,i-1}) + g_i(x_i^{k+1})\big). 
\end{split}
\end{equation}
As \eqref{assume:partial} holds, there exists $\mathbf s_i^{k+1}\in\partial u_i(x_i^{k+1},x^{k,i-1})$ and $\mathbf t_i^{k+1}\in\partial g_i(x_i^{k+1})$ such that
\begin{equation}
\label{eq:stik}
\mathcal G_i^k(x_i^k - x_i^{k-1}) - \mathcal A_i^*\big(\omega^k+\beta(\mathcal A x^{k,i-1} + \mathcal By^k-b) \big)
-\kappa_i\beta(x^{k+1}_i-x_i^k) = \mathbf s_i^{k+1} + \mathbf t_i^{k+1}
\end{equation}
As \eqref{eq:l1} holds, there exists $\xi_i^{k+1}\in\partial_{x_i} f(x^{k+1})$ such that
\begin{equation}
\label{eq:assumpl}
\|\xi_i^{k+1} - \mathbf s_i^{k+1}\| \leq L_i\|x^{k+1} - x^{k,i-1}\|.
\end{equation}

Denote $\tau^{k+1}_i:=  \xi_i^{k+1} + \mathbf t_i^{k+1} \in \partial_{x_i} F(x^{k+1})$ (as \eqref{assume:partial} holds). 
Then, from \eqref{eq:stik} we have 
\begin{equation}
\label{tauik}
\tau^{k+1}_i= \xi_i^{k+1} +  \mathcal G_i^k(x_i^k - x_i^{k-1}) - \mathcal A_i^*\big(\omega^k+\beta(\mathcal A x^{k,i-1} + \mathcal By^k-b) \big)
-\kappa_i\beta(x^{k+1}_i-x_i^k) - \mathbf s_i^{k+1}. 
\end{equation}
On the other hand, we note that  
\begin{equation}
\label{Sub3}
\partial_{x_i} \mL(x^{k+1},y^{k+1},\omega^{k+1})= \partial_{x_i} F(x^{k+1} ) + \mathcal A_i^*\big(\omega^{k+1} + \beta(\mathcal A x^{k+1} +\mathcal  By^{k+1}-b)  \big).
\end{equation} 
Let $d_i^{k+1}:= \tau_i^{k+1}+  \mathcal A_i^*\big(\omega^{k+1} + \beta(\mathcal A x^{k+1} + \mathcal By^{k+1}-b)  \big) \in \partial_{x_i} \mL(x^{k+1},y^{k+1},\omega^{k+1}) $. 
From~\eqref{tauik}, 
\begin{equation}
\begin{split}
\|d_i^{k+1}\|& = \Big\|\xi_i^{k+1} +  \mathcal G_i^k(x_i^k - x_i^{k-1}) - \mathcal A_i^*\big(\omega^k+\beta(\mathcal A x^{k,i-1} + \mathcal By^k-b) \big)
-\kappa_i\beta(x^{k+1}_i-x_i^k) \\
&\qquad\qquad- \mathbf s_i^{k+1}+  \mathcal A_i^*\big(\omega^{k+1} + \beta(\mathcal A x^{k+1} + \mathcal By^{k+1}-b)  \big) \Big\|
\end{split}
\end{equation}
Together with \eqref{eq:assumpl} we obtain 
\begin{equation}
\begin{split}
\|d_i^{k+1}\|& \leq a^k_i\|\Delta x_i^k\| + \beta\|\mathcal A_i^* A\| \|x^{k+1}-x^{k,i-1}\| + \beta\|\mathcal A_i^*\mathcal B\|  \|\Delta y^{k+1}\|  + \|\mathcal A_i^*\| \| \Delta \omega^{k+1}\|\\
&\qquad \qquad + \kappa_i \beta \|\Delta x_i^{k+1}\| + L_i\|x^{k+1} - x^{k,i-1}\|.
\end{split}
\end{equation}

It follows from~\eqref{eq:y_update} that 
$$\mathcal B^*\omega^k + \nabla h(\hat y^k) + \beta \mathcal B^* (\mathcal A x^{k+1} +\mathcal  B y^{k+1} -b) + L_h (y^{k+1} - \hat y^k) = 0.
$$
Let  $d_y^{k+1}:=\nabla h(y^{k+1}) +\mathcal  B^*\big(\omega^{k+1} +\beta(\mathcal  A x^{k+1} + \mathcal B y^{k+1} -b )\big).$ Then 
$d_y^{k+1}\in \partial_y \mL(x^{k+1},y^{k+1},\omega^{k+1})$ and 
\begin{align*}
&\|d_y^{k+1}\| =\|\nabla h(y^{k+1}) - \nabla h(\hat y^{k}) +\mathcal  B^*(\omega^{k+1} - \omega^k) -  L_h (y^{k+1} - \hat y^k)\|\\
& \leq 2L_h \| y^{k+1} - \hat y^{k}\| + \|\mathcal B^*\| \| \Delta \omega^{k+1} \|  \leq 2 L_h (\| \Delta y^{k+1}\| + \delta_k \| \Delta y^{k}\|) + \|\mathcal B^*\| \| \Delta \omega^{k+1} \|.
\end{align*}

Let $d_\omega^{k+1}:=\mathcal A x^{k+1} + \mathcal B^{k+1} -b$. We have  $d_\omega^{k+1}\in \partial_\omega \mL(x^{k+1},y^{k+1},\omega^{k+1})$ and  $$d_\omega^{k+1}=(\omega^{k+1} - \omega^k)/(\alpha \beta) = \Delta \omega^{k+1}/(\alpha \beta).$$

(B) Let us now prove $F(x^{k_n})$ converges to $F(x^*)$. This implies $\mL(x^{k_n},y^{k_n},\omega^{k_n})$ converges to $\mL(x^*,y^*,\omega^*)$ since $\mL$ is differentiable in $y$ and $\omega$. 
We have
\begin{align*}
F(x^{k_n})= f(x^{k_n})+\sum_{i=1}^s g_i(x_i^{k_n}) =u_s(x_s^{k_n},x^{k_n})  +\sum_{i=1}^s g_i(x_i^{k_n}).
\end{align*}
So $F(x^{k_n})$ converges to $u_s(x_i^*,x^*) +\sum_{i=1}^s g_i(x_i^*)=F(x^*)$.

We now proceed to prove the global convergence. Denote $\mathbf z= (x,y,\omega)$, $\tilde{\mathbf{z}}= (\tilde x, \tilde y, \tilde \omega)$, and $\mathbf z^k= (x^k,y^k,\omega^k)$. 
We consider the following auxiliary function
\begin{equation*}
\bar\mL(\mathbf z, \tilde{\mathbf{z}})=\mL(x,y,\omega) + \sum_{i=1}^s \frac{\eta_i + C_x \eta_i}{2}\|x_i - \tilde x_i \|^2 + \frac{(1+C_y) \mu}{2} \|y-\tilde y\|^2 + \frac{\alpha_1}{\beta} \|B^* (\omega - \tilde \omega)\|^2.
\end{equation*}
The auxiliary sequence $ \bar\mL (\mathbf z^k, \mathbf z^{k-1})$ has the following properties.
\begin{enumerate}
\item \textbf{Sufficient decreasing property}. From~\eqref{eq:recursive3} we have
\begin{align*}
&\bar\mL (\mathbf z^{k+1}, \mathbf z^{k}) + \sum_{i=1}^s \frac{\eta_i- C_x  \eta_i}{2}\big(  \|x_i^{k+1} -x_i^k \|^2 
+  \|x_i^{k} -x_i^{k-1} \|^2\big) \\
&\qquad+ \frac{(1-C_y)\mu}{2} \big(  \|y^{k+1} -y^k \|^2 +  \|y^{k} -y^{k-1} \|^2\big)\leq \bar\mL (\mathbf z^k, \mathbf z^{k-1}). 
\end{align*}
\item \textbf{Boundedness of subgradient}. In the proof (A) above, we have proved that 
$$\|d^{k+1}\| \leq a_1 (\|x^{k+1}-x^k\|+\|x^k-x^{k-1}\| + \|y^{k+1}-y^k\| + \|\omega^{k+1}-\omega^k\|)$$ for some constant $a_1$ and $d^{k+1} \in \partial \mL(\mathbf z^{k+1})$. On the other hand, as we use $\alpha=1$, from~\eqref{eq:deltaw} we obtain 
 \begin{equation}
 \label{Sub5}
\begin{split}
&\sqrt{\sigma_{\mathcal B}}\|\omega^{k+1}-\omega^k\|\leq \|B^*(\omega^{k+1}-\omega^k)\| = \|\Delta z^{k+1}\|\\
& =\|\nabla h(y^{k}) - \nabla h(y^{k-1}) + L_h(\Delta y^{k+1} - \Delta y^k) \|\leq 2L_h\|y^{k}-y^{k-1}\| + L_h\|y^{k+1}-y^{k}\|.  
\end{split}
\end{equation} 
Hence, 
$$\|d^{k+1}\| \leq a_2 (\|x^{k+1}-x^k\|+\|x^k-x^{k-1}\| + \|y^{k+1}-y^k\| + \|y^{k}-y^{k-1}\|)$$ for some constant $a_2$. 
Note that 
$$\partial \bar\mL(\mathbf z, \tilde{\mathbf{z}})=\partial \mL(\mathbf z, \tilde{\mathbf{z}}) + \partial \Big(\sum_{i=1}^s \frac{\eta_i + C_x \eta_i}{2}\|x_i - \tilde x_i \|^2 + \frac{(1+C_y) \mu}{2} \|y-\tilde y\|^2 + \frac{\alpha_1}{\beta} \|B^* (\omega - \tilde \omega)\|^2 \Big).$$ 
Hence, it is not difficult to show that 
$$\|\mathbf d^{k+1}\| \leq a_3 (\|x^{k+1}-x^k\|+\|x^k-x^{k-1}\| + \|y^{k+1}-y^k\| + \|y^{k}-y^{k-1}\|)$$ for some constant $a_3$ and $\mathbf d^{k+1} \in \partial \bar\mL(\mathbf z^{k+1},\mathbf z^{k})$.
\item \textbf{KL property}. Since $F(x) + h(y)$ has KL property, then $ \bar\mL(\mathbf z, \tilde{\mathbf{z}})$ also has K{\L} property. 
\item \textbf{A continuity condition}. Suppose $\mathbf z^{k_n} $ converges to $(x^*,y^*,\omega^*)$. In the proof (B) above, we have proved that $\mL (\mathbf z^{k_n})$ converges to $\mL(x^*,y^*,\omega^*)$. Furthermore, from Proposition~\ref{prop:subsequential} we proved that $\|\mathbf z^{k+1}-\mathbf z^{k} \|$ goes to 0. Hence we have $\mathbf z^{k_n-1} $ converges to $(x^*,y^*,\omega^*)$.  So, $\bar\mL (\mathbf z^{k+1}, \mathbf z^{k})$  converges to $\bar\mL (\mathbf z^*, \mathbf z^*)$. 
\end{enumerate}
Using the same technique as in \cite[Theorem 1]{Bolte2014}, see also \cite{Hien_ICML2020,Ochs2019}, we can prove that 
$$\sum_{k=1}^\infty \big(\|x^{k+1}-x^k\|+\|x^k-x^{k-1}\| + \|y^{k+1}-y^k\| + \|y^{k}-y^{k-1}\| \big)<\infty.
$$
which implies $\{(x^k,y^k)\} $ converges to $(x^*,y^*)$.  From~\eqref{Sub5} we obtain 
$$\sum_{k=1}^\infty \|\omega^{k+1}-\omega^k\| \leq \sum_{k=1}^\infty \big( \|y^{k+1}-y^k\| + \|y^{k}-y^{k-1}\| \big)<\infty.
$$
Hence, $\{\omega^k\}$ also converges to $\omega^*$.

\section{Additional Experiment for different values of $\alpha$}
\label{sec:addexp}

\begin{figure*}[h!]
\begin{center}
\begin{tabular}{cc}
\includegraphics[width=0.5\textwidth]{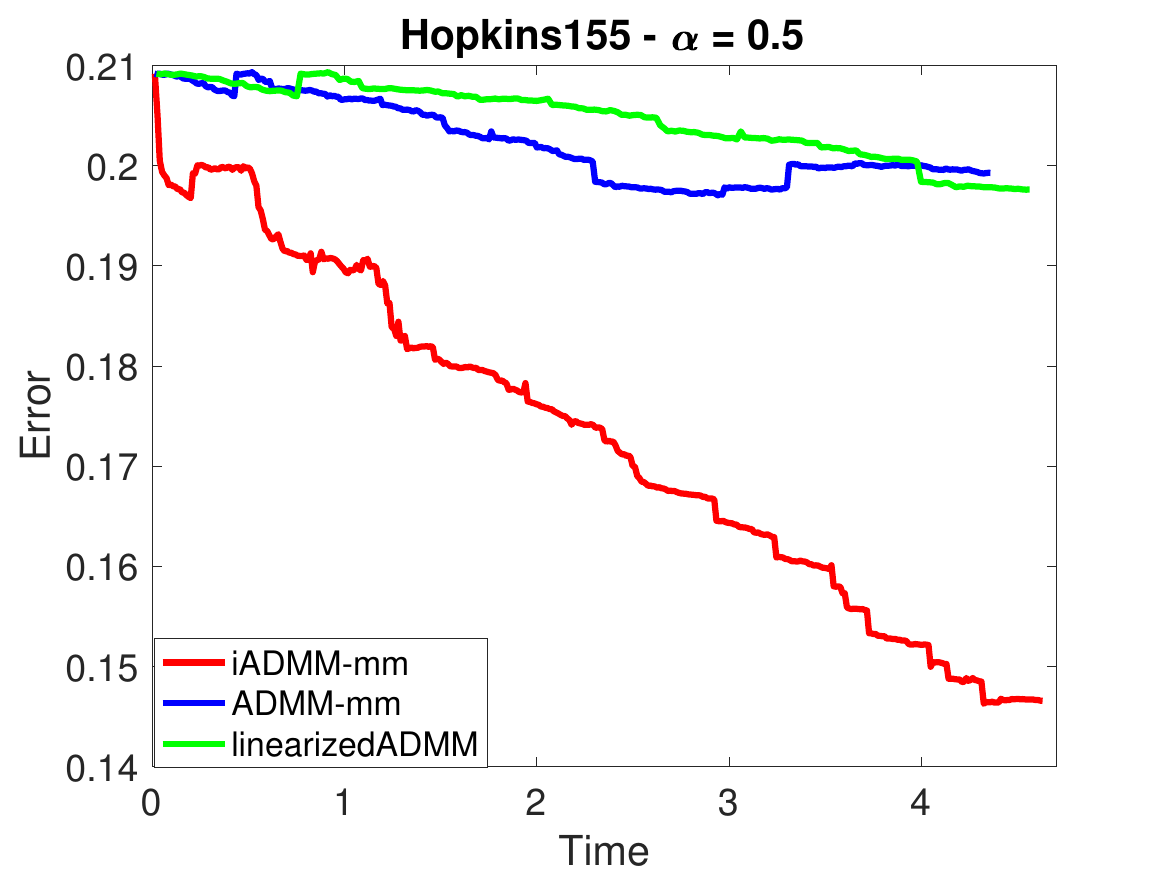}  & 
\includegraphics[width=0.5\textwidth]{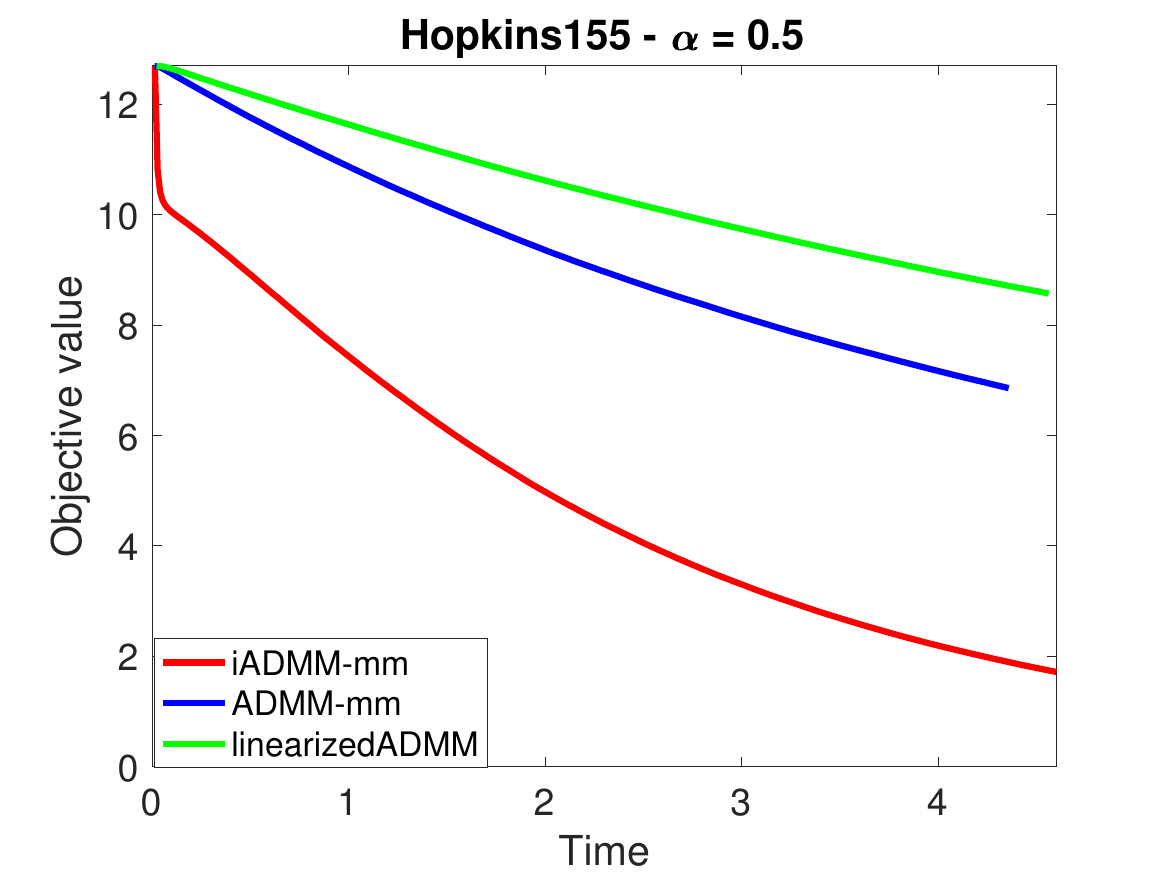} \\
\includegraphics[width=0.5\textwidth]{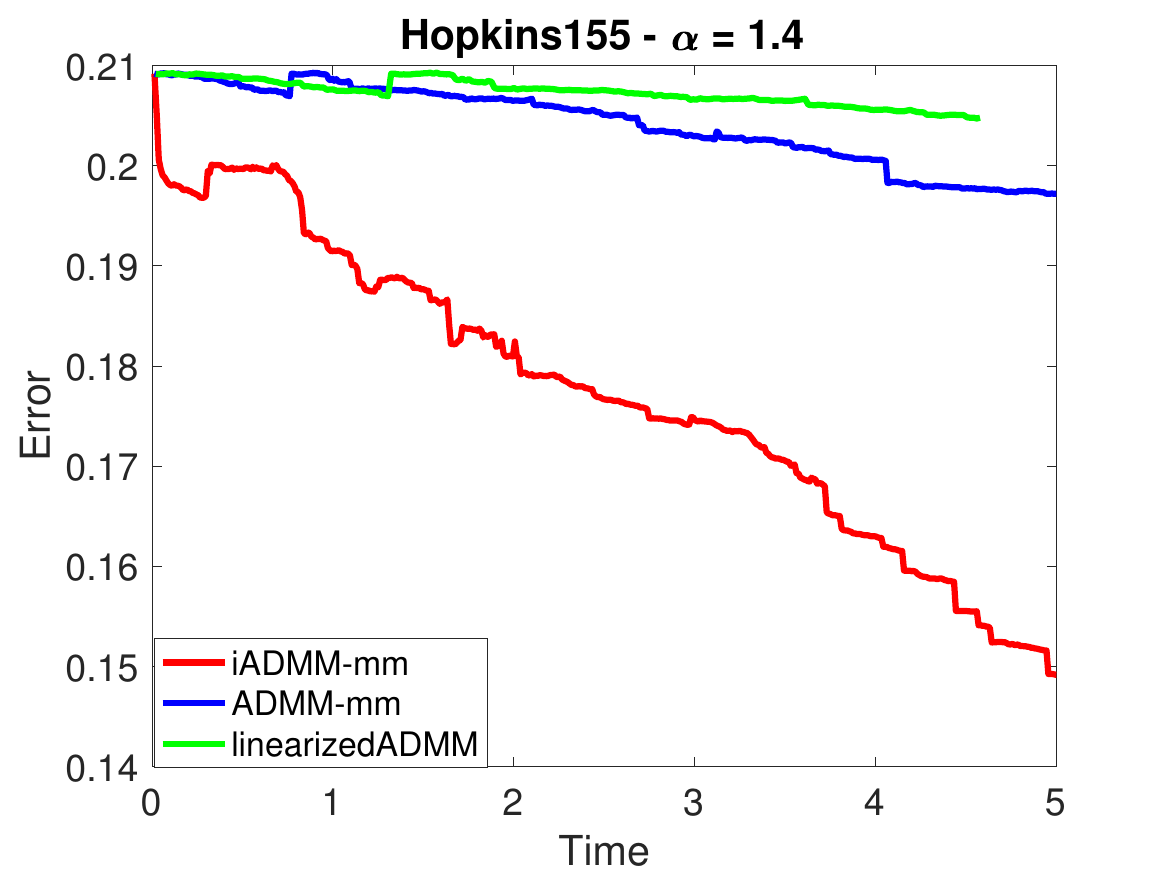}  & 
\includegraphics[width=0.5\textwidth]{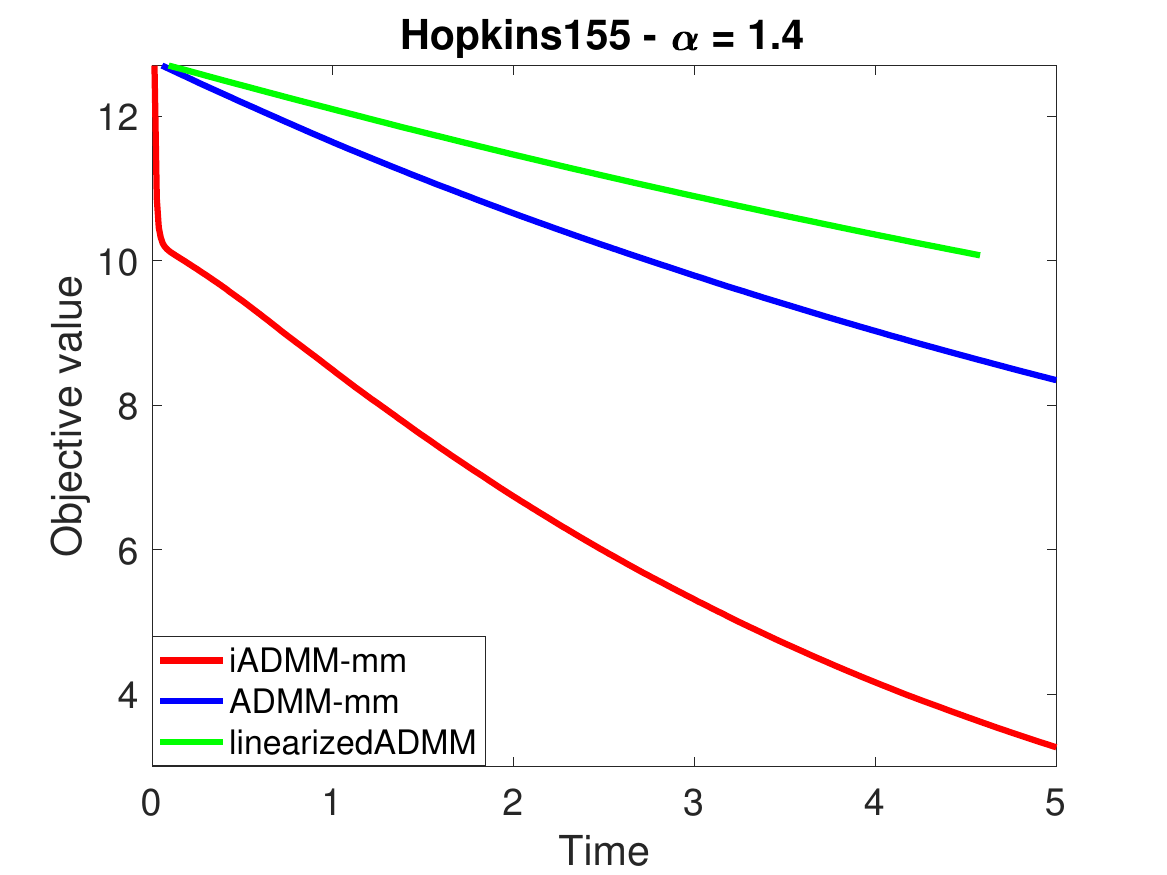} \\
\includegraphics[width=0.5\textwidth]{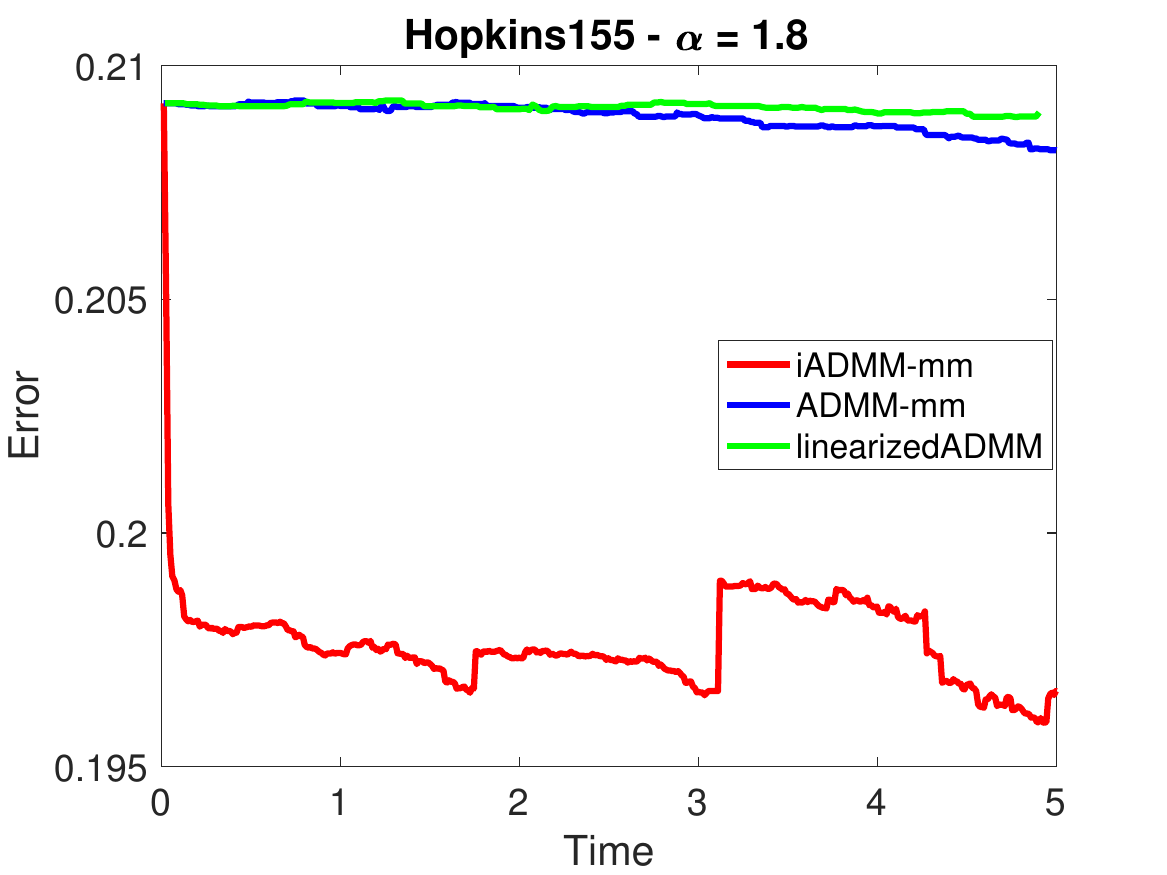}  & 
\includegraphics[width=0.5\textwidth]{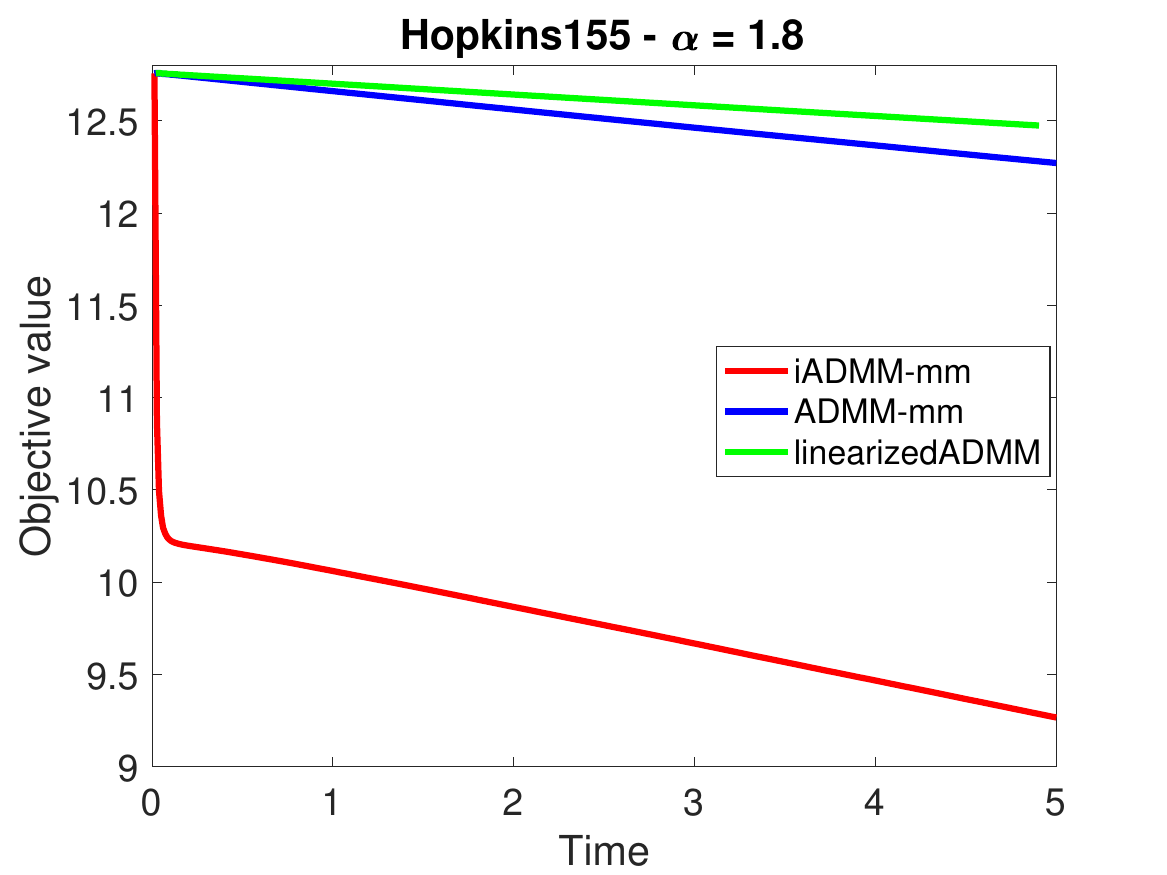}
\end{tabular}
\caption{Evolution of the average value of the segmentation error rate and the objective function value with respect to time on Hopkins155. 
\label{fig:Hopkins155}} 
\end{center}
\end{figure*}

\begin{figure*}[h!]
\begin{center}
\begin{tabular}{cc}
\includegraphics[width=0.5\textwidth]{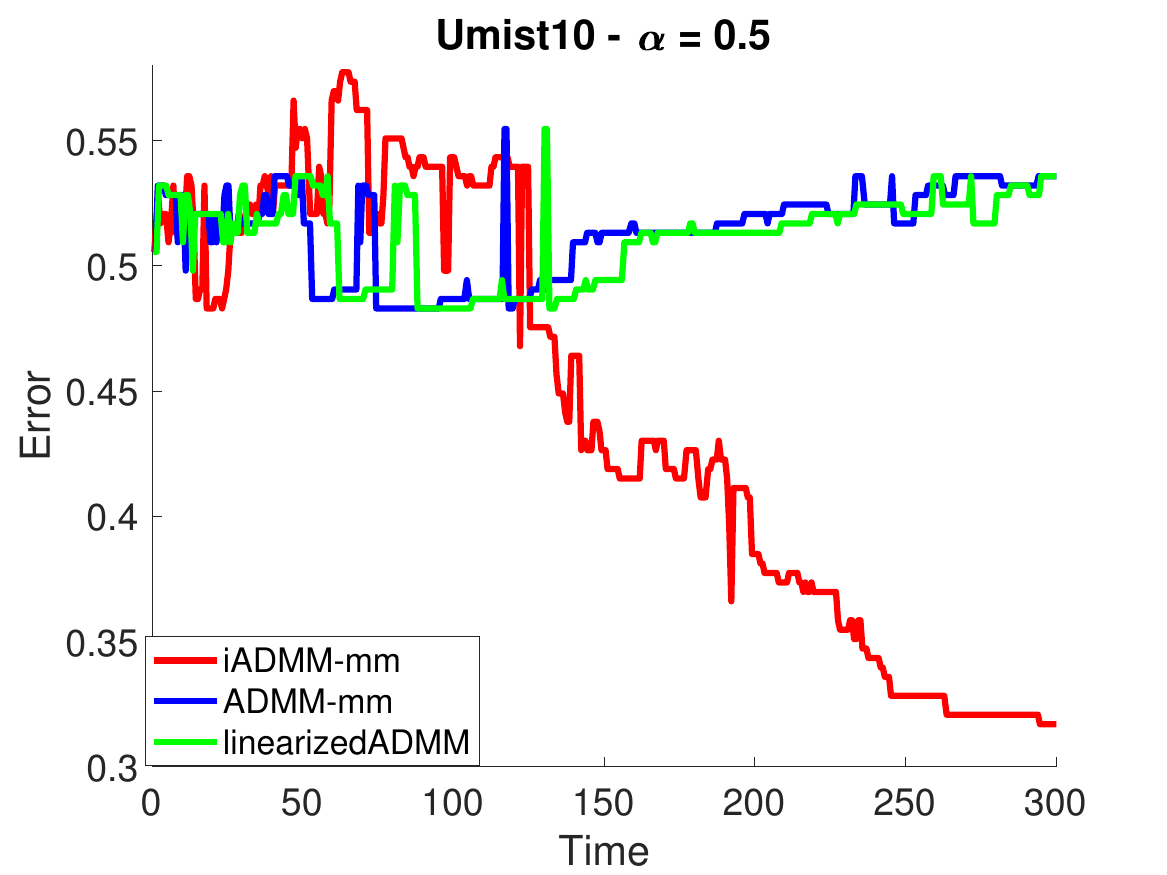}  & 
\includegraphics[width=0.5\textwidth]{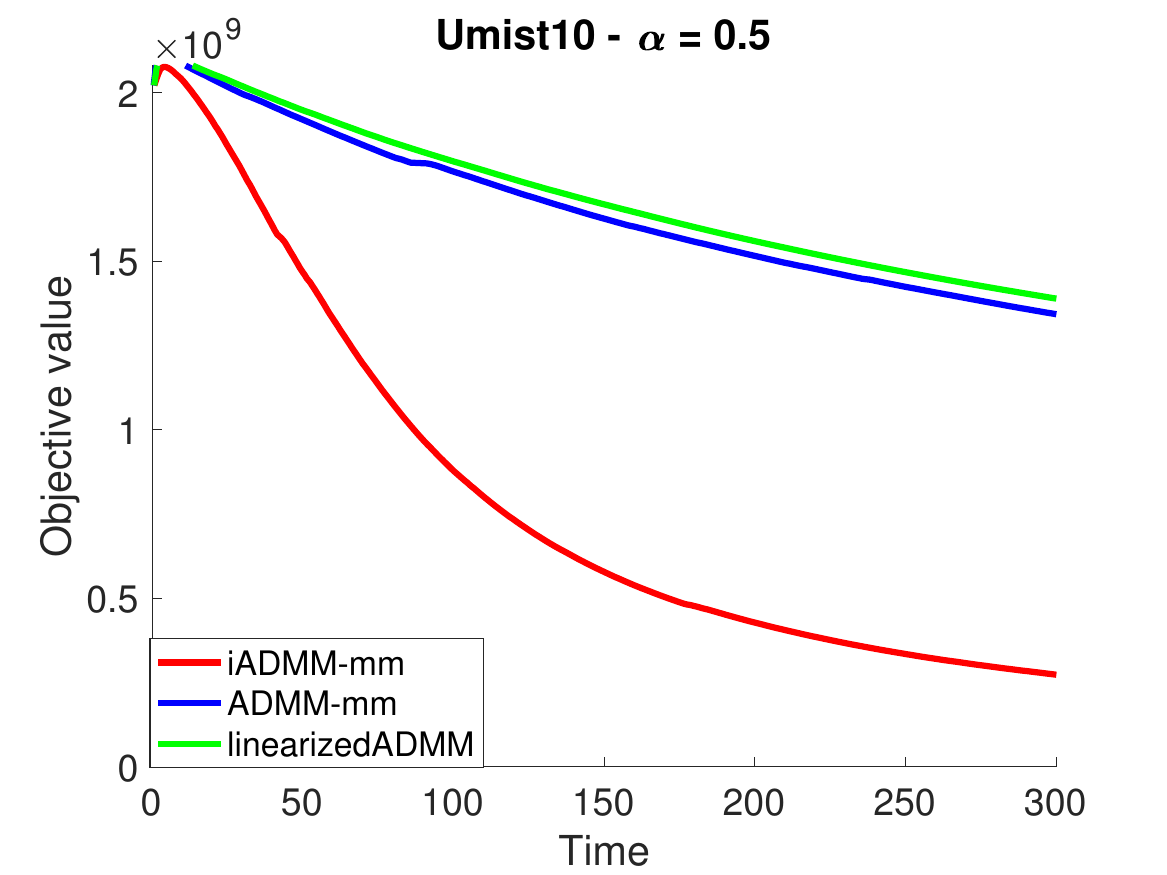} \\
\includegraphics[width=0.5\textwidth]{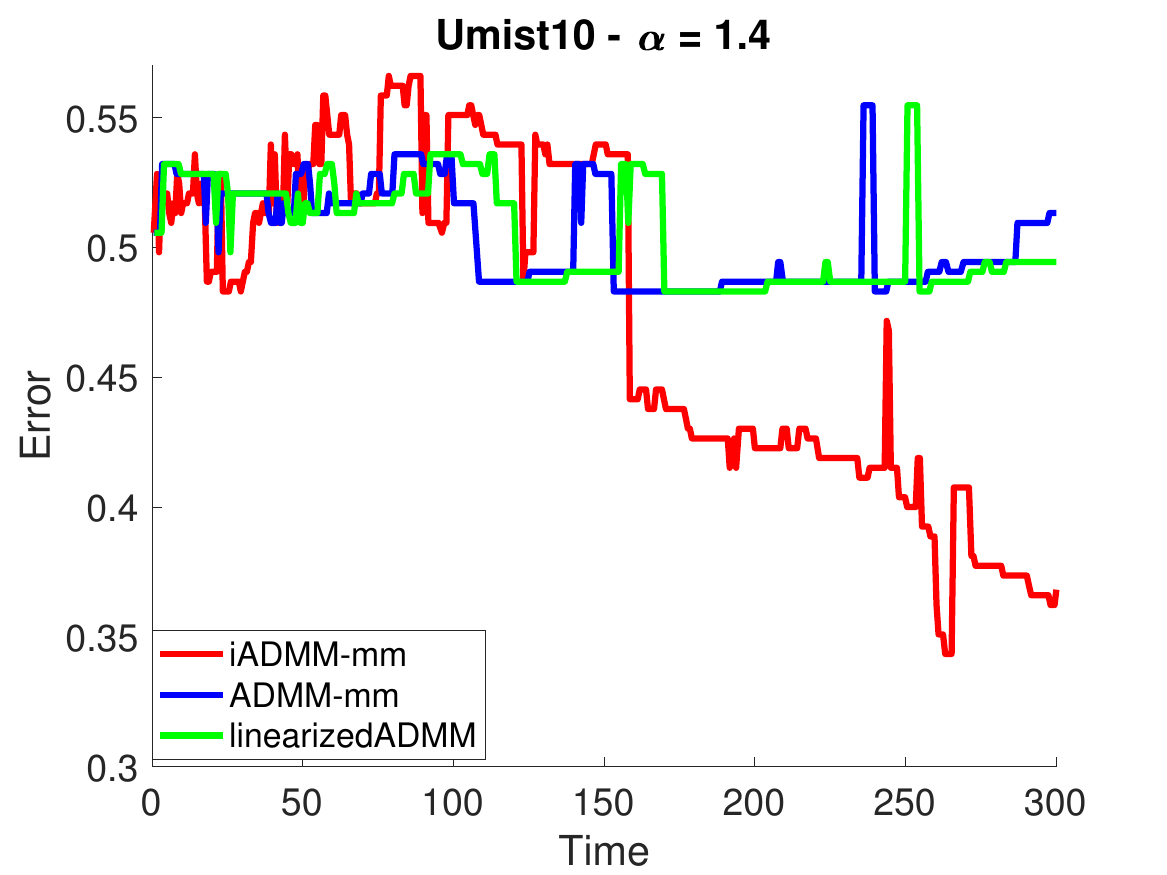}  & 
\includegraphics[width=0.5\textwidth]{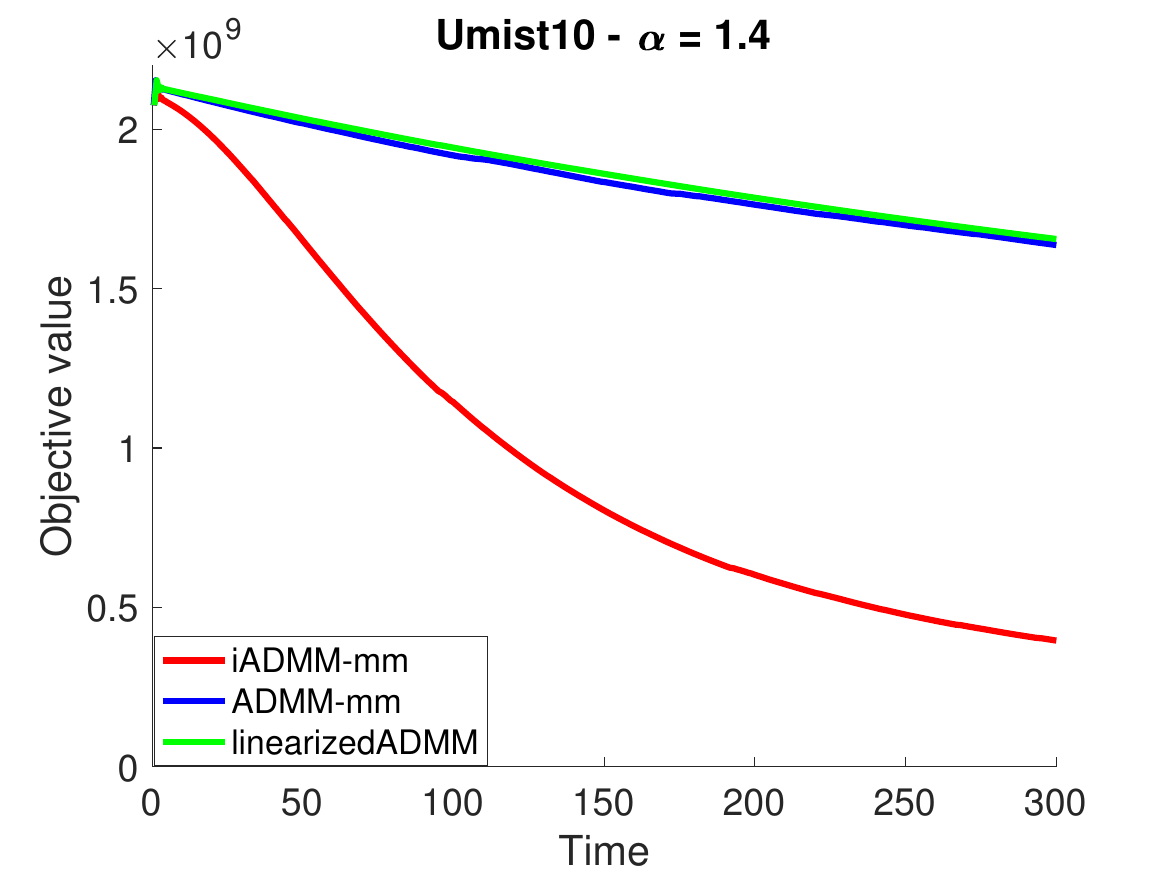} \\
\includegraphics[width=0.5\textwidth]{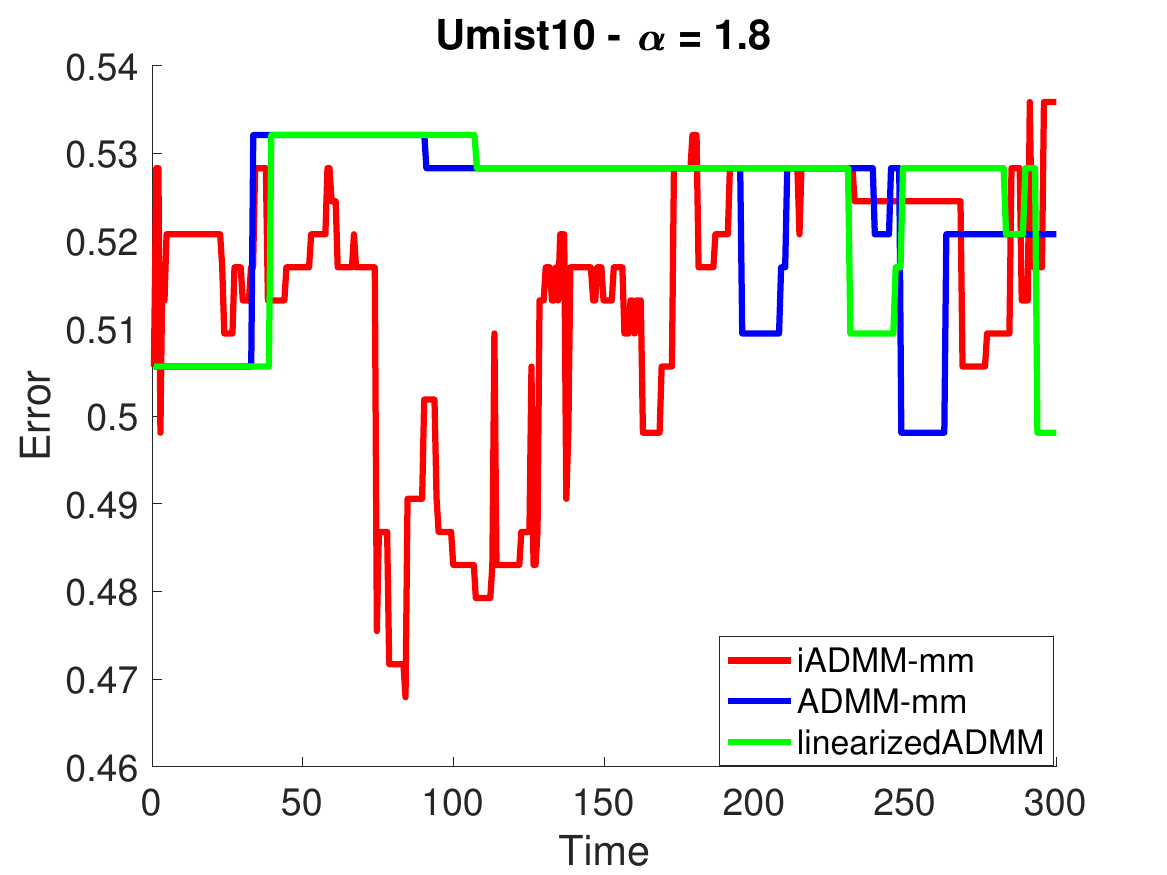}  & 
\includegraphics[width=0.5\textwidth]{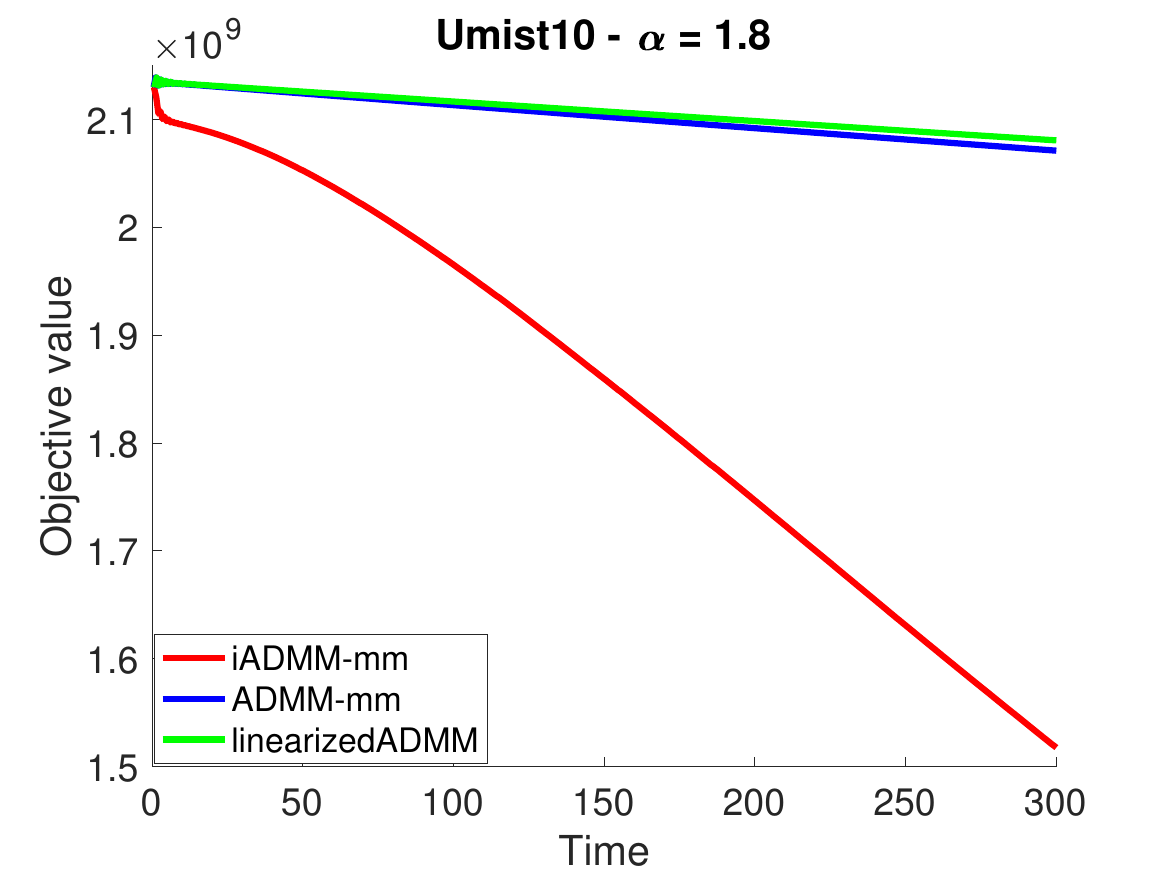}
\end{tabular}
\caption{Evolution of the segmentation error rate and the objective function value with respect to time on Umist10. 
\label{fig:Umist10}} 
\end{center}
\end{figure*}

\begin{figure*}[h!]
\begin{center}
\begin{tabular}{cc}
\includegraphics[width=0.5\textwidth]{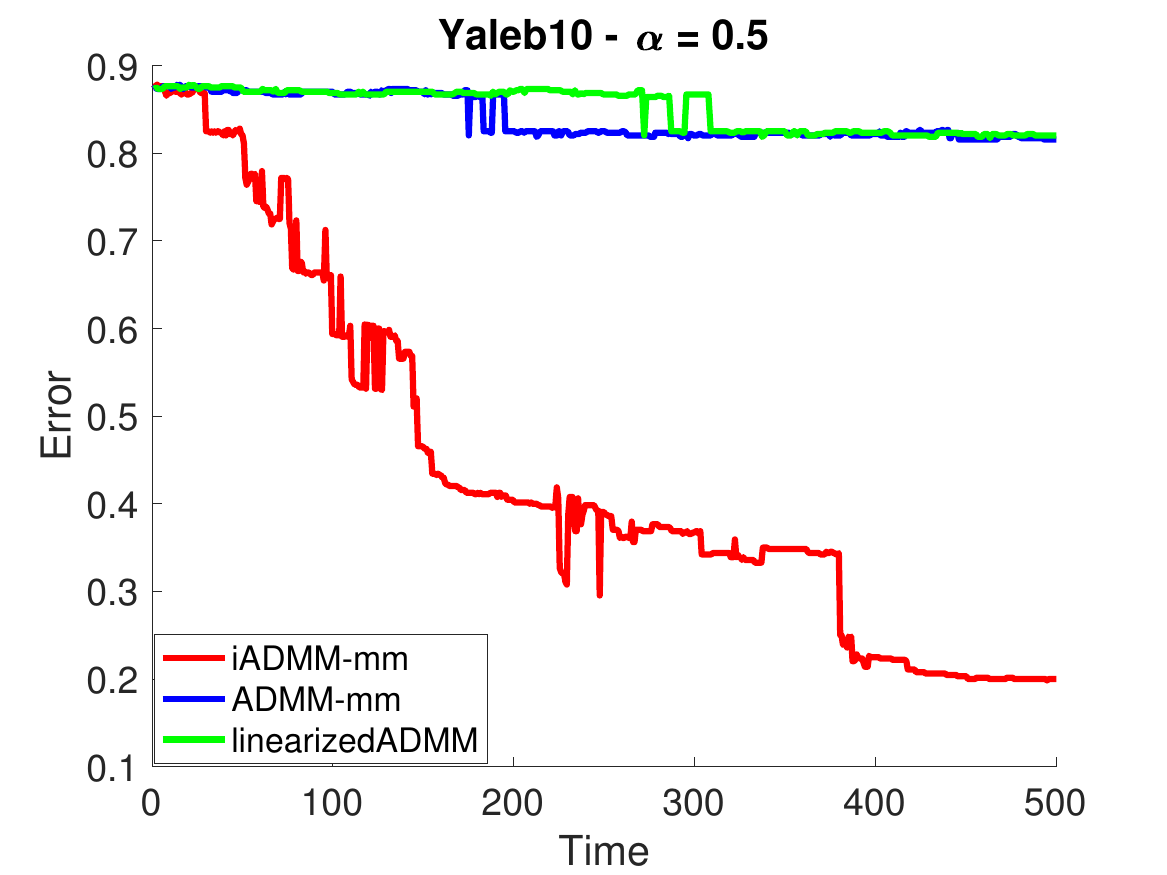}  & 
\includegraphics[width=0.5\textwidth]{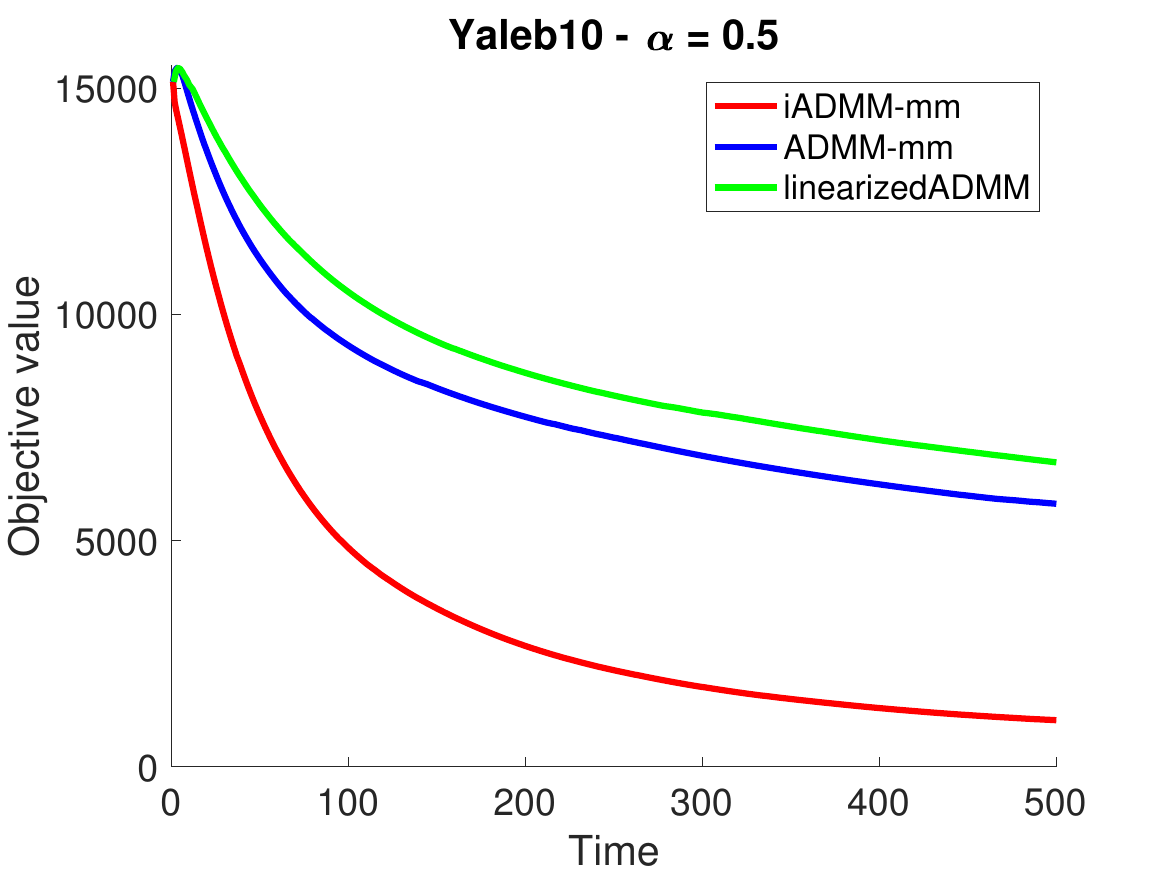} \\
\includegraphics[width=0.5\textwidth]{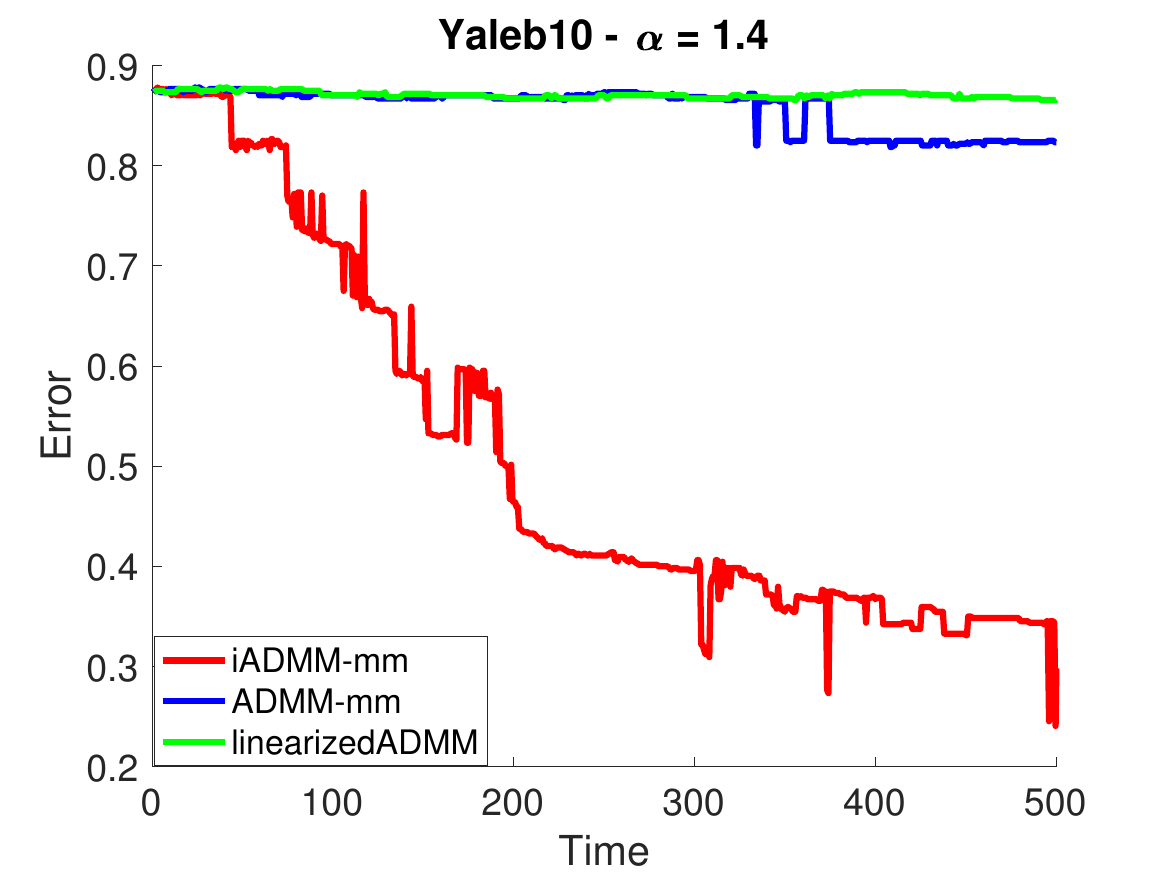}  & 
\includegraphics[width=0.5\textwidth]{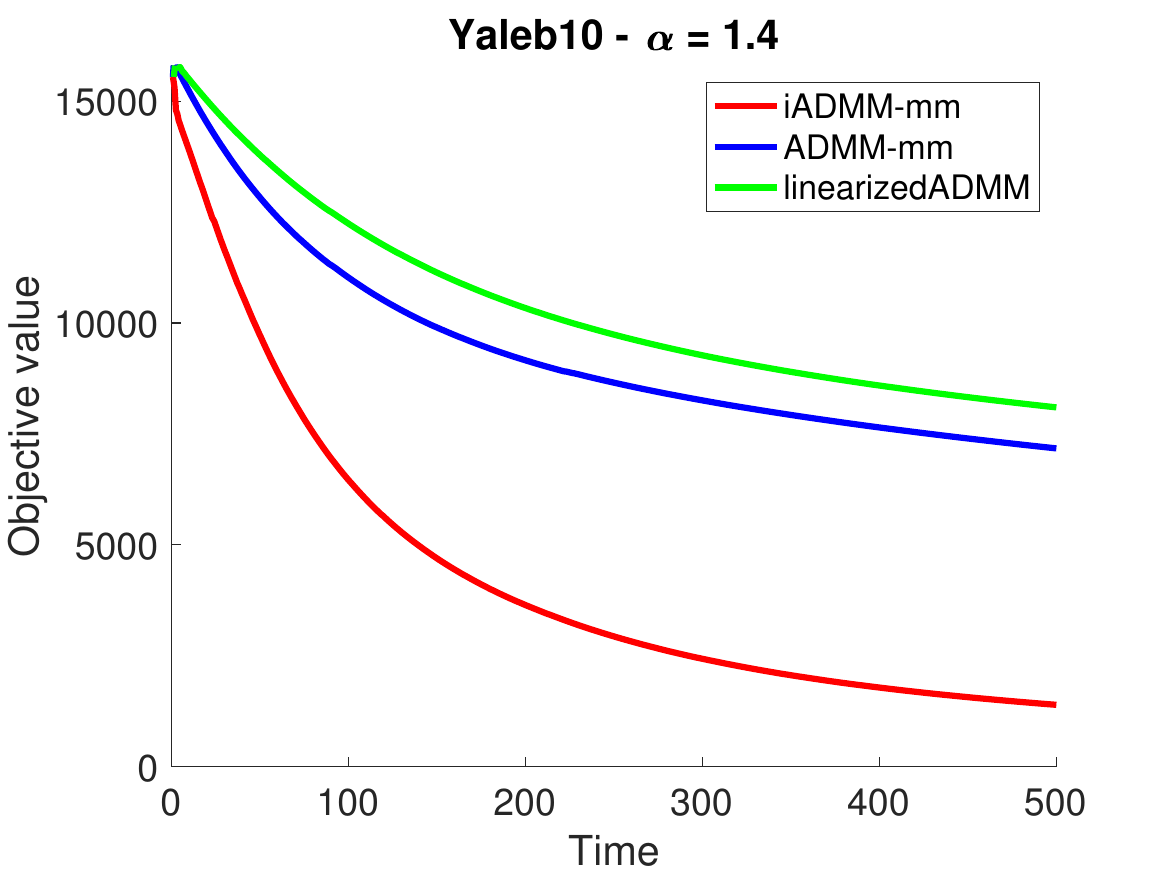} \\
\includegraphics[width=0.5\textwidth]{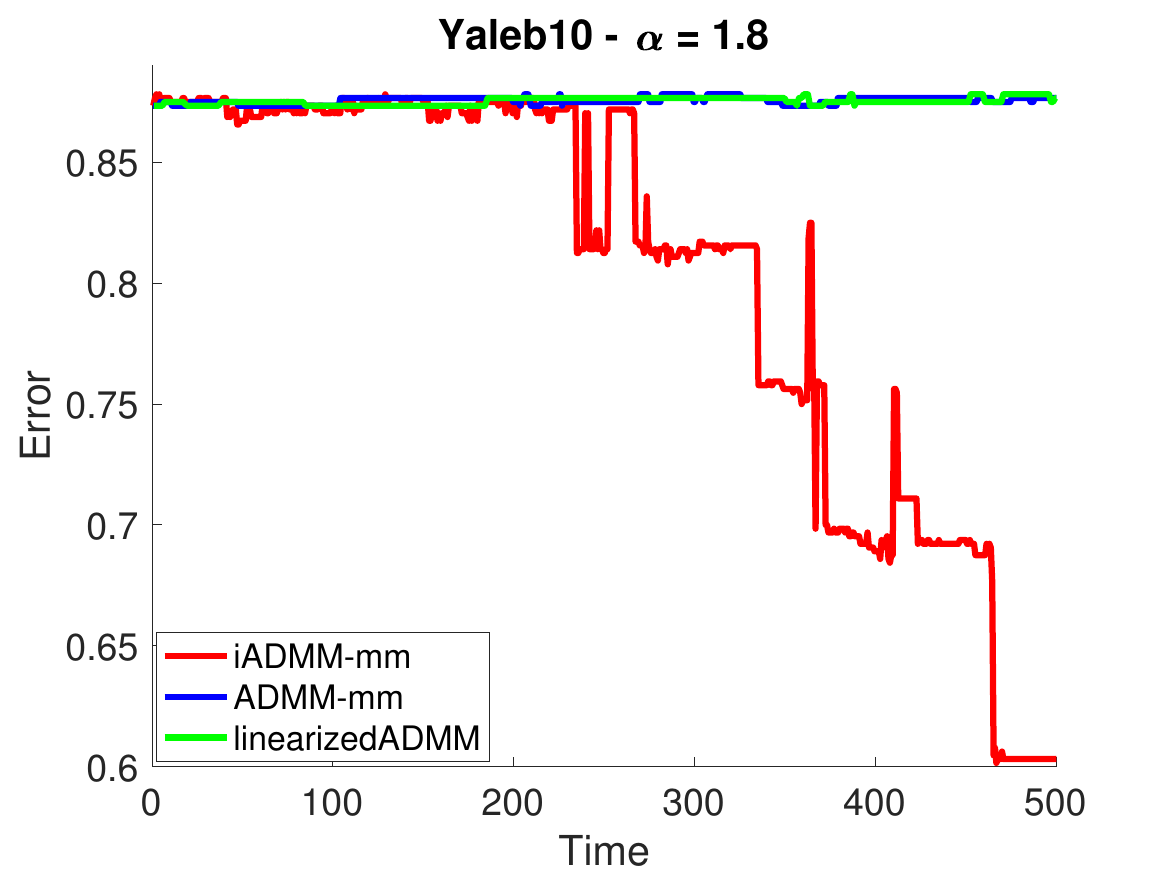}  & 
\includegraphics[width=0.5\textwidth]{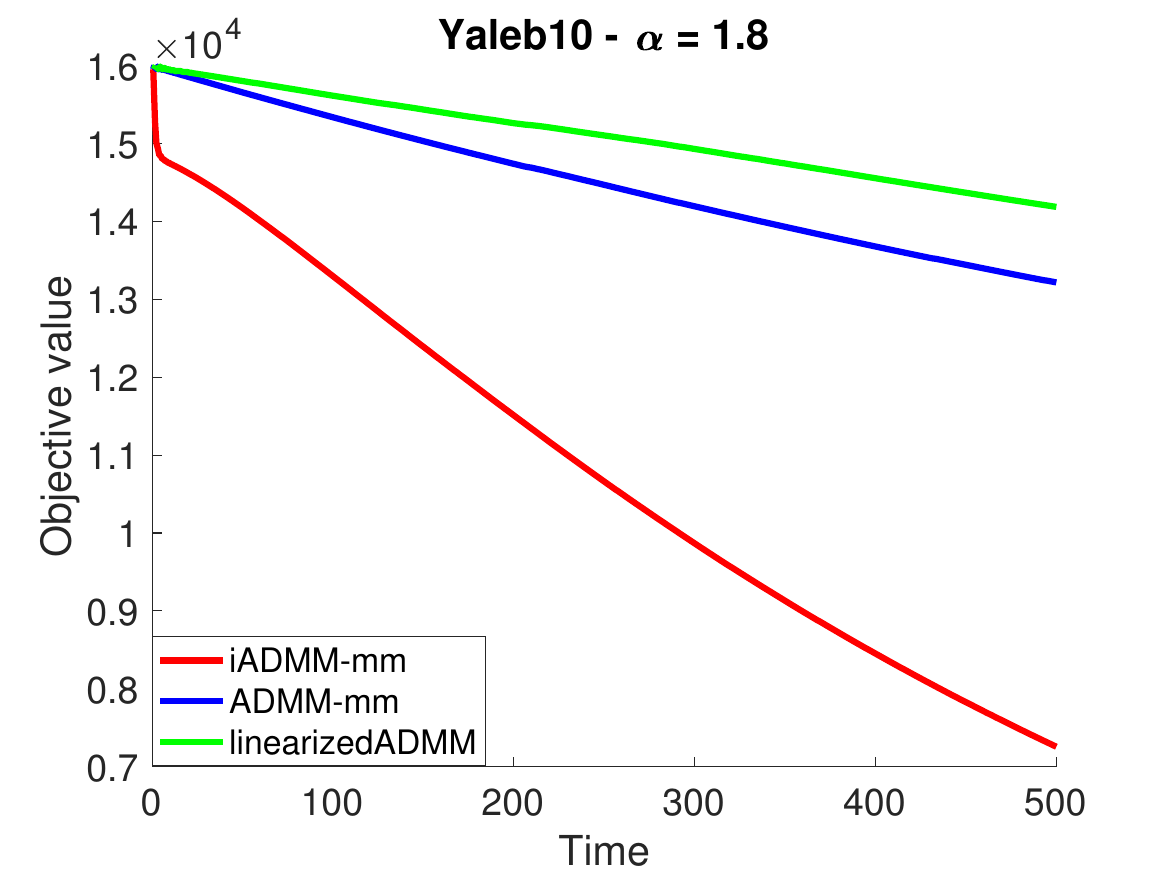}
\end{tabular}
\caption{Evolution of the segmentation error rate and the objective function value with respect to time on Yaleb10. 
\label{fig:Yaleb10}} 
\end{center}
\end{figure*}

In this experiment, we rerun the experiments from Section~\ref{sec:numerical} with other values for $\alpha$, namely $0.5, 1.4$ and $1.8$; see Figures~\ref{fig:Hopkins155}-\ref{fig:Yaleb10} (on pages~\pageref{fig:Hopkins155}-\pageref{fig:Yaleb10}). The penalty parameter $\beta$ is computed by $\beta = 2(2 + C_y)\alpha_2/C_y$, where $C_y = 1 - 10^{-6}$ and $ \alpha_2=\frac{3\alpha}{(1-|1-\alpha|)^2}$. 
Although the segmentation errors and objective function values differ for different values of $\alpha$, we observe  that, in all cases, iADMM-mm outperforms ADMM-mm which outperforms linearizedADMM. This confirms our observations from Section~\ref{sec:numerical}. \revise{On the other hand, we observe that the performances of ADMM-mm and linearizedADMM are similar for different values of $\alpha$; however, the performances of iADMM-mm (that is,  ADMM-mm with inertial terms) for $\alpha = 0.5$ and $\alpha = 1.4$ are slightly worse than for $\alpha=1$, and the value $\alpha = 1.8$ leads to significantly worse performances for iADMM-mm. 
It is known that, in the convex setting, the ADMM variants often perform better for $\alpha > 1$. However, in our experiments, $\alpha=1$ provides the best performance for iADMM-mm. A possible reason is that the global convergence of iADMM-mm has been established only for the case $\alpha=1$ (see Theorem~\ref{thrm:global}) while $\alpha\in (0,2)$ only guarantees a subsequential convergence (see Theorem~\ref{thrm:subsequential}).}

\section{Additional experiments for a regularized nonnegative matrix factorization problem} \label{appendix:NMF} 

{In the previous example, the function $f(X,Y) = \lambda_1 \|X\|_* + r_2(Y )$ was separable while our framework allows non-separable functions; see~\eqref{model} and the discussion that follows. To illustrate the use and effectiveness of iADMM on a non-separable case, let us} consider the following regularized nonnegative matrix factorization (NMF) problem
\begin{equation}
\label{NMFprob}
    \min_{W \in\mathbb R^{n\times r}_+ ,H \in\mathbb R^{r\times m}_+} \nicefrac{1}{2} \|X-WH\|^2 + c_1 \|W\|_F^2 + c_2 \|H\|_F^2, 
\end{equation}
where $X\in \mathbb R^{n\times m}$ is a given nonnegative matrix, and $c_1>0$ and $c_2>0$ are regularized parameters. Problem~\eqref{NMFprob} can be rewritten in the form of \eqref{model} as follows: 
\begin{equation}
    \label{NMF-rewritte}
    \begin{split}
           \min_{W \in\mathbb R^{n\times r}_+ ,H \in\mathbb R^{r\times m}_+} & \nicefrac{1}{2} \|X-W  H\|^2  + c_1 \|W\|_F^2 + c_2 \|Y\|_F^2, \\
         &  {\rm{such \,that}}\quad H -Y = 0.
    \end{split}
\end{equation}
In this case, $x_1=W$, $x_2=H$, $y=Y$, $f(W,H)=\frac{1}{2} \|X-W  H\|^2  +  \beta \|W\|_F^2 $, $g_1(W)$ and  $g_2(H)$ are indicator functions of $\mathbb R^{n\times r}_+$ and $\mathbb R^{r\times m}_+$ respectively, $h(Y)=c_2 \|Y\|_F^2$, $\mathcal A_1=0$, $\mathcal A_2=\mathcal I$,  $ \mathcal{B}= -\mathcal I$ (where $\mathcal I$  is identity operator),  and $b=0$. 
As $W\mapsto f(W,H)$ is $L_W$-Lipschitz smooth and $H\mapsto f(W,H)$ is $L_H$-Lipschitz smooth, where $L_W=\|H H^\top\|+2 c_1$ and $L_H=\|W^\top W\|$, we use the Lipschitz gradient surrogate for block $W$ and $H$ as in \eqref{Lipschitz_grad_surrogate}, and apply the inertial term as in the footnote~\ref{footnote_general} (that is, we apply inertial terms that also lead to the extrapolation for the block surrogate of $f$).
The augmented Lagrangian for \eqref{NMF-rewritte} is
$$\mL(W,H,Y,\omega)=f(W,H)+h(Y)+\langle H-Y,\omega
\rangle +\frac{\beta}{2}\|H-Y\|^2.$$
Applying iADMM for solving \eqref{NMF-rewritte}, the update of $W$ is 
\begin{equation}
    \begin{split}
\label{update_W}
    W^{k+1}&\in\arg\min_{W\in \mathbb R^{n\times r}_+} \langle -(X-\bar W^k H^k)(H^k)^\top+2 c_1 \bar W^k,W\rangle + \frac{L_W(H^k)}{2}\|W-\bar W^k\|^2  \\
    &=\max \Big\{\bar W^k -\frac{1}{L_W(H^k)}\big(-(X-\bar W^k H^k)(H^k)^\top+2 c_1 \bar W^k\big),0\Big\},
    \end{split}
\end{equation}
where $\bar W^k=W^k + \zeta_1^k (W^k-W^{k-1})$. Note that we have used extrapolation for the surrogate of $W\mapsto f(W,H)$. 
The update of $H$ is 
\begin{equation}
    \begin{split}
\label{update_H}
H^{k+1}&\in\arg\min_{H\in \mathbb R^{r\times m}_+} \langle -(W^{k+1})^\top(X- W^{k+1} \bar H^k)+\omega^k+\beta(\bar H^k - Y^k),H\rangle \\
&\qquad\qquad \qquad+ \frac{\beta+L_H(W^{k+1})}{2}\|H-\bar H^k\|^2  \\
&=\max \Big\{\bar H^k-\frac{1}{\beta+L_H(W^{k+1})}\big(-(W^{k+1})^\top(X- W^{k+1} \bar H^k)+\omega^k+\beta(\bar H^k - Y^k)\big) ,0\Big\},
    \end{split}
\end{equation}
where $\bar H^k=H^k + \zeta_2^k (H^k-H^{k-1})$. We do not use extrapolation for $Y$ (that is, $\delta_k=0$), and simply choose $\alpha=1$. The update of $Y$ is  
\begin{equation}
    \begin{split}
\label{update_Y}
    Y^{k+1} &\in \arg\min_Y  \langle -\omega^k+ 2 c_2 Y^k,Y\rangle  + \frac{\beta}{2} \|Y-H^{k+1}\|^2 + c_2 \|Y-Y^k\|^2 \\
    &= \frac{1}{\beta + 2 c_2}(\beta H^{k+1} + \omega^k ), 
    \end{split}
\end{equation}
while the update of $\omega$ is 
$$
\omega^{k+1}= \omega^k + \beta (H^{k+1}-Y^{k+1}). 
$$

\paragraph{Choosing parameters} By   Proposition~\ref{prop:Lsmooth-strong}, 
the update of $W$ in \eqref{update_W} implies that
 Inequality~\eqref{NSDP} is satisfied: 
 $$\mL(W^{k+1},H^k,Y^k,\omega^k)+\eta^k_1 \|W^{k+1}-W^k\|^2 \leq \mL(W^{k},H^k,Y^k,\omega^k)+\gamma^k_1 \|W^{k}-W^{k-1}\|^2,
$$
where $$\eta^k_1=\frac{L_W(H^k)}{2}, \quad  \gamma_1^k= \frac{L_W(H^k)}{2} (\zeta_1^k)^2.
$$
Note that we use $\eta^k_1$ instead of $\eta_1$ as this value varies along with the update of $H$ (because we used the extrapolation for the surrogate of $W\mapsto f(W,H)$).    
Similarly,  the update of $H$ in \eqref{update_H} implies that Inequality \eqref{NSDP} is satisfied: 
$$\mL(W^{k+1},H^{k+1},Y^k,\omega^k)+\eta^k_2 \|H^{k+1}-H^k\|^2 \leq \mL(W^{k+1},H^k,Y^k,\omega^k)+\gamma^k_2 \|H^{k}-H^{k-1}\|^2,
$$ 
where
$$\eta_2^k=\frac{L_H(W^{k+1})+\beta}{2}, \gamma_2^k=\frac{L_H(W^{k+1})+\beta}{2} (\zeta_2^k)^2.
$$
Because of the update of $Y$ in \eqref{update_Y}, the inequality in Proposition~\eqref{eq:NSDP-y} is satisfied: 
$$\mL(W^{k+1},H^{k+1},Y^{k+1},\omega^k)+\eta_y \|Y^{k+1}-Y^k\|^2 \leq \mL(W^{k+1},H^{k+1},Y^k,\omega^k)+\gamma^k_y \|Y^{k}-Y^{k-1}\|^2,
$$ 
where $\eta_y=c_2$ and $\gamma_y^k=0$. Following the same rationale that leads to  Theorem \ref{thrm:subsequential}, we obtain, as in \eqref{requirement}, 
$$\gamma_i^k \leq C_x \eta_i^{k-1}, \frac{2\alpha_2(2c_2)^2}{\beta} \leq C_y (\eta_y-\frac{\alpha_2 (2c_2)^2}{\beta}),
$$ 
where $\alpha_2=\frac{3\alpha}{\sigma_{\mathcal B}(1-|1-\alpha|)^2}=3$ and $0<C_x, C_y<1$. 
In our experiments, we choose 
$$\zeta_1^k=\min\Big\{\frac{a_{k-1}-1}{a_k},\sqrt{C_x\frac{L_W(H^{k-1})}{L_W(H^k)}} \Big\}, \zeta_2^k=\min\Big\{\frac{a_{k-1}-1}{a_k},\sqrt{C_x\frac{L_H(W^{k+1})+\beta}{L_H(W^k)+\beta}} \Big\},  
$$
where  $a_0=1$, 
$a_k=\frac12(1+\sqrt{1+4a_{k-1}^2})$, and $\beta \geq 4 c_2 \frac{(6+3C_y)}{C_y}$.

\paragraph{Experiments} We will compare iADMM with 
(i) ADMM (that is iADMM without using the inertial terms: $\zeta_1^k=\zeta_2^k=0$), 
and 
(ii) TITAN - the inertial block majorization minimization proposed in  \cite{Titan2020} that directly solves Problem~\eqref{NMFprob} and {competes favorably with the state of the art on the NMF problem  (see~\cite{Hien_ICML2020} which is a special case of TITAN)}. In our implementation, we use Lipschitz gradient surrogate for $W$ and $H$ and use default parameter setting for TITAN. 

 In the following experiments, we set the parameters $c_1$ and $c_2$ of Problem \eqref{NMFprob} to be $c_1=0.001$ and $c_2=0.01$. 

\begin{figure*}[h!]
\begin{center}
\begin{tabular}{cc}
\includegraphics[width=0.45\textwidth]{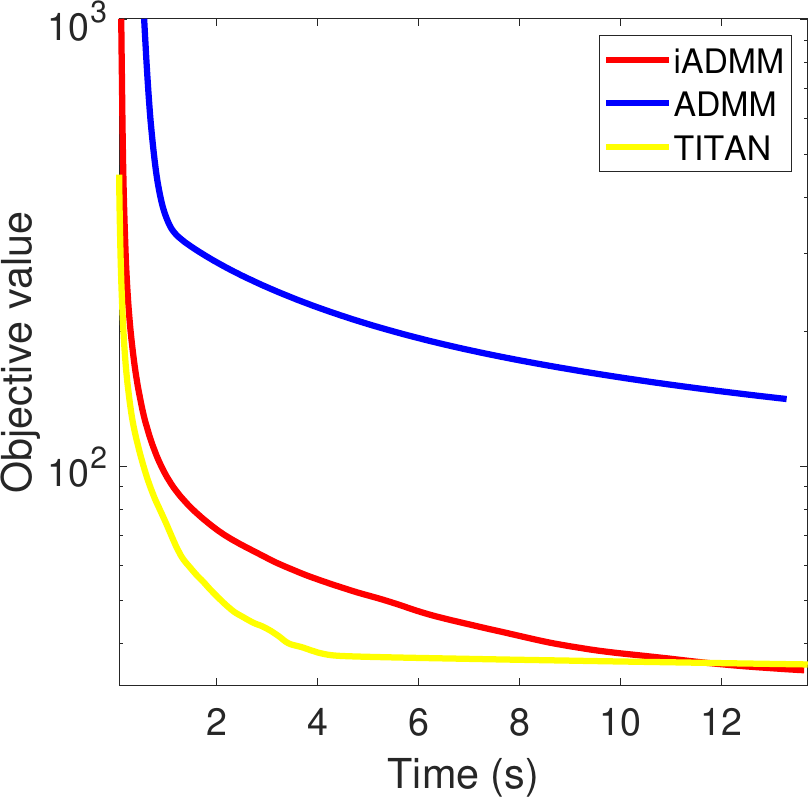}  & 
\includegraphics[width=0.45\textwidth]{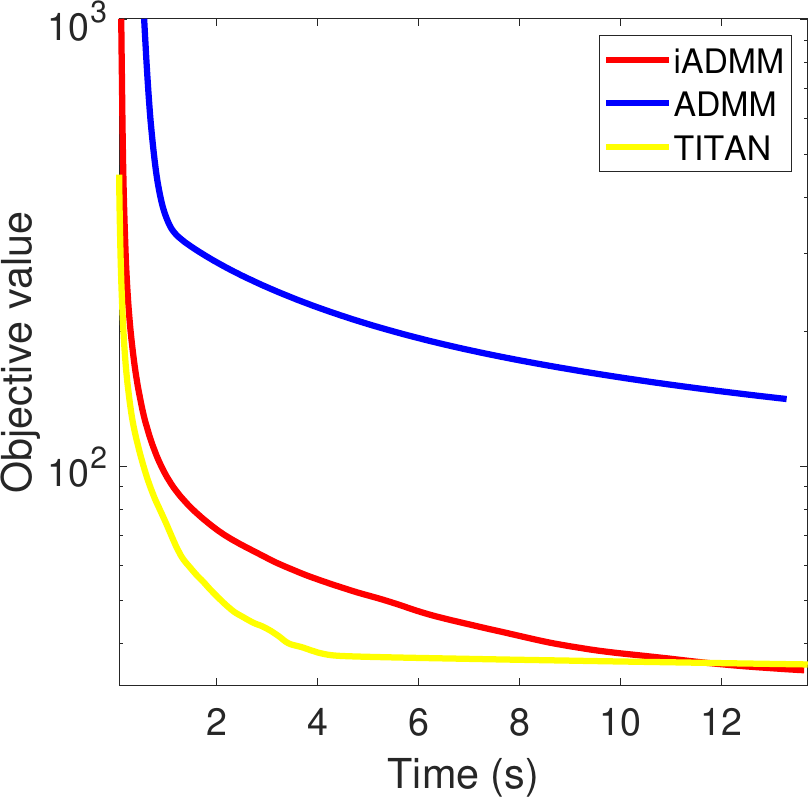} 
\end{tabular}
\caption{Evolution of the average value of the objective function value of Problem~\eqref{NMFprob} with respect to time on synthetic data sets with $ (n,m,r)= (500,200,20)$ (left) and $ (n,m,r)= (500,500,20)$ (right). 
\label{fig:nmf-synthetic}} 
\end{center}
\end{figure*}

In the first experiment, we generate 2 synthetic low-rank data sets $X$ with $(n,m,r)=(500,200,20)$ and $(n,m,r)=(500,500,20)$: we generate $U$ and $V$ by using the MATLAB command \texttt{rand(n,r)} and \texttt{rand(r,m)} respectively, and then let \texttt{X=U*V}. For each data set, we run each algorithm  with the same 30 random initial points \texttt{$W_0$=rand(n,r)},  \texttt{$H_0$=rand(r,m)} (for iADMM and ADMM we let \texttt{$Y_0$=$H_0$} and \texttt{$\omega_0$=zeros(r,m)}), and for each initial point we run each algorithm for 15 seconds. We report the evolution of the average objective function values of Problem~\eqref{NMFprob} with respect to time in 
Figure~\ref{fig:nmf-synthetic} and the mean $\pm$ std of the final objective function values in Table~\ref{tab:nmf-synthetic}.   We observe that iADMM outperforms ADMM which illustrates the acceleration effect. 
Among the algorithms, TITAN converges the fastest, but only slightly faster than iADMM.  However, iADMM provides the best final objective function values on average. 

\begin{figure*}[h!]
\begin{center}
\begin{tabular}{cc}
\includegraphics[width=0.47\textwidth]{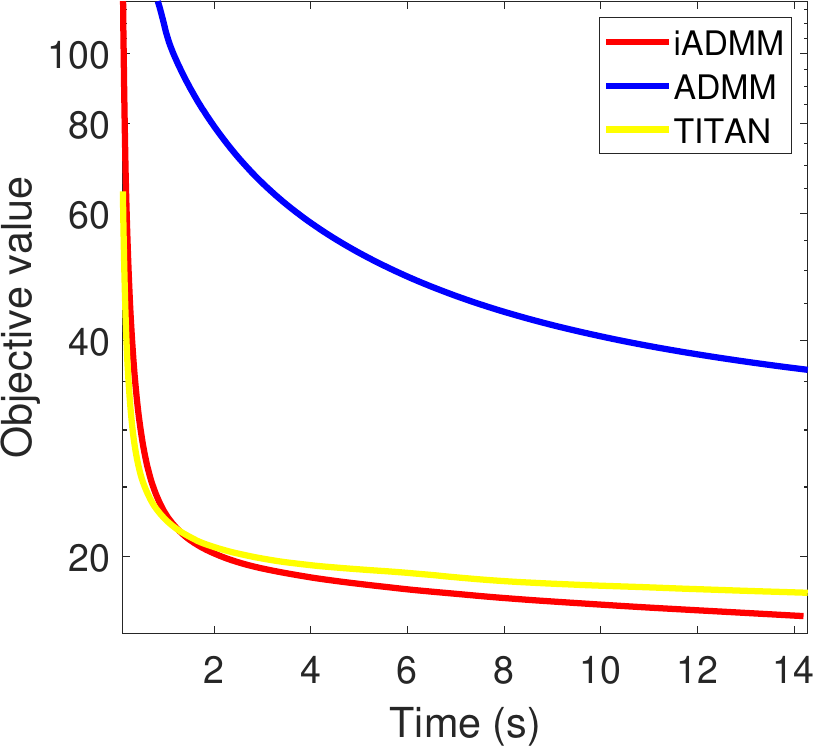}  & 
\includegraphics[width=0.46\textwidth]{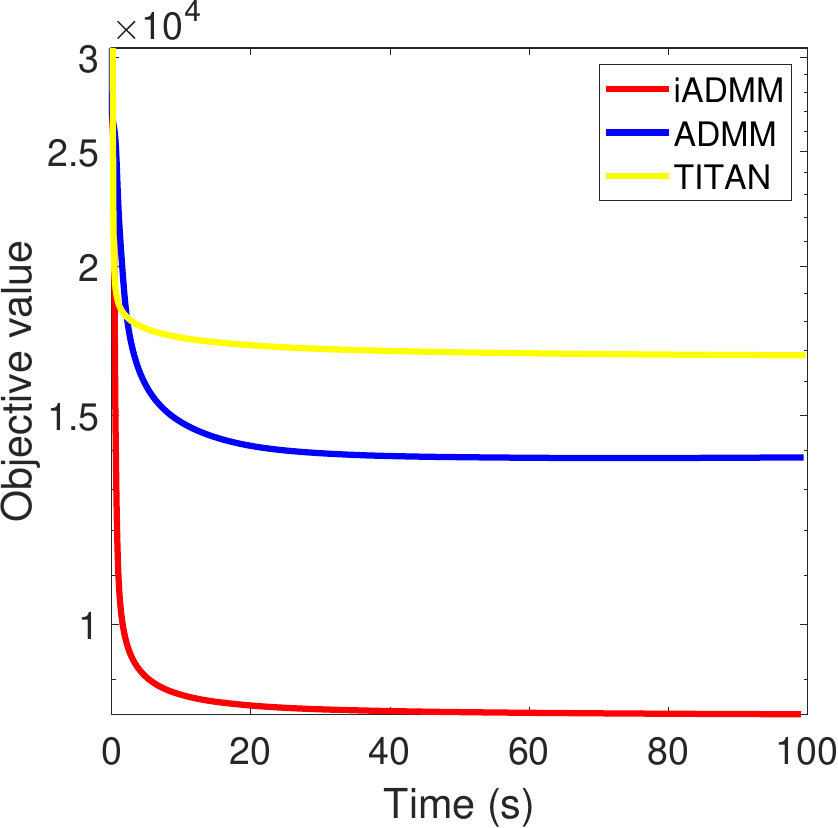} \\
\includegraphics[width=0.44\textwidth]{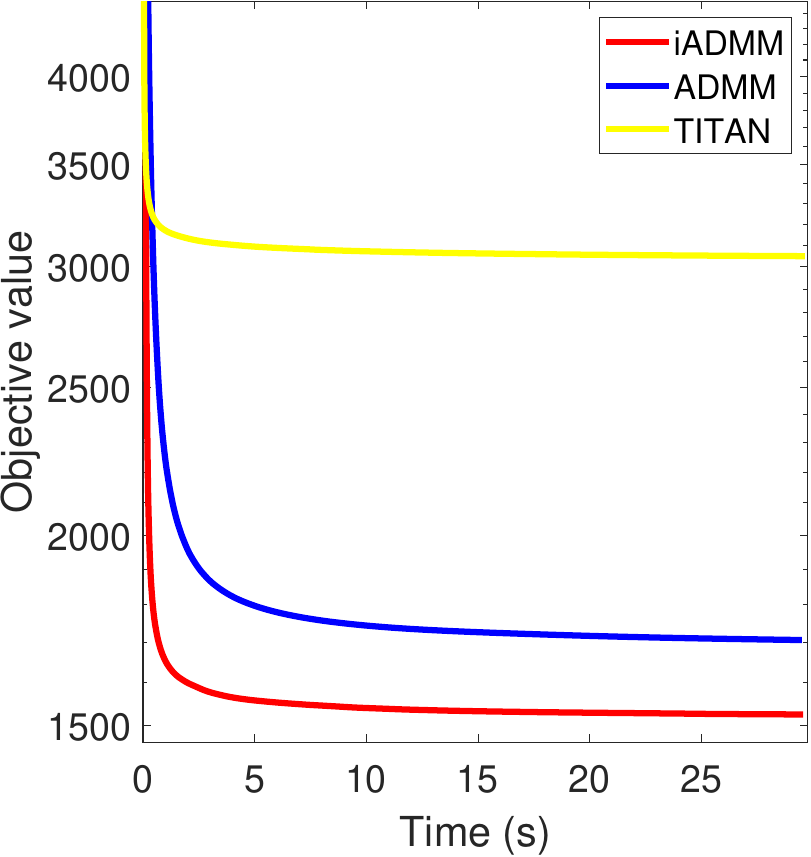}  & 
\includegraphics[width=0.44\textwidth]{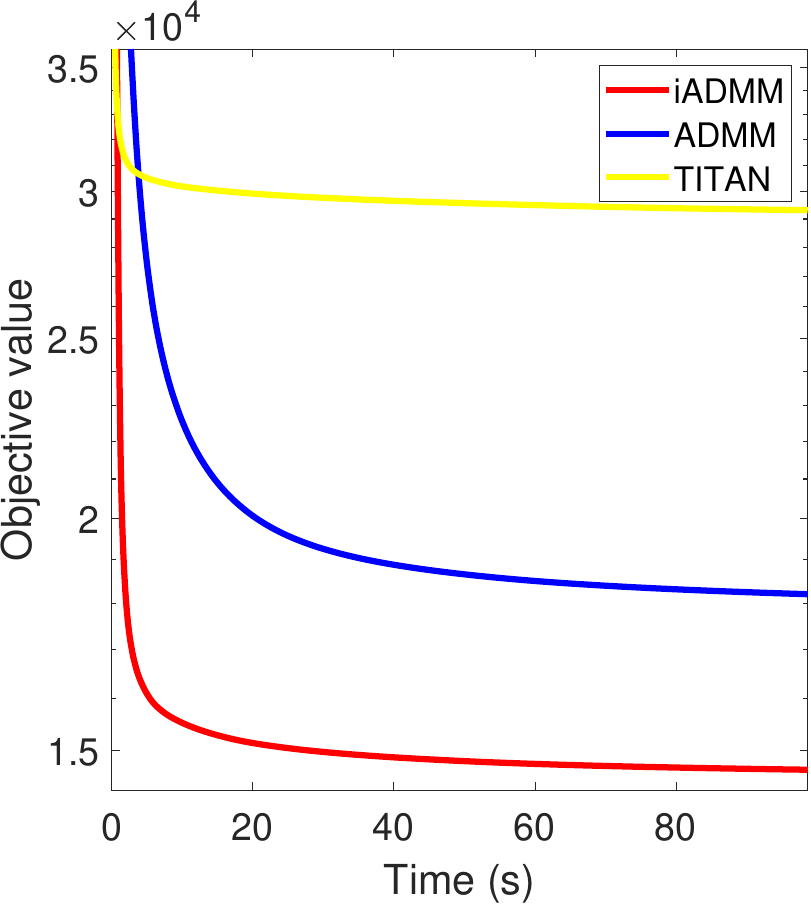}
\end{tabular}
\caption{Evolution of the average value of the objective function value of Problem~\eqref{NMFprob} with respect to time on the image data sets CBCL (top left), ORL (top right), Frey (bottom left) and Umist (bottom right). 
\label{fig:nmf-image}} 
\end{center}
\end{figure*}

In the second experiment, we test the algorithms on 4 image data sets  CBCL\footnote{\url{http://cbcl.mit.edu/software-datasets/heisele/facerecognition-database.html}} (2429 images of dimension $19 \times 19$), ORL\footnote{\url{https://cam-orl.co.uk/facedatabase.html}} (400 images of dimension $92 \times 112$), 
Frey\footnote{\url{https://cs.nyu.edu/~roweis/data.html}}(1965 images of dimension $28 \times 20$), and
 Umist\footnote{\url{https://cs.nyu.edu/~roweis/data.html}} (565 images of dimension $92 \times 112$). For each data set, we run each algorithm with the same 20 random initial points. We run each algorithm 100 seconds for the data sets Umist and ORL and 30 seconds for the data sets CBCL and Frey. We draw the evolution of the average objective functions values with respect to time in 
 Figure~\ref{fig:nmf-image} and the mean $\pm$ std of the final objective function values in Table~\ref{tab:nmf-image}. 
 
 Once again we observe that although iADMM converges slighly slower than TITAN, iADMM always produces the best objective function values among the three algorithms. On the other hand, ADMM also outperforms TITAN in term of the final objective function values. {This means that, for some reason, ADMM and iADMM are able to avoid spurious local minima more effectively than TITAN.} 

\begin{center}  
 \begin{table}[h!] 
 \caption{Mean and standard deviation of the objective function value over 30 random initializations on the synthetic data sets. 
 The best result is highlighted in bold.  \label{tab:nmf-synthetic}}
 \begin{center}  
 \begin{tabular}{|c|c|c|c|} 
 \hline $(n,m,r)$ & iADMM & ADMM & TITAN \\ 
 \hline 
$(500,200,20)$&     $\mathbf{16.443  \pm 5.015\times 10^{-1}}$ &   $35.873  \pm 2.299$   &   $17.751  \pm 1.092$ \\ 
$(500,500,20)$&  $\mathbf{34.289  \pm 4.492}$  & $135.037  \pm 6.409$ &   $35.799  \pm 1.525$ \\ 
  \hline 
  \end{tabular} 
  \end{center} 
 \end{table} 
 \end{center} 
 \begin{center}  
 \begin{table}[h!] 
 \caption{Mean and standard deviation of the objective function value over 20 random initializations on the image data sets. 
 The best result is highlighted in bold.  \label{tab:nmf-image}}
 \begin{center}  
 \begin{tabular}{|c|c|c|c|} 
 \hline Data set & iADMM & ADMM & TITAN \\ 
 \hline 
 CBCL & $\mathbf{1\,659.323  \pm 2.514}$ &  $1\,800.626  \pm 1.156\times 10^{1}$ &  $3\,321.104  \pm 7.271$\\
 ORL&  $\mathbf{8\,409.274  \pm 7.688}$  &  $13\,825.606  \pm 1.312\times 10^{2}$ &    $16\,844.426  \pm 1.439 \times 10^{1}$ \\ 
 Frey & $\mathbf{1\,525.242  \pm 3.555}$ & $1\, 706.385  \pm 7.380$ & $3\,048.246  \pm 4.737$ \\
Umist&     $\mathbf{14\,635.222  \pm 1.444 \times 10^{1}}$ &   $18\,195.557  \pm 1.059 \times 10^{2}$   &   $29\,316.019  \pm 3.433 \times 10^{1}$\\ 
  \hline 
  \end{tabular} 
  \end{center} 
 \end{table} 
 \end{center} 
 
\end{document}